\documentclass[11pt,twoside, leqno]{article}

\usepackage{extarrows}
\usepackage{graphics}
\usepackage{epsfig}
\usepackage{color}
\usepackage{verbatim}
\usepackage{hyperref}
\usepackage[displaymath,mathlines,pagewise]{lineno}
\usepackage{amssymb,amsmath,amsfonts,amsthm,color,mathrsfs}
\usepackage[Symbol]{upgreek}
\usepackage{txfonts}
\usepackage[nottoc,notlot,notlof]{tocbibind}
\usepackage[active]{srcltx}
\usepackage{picinpar,graphicx} 
\usepackage{bbm}
\usepackage{indentfirst}     

\usepackage{tikz}
\usetikzlibrary{positioning} 
\usepackage{xcolor} 

\usepackage{citeref} 

\usepackage{latexsym}
\usepackage{amsmath}
\usepackage{amssymb}
\usepackage{arydshln}

\allowdisplaybreaks
\def\rr{{\mathbb R}}

\def\zz{{\mathbb Z}}
\def\nn{{\mathbb N}}

\def\cc{{\mathbb C}}

\def\fz{\infty}

\def\supp{{\mathop\mathrm{\,supp\,}}}
\def\dist{{\mathop\mathrm {\,dist\,}}}

\def\hs{\hspace{0.15cm}}

\def\r{\right}
\def\lf{\left}

\def\eqref#1{(\ref{#1})}

\newtheorem{theorem}{Theorem}[section]

\newtheorem{lemma}[theorem]{Lemma}
\newtheorem{corollary}[theorem]{Corollary}
\newtheorem{proposition}[theorem]{Proposition}
\newtheorem{example}[theorem]{Example}
\theoremstyle{definition}
\newtheorem{remark}[theorem]{Remark}
\newtheorem{definition}[theorem]{Definition}
\newtheorem{assumption}[theorem]{Assumption}

\numberwithin{equation}{section}

\def\nn{{\mathbb N}}
\def\bx{{\mathbb X}}
\def\fz{\infty}
\def\r{\right}
\def\lf{\left}

\def\rr{{\mathbb R}}

\def\zz{{\mathbb Z}}

\begin{document}
\arraycolsep=1pt

\title{\Large\bf
 Local Hardy Spaces Associated with Ball Quasi-Banach Function Spaces and Non-negative Self-adjoint Operators on Spaces of Homogeneous Type and Their Applications
\footnotetext{\hspace{-0.35cm}
{\it 2020 Mathematics Subject Classification}.
{Primary 42B30; Secondary  42B25, 46E30.}
\endgraf{\it Key words and phrases.} Ball quasi-Banach function space, local Hardy space, atom, non-negative self-adjoint operator, local Riesz transforms.
}}

\author{ {Xiong Liu, Wenhua Wang, Tiantian Zhao and Hongliang Zuo}\\ }
\date{ }
\maketitle


\begin{center}
\begin{minipage}{13cm}\small
{
\noindent
{\bf Abstract.}
Let $(\mathbb X,d,\mu)$ be a space of homogeneous type in the sense of
Coifman--Weiss, let $X$ be a ball quasi-Banach function space on
$\mathbb X$ under suitable maximal-function and associate-space assumptions, and let $L$ be a non-negative self-adjoint operator on
$L^2(\mathbb X)$. Assume that, for every $t>0$, the semigroup
$e^{-tL}$ admits an integral kernel satisfying a Gaussian upper bound. In
this paper, we introduce and systematically
study the local Hardy space $h_L^X(\mathbb X)$ associated with $X$ and
$L$, defined in terms of a local Lusin area function together with an
appropriate low-frequency term. As applications of this theory, we establish the
boundedness of the local Riesz transform $\nabla(L+\kappa I)^{-1/2}$
from $h_L^X(\mathbb R^d)$ into the corresponding $X$-valued vector
function space for second-order divergence-form elliptic operators.
We also obtain a H\"{o}rmander-type spectral multiplier theorem for
$F(L+I)$ on $h_L^X(\mathbb X)$. Finally, the abstract results are
applied to local Orlicz-Hardy spaces, local variable Hardy spaces,
and local mixed-norm Hardy spaces. This theory develops
Goldberg's original local Hardy space theory [Duke Math. J. {\bf 46}
(1979), 27-42; MR0523600] to the setting of ball quasi-Banach
function spaces and non-negative self-adjoint operators on spaces of
homogeneous type. To the best of our knowledge,
several of the results obtained in this paper are new even in the
Euclidean setting $\mathbb X:=\mathbb R^d$.
}
\end{minipage}
\end{center}

\tableofcontents


\section{Introduction}\label{sec1}
The real-variable theory of Hardy spaces and BMO spaces on Euclidean spaces
has played a fundamental role in harmonic analysis and partial differential
equations; see, for instance, \cite{cw77,fs72,g79,g14,s60}. This theory goes
back, in particular, to the work of Stein--Weiss \cite{s60} on Hardy spaces
and John--Nirenberg \cite{jn61} on BMO spaces. Hardy spaces and BMO spaces
serve as effective endpoint substitutes for Lebesgue spaces and are widely
used, through interpolation and duality, in the study of \(L^p\)-boundedness
\((1<p<\infty)\) of singular integral operators.

Motivated by various applications, many function spaces more general than
Lebesgue spaces have been introduced and systematically studied, including
Orlicz spaces, Morrey spaces, weighted Lebesgue spaces, mixed-norm Lebesgue
spaces, and variable Lebesgue spaces. Among these, quasi-Banach function
spaces provide a useful abstract framework which covers,
for example, variable
Lebesgue spaces and Orlicz spaces. However, several important spaces, such as
weighted Lebesgue spaces, Morrey spaces, and Herz spaces, are not necessarily
quasi-Banach function spaces in the usual sense; see, for instance,
\cite{shyy17}. Recently,
 Sawano et al.
\cite{shyy17} introduced ball quasi-Banach function spaces on \(\rr^d\) and
developed the real-variable theory of Hardy-type spaces \(H^X(\rr^d)\)
associated with such spaces. Compared with quasi-Banach function spaces, ball
quasi-Banach function spaces are more flexible and include a wider class of
function spaces. For further developments on ball quasi-Banach function
spaces and related Hardy-type spaces, we refer the reader to
\cite{syy2022,syy22,wyy20,wyyz21,yhyy23}.

On the other hand, many problems in harmonic analysis naturally go beyond the
Euclidean setting. In the early 1970s, Coifman and Weiss introduced spaces of
homogeneous type \(\bx\) \cite{cw77}, which provide a general metric-measure
framework including, among others, Lie groups with polynomial growth,
Riemannian manifolds with non-negative Ricci curvature, and Euclidean spaces
equipped with non-standard measures such as Bessel measures. During the past
several decades, the theory of singular integrals and function spaces on
spaces of homogeneous type has undergone substantial development; see, for
example, \cite{dy03,dy05,gly08,HS1994,HHLLYY} and the references therein.

Another important line of research concerns Hardy spaces associated with
operators. Let \(L\) be a non-negative self-adjoint operator on \(L^2(\bx)\).
The operator \(L\) generates the semigroup \(\{e^{-tL}\}_{t>0}\), and Hardy
spaces associated with \(L\) are usually defined under suitable off-diagonal
or heat kernel estimates for this semigroup. Auscher et al. \cite{adm05} and
Duong--Yan \cite{dy03,dy05} studied Hardy spaces associated with operators
whose heat kernels satisfy pointwise Gaussian upper bounds. Since such bounds
may fail for some important operators, for instance second-order divergence-form elliptic operators with complex coefficients, Hofmann et al.
\cite{hlmmy07} developed Hardy spaces under Davies--Gaffney type estimates
instead of pointwise kernel bounds. More recently, in order to place various
operator-adapted Hardy spaces into a common framework, Lin et al.
\cite{lyyy23} studied Hardy-type spaces \(H_L^X(\bx)\) associated
simultaneously with a non-negative self-adjoint operator \(L\) and a ball
quasi-Banach function space \(X(\bx)\). For further results on Hardy-type
spaces associated with operators, we refer the reader to
\cite{abdr20,bl11,dy05b,hlmmy07,jy10,sy10} and the references therein.

The distinction between global and local Hardy spaces goes back to
Goldberg \cite{g79}. One of the characteristic features of the classical
local Hardy space is that cancellation is required only for atoms supported
on small balls, whereas atoms supported on balls of unit or larger radius
need not satisfy any moment condition. This inhomogeneous structure makes
local Hardy spaces particularly suitable for problems governed primarily by
small-scale behavior and for operators containing a non-trivial
low-frequency component.

Recently,
local Hardy spaces associated with operators have been studied
in several particular settings. For example, Jiang, Yang and Zhou \cite{jyz09}
investigated localized Hardy spaces associated with operators, while Cao,
Mayboroda and Yang \cite{cmy17} considered inhomogeneous higher-order
elliptic operators. Almeida, Betancor, Dalmasso and Rodr\'{\i}guez-Mesa
\cite{abdr20} developed a local theory in the variable-exponent setting,
and more recent weighted results can be found in \cite{csyyy25}.
These works provide important operator-adapted local theories, but they are
mainly formulated for particular underlying function spaces or particular
classes of operators. A unified local theory simultaneously adapted to a
general ball quasi-Banach function space, a non-negative self-adjoint
operator, and a doubling metric measure space has not yet been systematically
developed.

The aim of the present paper is to establish such a theory. Let
\((\mathbb X,d,\mu)\) be a doubling metric measure space, let \(X(\mathbb X)\)
be a ball quasi-Banach function space satisfying suitable vector-valued
maximal and associate-space assumptions, and let \(L\) be a non-negative
self-adjoint operator on \(L^2(\mathbb X)\) whose heat kernel satisfies a
Gaussian upper bound. For \(f\in L^2(\mathbb X)\), we define the local
Hardy-type quasi-norm by
$$
 \|f\|_{h_L^X(\mathbb X)}
 :=
 \lf\|S_L^{\rm loc}(f)\r\|_{X(\mathbb X)}
 +
 \lf\|S_I(e^{-L}f)\r\|_{X(\mathbb X)}.
$$
Here \(S_L^{\rm loc}\) is obtained by restricting the Lusin area function
associated with \(L\) to \(0<t<1\), and \(S_I\) is the global area function
associated with the identity operator. The first term measures the
small-scale component of \(f\), while the second term recovers the
low-frequency information which is not detected by the truncated square
function. In this sense, the resulting space may be viewed as an
operator-adapted analogue of Goldberg's local Hardy space in the setting of
ball quasi-Banach function spaces.

Several difficulties arise in developing this theory. Firstly, a general
doubling metric measure space has neither a Fourier transform nor a
translation or dilation structure. Second, a ball quasi-Banach function
space need not be locally convex, so that the usual duality and interpolation
arguments are not directly available. Third, the small and large-scale
components obey different cancellation principles: molecules associated with small balls require
\(L\)-adapted cancellation, whereas no cancellation should be imposed at
unit or larger scales. Finally, the truncated square function alone does not
control the low-frequency part of a function.

To overcome these difficulties, we use atomic decompositions of global
and local tent spaces, Gaussian off-diagonal estimates, vector-valued
Hardy--Littlewood maximal inequalities, and the spectral functional calculus
of \(L\). We introduce the auxiliary Hardy space \(H_I^X(\mathbb X)\)
associated with the identity operator and establish its atomic and molecular
characterizations. A Calder\'on-type reproducing formula then separates the
small-scale part from the low-frequency part and leads naturally to two
classes of local molecules: type \({\rm(II)}\) molecules on balls of radius
less than one, which satisfy \(L\)-adapted cancellation conditions, and type
\({\rm(I)}\) molecules on balls of radius at least one, for which no
cancellation is required.

Our first main result establishes a molecular characterization of the local
Hardy space \(h_L^X(\mathbb X)\).  Thus every element of \(h_L^X(\mathbb X)\)
admits a representation in terms of type \({\rm(I)}\) and type
\({\rm(II)}\) local molecules satisfying the prescribed coefficient
condition. Conversely, every such molecular representation determines an
element of \(h_L^X(\mathbb X)\). Here type \({\rm(II)}\) molecules, which
are associated with balls of radius less than one, satisfy
\(L\)-adapted cancellation conditions, whereas type \({\rm(I)}\)
molecules, associated with balls of radius at least one, require no
cancellation.

Our second main result relates the local Hardy space to global Hardy spaces
associated with \(L\) and its shifted operator. Without imposing any
spectral-gap condition on \(L\), we prove that
\[
h_L^X(\mathbb X)=H_{L+I}^X(\mathbb X)
\]
with equivalent quasi-norms. The shift by the identity is natural in the
local setting because
\[
e^{-t(L+I)}=e^{-t}e^{-tL},\qquad t>0,
\]
and the additional exponential decay controls the large-time, or
low-frequency, component of the theory. If, in addition,
$
\inf\sigma(L)>0,
$
then \(L\) itself has a spectral gap, and we further obtain
\[
h_L^X(\mathbb X)=H_L^X(\mathbb X)
\]
with equivalent quasi-norms. Consequently, \(H_{L+I}^X(\mathbb X)\)
provides a global realization of \(h_L^X(\mathbb X)\) in general, while
\(H_L^X(\mathbb X)\) plays the same role whenever \(L\) has a positive
spectral gap.

As applications, we consider the real symmetric second-order
divergence-form elliptic operator
$
L=-\operatorname{div}(A\nabla)
$
on \(\mathbb R^d\), where \(A\) is bounded and uniformly elliptic. For
every fixed \(\kappa>0\), we prove that the local Riesz transform
\[
\mathcal R_{L,\kappa}^{\mathrm{loc}}
:=
\nabla(L+\kappa I)^{-1/2}
\]
extends to a bounded operator
\[
\mathcal R_{L,\kappa}^{\mathrm{loc}}:
h_L^X(\mathbb R^d)
\longrightarrow
X(\mathbb R^d;\mathbb C^d).
\]
We also use the identification
\(h_L^X(\mathbb X)=H_{L+I}^X(\mathbb X)\) to transfer a global
H\"ormander-type spectral multiplier theorem for \(L+I\) to the local
Hardy space \(h_L^X(\mathbb X)\). In particular, suitable bounded Borel
multipliers \(F\) give rise to bounded operators \(F(L+I)\) on
\(h_L^X(\mathbb X)\). Finally, the abstract results can be applied to some special example, for example, local
Orlicz-Hardy spaces, local variable Hardy spaces, and local mixed-norm
Hardy spaces.

The remainder of this paper is organized as follows.

In Section \ref{s2}, we recall some necessary preliminaries, including the
basic notions and properties of spaces of homogeneous type, ball
quasi-Banach function spaces, global tent spaces, and local tent spaces.

In Section \ref{s3}, we introduce the local Hardy space \(h_L^X(\bx)\)
associated with the ball quasi-Banach function space \(X(\bx)\) and the
non-negative self-adjoint operator \(L\). We also state the main results of
this paper, including the molecular characterization of \(h_L^X(\bx)\), the
coincidence between \(h_L^X(\bx)\) and \(H_L^X(\bx)\) under the spectral gap
condition \(\inf \sigma(L)>0\), and the identification of \(h_L^X(\bx)\) with
\(H_{L+I}^X(\bx)\).

Section \ref{s4} is devoted to the auxiliary Hardy space \(H_I^X(\bx)\) associated
with the identity operator \(I\). In this section, we establish the atomic and
molecular characterizations of \(H_I^X(\bx)\), which will be used to handle
the low-frequency part appearing in the definition of \(h_L^X(\bx)\).

In Section \ref{s5}, we prove the main results stated in Section \ref{s3}. More
precisely, we first prove the molecular characterization of \(h_L^X(\bx)\),
and then establish the equivalence between the local Hardy space
\(h_L^X(\bx)\) and the global Hardy spaces \(H_L^X(\bx)\) and
\(H_{L+I}^X(\bx)\) under the corresponding assumptions.

In Section \ref{s6x}, we present two applications of the local Hardy
spaces developed in this article.
First, we study local Riesz transforms in the Euclidean setting. Let
$
L=-\operatorname{div}(A\nabla)
$
on \(\mathbb R^d\), where \(A\) is real symmetric, bounded, and uniformly
elliptic. For every fixed \(\kappa>0\), we prove that the local Riesz
transform
\[
\mathcal R_{L,\kappa}^{\rm loc}
:=
\nabla(L+\kappa I)^{-1/2}
\]
extends uniquely to a bounded operator
\[
\mathcal R_{L,\kappa}^{\rm loc}:
h_L^X(\mathbb R^d)
\longrightarrow
X(\mathbb R^d;\mathbb C^d).
\]
The proof relies on the molecular characterization of
\(h_L^X(\mathbb R^d)\), the \(L^2\)-boundedness of
\(\nabla(L+\kappa I)^{-1/2}\), and suitable off-diagonal estimates for the
gradient of the heat semigroup generated by \(L+\kappa I\). These estimates
yield the required annular decay for the images of both type
\({\rm(I)}\) and type \({\rm(II)}\) local molecules, from which the desired
boundedness follows.
Second, we establish a H\"{o}rmander-type spectral multiplier theorem on
\(h_{L}^{X}(\rr^d)\), showing that spectral multipliers satisfying
an appropriate localized smoothness condition induce bounded operators
on \(h_{L}^{X}(\mathbb{R}^d)\).

Finally, in Section \ref{s6}, we apply the abstract theory to several concrete
classes of Hardy-type function spaces. These include local Orlicz-Hardy spaces in
Subsection \ref{OHS}, local variable Hardy spaces in Subsection \ref{VHS}, and
local mixed-norm Hardy spaces in Subsection \ref{MHS}.

{\bf Notation:}
Throughout the paper, we set
$
 \nn:=\{1,2,\ldots\},
 \zz_+:=\{0\}\cup\nn,
 \rr_+:=(0,\infty).
$
We denote by \(C\) a positive constant independent of the main parameters; its
value may change from line to line. We write \(A\lesssim B\) if
\(A\leq CB\) for some constant \(C>0\), and write \(A\sim B\) if
\(A\lesssim B\) and \(B\lesssim A\). For any \(\lambda>0\) and any ball
\(B:=B_d(x_B,r_B)\subset\bx\), we denote by \(\lambda B\) the ball centered
at \(x_B\) with radius \(\lambda r_B\).
For a locally integrable function \(f\) and a ball \(B\subset\bx\), we write
\[
 f_B:=\fint_B f(x)\,d\mu(x)
 :=\frac1{\mu(B)}\int_B f(x)\,d\mu(x).
\]
For \(1\leq p\leq\infty\), \(p'\) denotes the conjugate exponent of \(p\),
that is, \(1/p+1/p'=1\), with the usual convention when \(p=1\) or
\(p=\infty\). For a measurable set \(E\subset\bx\), we let
\(E^\complement:=\bx\setminus E\) and denote by \({\bf 1}_E\) its
characteristic function.

\section{Preliminaries}\label{s2}

In this section, we will recall some elementary definitions and properties of  spaces of homogeneous type, ball quasi-Banach function spaces and tent spaces.
\subsection{Spaces of homogeneous type}\label{homogene}
A function $d:\bx\times\bx\rightarrow[0,\,\fz)$ is said to be a quasi-metric function if it satisfies
\begin{enumerate}
  \item  $d(x,\,y)=d(y,\,x)\geq 0$, for any $x,\,y\in\bx$;
  \item $d(x,\,y)=0$ if and only if $x=y$;
  \item There exists a constant $A_d\in[1,\infty)$ such that, for any $x,\,y,\,z\in\bx$,
 \begin{align}\label{noimp}d(x,\,y)\leq A_d[d(x,\,z)+d(z,\,y)].\end{align}
\end{enumerate}
Moreover, $\mu$ is said to be a doubling measure if there is a constant $C>0$ such that, for any $x\in\bx$ and $r>0,$
\begin{align}\label{e2.1}
\mu(B(x,\,2r))\leq C\mu(B(x,\,r)).
\end{align}

In what follows, we say that $(\bx,\, d,\, \mu)$ is a {\it space of
homogeneous type} in the sense of Coifman and Weiss \cite{cw77} if
$\bx$ is a metric space with a quasi-metric function $d$ and a non-negative, Borel, doubling measure $\mu$. Note that the constant $A_d$ on the right-hand side of inequality \eqref{noimp} plays no essential role in our proof below, so we shall only present our proof for the case of $A_d=1$.

In what follows, for any ball $B(x,\,r):=\{y\in\bx:d(x,\,y)<r\}\subset\bx$, we define the volume functions
\begin{align*}
V_r(x):=\mu(B(x,\,r)) \ \ \mathrm{and} \ \ V(x,\,y):=\mu(B(x,\,d(x,\,y))).
\end{align*}
From \eqref{e2.1}, we have the following strong homogeneity property:
\begin{align}\label{e2.2}
V_{\lambda r}(x)\leq C\lambda^dV_{r}(x)
\end{align}
for some constants $C,d>0$  independent of all $\lambda\geq 1$ and $x\in\bx$. The smallest constant $d$ satisfying the above inequality is said to be the homogeneous dimension of $\bx$.

The following lemma states some basic  properties of spaces of
homogeneous type, it is essential to obtain the main results.
\begin{lemma}[\rm\cite{ghl09}]\label{l2.0}
Let $(\bx,\, d,\, \mu)$ be a space of
homogeneous type. Then we have
\begin{enumerate}
\item[\rm{(i)}] For any $x\in\bx$ and $r>0$, one has
$V(x,\,y)\thicksim V(y,\,x)$ and
$$V_r(x)+V_r(y)+V(x,\,y)\thicksim V_r(y)+V(x,\,y)\thicksim V_r(x)+V(x,\,y)
\thicksim \mu(B(x,\,r+d(x,\,y))).$$
\item[\rm{(ii)}]  There exist two constants $C>0$ and $0\leq\gamma\leq d$ such that
$$V_{r_1}(x)\lesssim \lf[1+\frac{r_1+d(x,\,y)}{r_2}\r]^{\gamma}V_{r_2}(y).$$
Here all the implicit constants are independent of $x,\,y\in\bx$ and $r>0$.
\end{enumerate}
\end{lemma}

\subsection{Ball quasi-Banach function spaces}\label{nclp}
In this subsection, we recall the definition of ball quasi-Banach function spaces, which was originally studied in \cite[Definition 2.3]{shyy17}. In what follows,
we always use the symbol $\mathfrak{M}(\mathbb{X})$ to
denote the set of all measurable functions on $\bx$. For any
$x\in\bx$ and $r\in(0,\,\infty)$, let $B(x,\,r):=\{y\in \bx :d(x,y)<r\}$ and
\begin{align*}
\mathfrak{{B}}:=\{B(x,\,r): x\in\bx \ \ {\text{and}}\  \ r\in(0,\,\infty)\}.
\end{align*}

 \begin{definition}\label{d2.1xx}
 A quasi-Banach space $X(\bx)\subset\mathfrak{M}(\bx)$ is called a {\it ball quasi-Banach function space} (for short, BQBF space) on $(\bx,\,d,\,\mu)$ if it
satisfies
\begin{enumerate}
\item[\rm{(i)}]
$\|f\|_{X(\bx)}=0$ implies that $f=0$ almost everywhere;
\item[\rm{(ii)}] $|g|\leq|f|$ almost everywhere implies that $\|g\|_{X(\bx)}\leq\|f\|_{X(\bx)}$;
\item[\rm{(iii)}] $0\leq f_m\uparrow f$ as $m\rightarrow\infty$ almost everywhere implies that $\|f_m\|_{X(\bx)}\uparrow\|f\|_{X(\bx)}$ as $m\rightarrow\infty$;
\item[\rm{(iv)}] $B\in\mathfrak{B}$ implies that ${\bf 1}_B\in {X(\bx)}$.
\end{enumerate}
 \end{definition}

Furthermore, a ball quasi-Banach function space ${X(\bx)}$ is called a ball Banach function
space, for short, BBF space, if the norm of ${X(\bx)}$ satisfies that for any $f,\,g\in {X(\bx)}$,
\begin{align*}
\|f+g\|_{{X(\bx)}}\leq\|f\|_{{X(\bx)}}+\|g\|_{{X(\bx)}}
\end{align*}
and, for any $B\in\mathfrak{B}$, there exists a positive constant $C_{(B)}$, depending on $B$, such that,
for any $f\in {X(\bx)}$,
\begin{align*}
\int_{B}|f(x)|d\mu(x) \leq C_{(B)}\|f\|_{X(\bx)}.
\end{align*}

Notice that, in Definition \ref{d2.1xx}, if we replace any ball $B$ by any bounded measurable set $E$,
we  obtain its another equivalent definition.

\begin{definition}
 For any given ball Banach function space ${X(\bx)}$, its {\it associate space (also called the K\"{o}the
dual space)} $X'(\bx)$ is defined by setting
$$X':=\lf\{f\in\mathfrak{M}(\bx):\|f\|_{X'(\bx)}<\infty\r\},$$
where, for any $f\in\mathfrak{M}(\bx)$,
$\|f\|_{X'(\bx)}:=\sup\lf\{\|fg\|_{L^1(\bx)}:g\in X(\bx), \|g\|_{X(\bx)}=1\r\}$
and $\|\cdot\|_{X'(\bx)}$ is called the associate norm of $\|\cdot\|_{X(\bx)}$.
 \end{definition}


\begin{lemma}\rm{\cite{yhyy23}}
 Every ball Banach function space $X(\bx)$ coincides with its second associate
space $X''(\bx)$. In other words, a function $f$ belongs to $X(\bx)$ if and only if it belongs to $X''(\bx)$
and, in that case,
$$\|f\|_{X(\bx)}=\|f\|_{X''(\bx)}.$$
\end{lemma}

We also need to recall the definition of the $p$-convexification and the concavity of ball quasi-Banach function spaces, which is a
part of \cite[Definition 2.6]{shyy17}.
\begin{definition} \label{d2.5x}
Let $X$ be a ball quasi-Banach function space and $p\in(0,\,\infty)$.
\begin{enumerate}
\item[\rm{(i)}] The {\it $p$-convexification} $X^p(\bx)$ of $X(\bx)$ is defined by setting
$$X^p(\bx):=\lf\{f\in\mathfrak{M}(\bx):|f|^p\in X(\bx)\r\},$$
equipped with the quasi-norm $\|f\|_{X^p(\bx)}:=\lf\||f|^p\r\|_{X(\bx)}^{\frac1p}$.
\item[\rm{(ii)}] The space $X$ is said to be {\it$p$-concave} if there exists a positive constant $C$ such that, for any $\{f_k\}_{k\in\nn}\subset X^{\frac1p}(\bx)$,
\begin{align*}
\sum_{k\in\nn}\|f_k\|_{X^{\frac1p}(\bx)}\leq C\lf\|\sum_{k\in\nn}\lf|f_k\r|\r\|_{X^{\frac1p}(\bx)}.
\end{align*}
In particular, when $C=1$, $X$ is said to be strictly $p$-concave.
\end{enumerate}
 \end{definition}

 For any $f\in \mathfrak{M}(\bx)$ and $x\in\bx$, the {\it Hardy-Littlewood maximal function} $M_{\mathrm{HL}}(f)$ is defined by
\begin{align*}
M_{\mathrm{HL}}(f)(x):=\sup_{r>0}\frac{1}{\mu(B(x,\,r))}\int_{B(x,\,r)}|f(y)|\,d\mu(y),
\end{align*}
where $B(x,\,r)\in\mathfrak{{B}}$.

In this paper, we need some essential assumptions on the ball quasi-Banach function space $X(\bx)$ as follows.
\begin{assumption}\label{a11}
 Let $X(\bx)$ be a BQBF space. Assume that there exists a positive constant $p$ such
that, for any given $t\in(0,\,p)$ and $u\in(1,\,\fz)$, there exists a positive constant $C$ such that, for any
$\{f_j\}_{j\in\nn}\subset \mathfrak{M}(\bx)$,
$$\lf\|\lf\{\sum_{j\in\nn}[M_{\mathrm{HL}}(f_j)]^u\r\}^{\frac1u}\r\|_{X^{\frac1t}(\bx)}
\leq C
\lf\|\lf\{\sum_{j\in\nn}|f_j|^u\r\}^{\frac1u}\r\|_{X^{\frac1t}(\bx)}.$$
\end{assumption}

 \begin{assumption}\label{a12} Let $X(\bx)$ be a BQBF space. Assume that there exist constants $s_0\in(0,\,\fz)$ and
$q_0\in(s_0,\,\fz)$ such that $X^{\frac{1}{s_0}}(\bx)$ is a BBF space and the Hardy-Littlewood maximal operator $M_{\mathrm{HL}}$ is
bounded on the $\frac{1}{({q_0}/{s_0})'}$-convexification of the associate space
$(X^{\frac{1}{s_0}})'(\bx)$, where $\frac{1}{({q_0}/{s_0})'}+\frac{1}{{q_0}/{s_0}}=1$.
\end{assumption}

\subsection{Global tent space $T^{X}_2(\mathbb{X})$ and local tent space $\mathbf{T}^{X}_2(\mathbb{X})$}\label{nclpx}
In this subsection, we will introduce the local tent space associated with the ball quasi-Banach function space $X(\bx)$. Firstly, we recall the definition of global tent space associated with the ball quasi-Banach function space $X(\bx)$, which is introduced by Lin et al. \cite{lyyy23}.
Let $F$ be a measurable complex-valued function on $\mathbb{X}\times\rr_+$. For any $x\in\mathbb{X}$, we define
$$\mathcal{A}(F)(x):=\lf(\int_0^{\fz}\int_{B(x,\,t)}|F(y,\,t)|^2\frac{d\mu(y)dt}
{V_t(x)t}\r)^{\frac{1}{2}}.$$
For any given BQBF space $X(\bx)$, the tent space $T^{X}_2(\mathbb{X})$ is defined to be the set of all $F$ such that $\mathcal{A}(F)\in X(\mathbb{X})$, equipped with the quasi-norm
$$\|F\|_{T^{X}_2(\mathbb{X})}:=
\lf\|\mathcal{A}(F)\r\|_{X(\mathbb{X})}.$$

For any open set $E\subset\bx$, define the tent $T(E)$ over $E$ by
$$T(E):=\lf\{(x,t)\in\bx\times\rr_+:\mathrm{dist}(x,E^{\complement})\geq t\r\}.$$
\begin{definition}
Let $X(\bx)$ be a BQBF space and $q\in(1,\,\fz)$. A measurable function $A:\bx\times\rr_+\rightarrow\cc$
is said to be a $(T^X_2,\,q)$-atom if there exists a ball $B:=B(x_B,r_B)\subset\bx$ such that
\begin{enumerate}
\item[\rm{(i)}]
$\supp (A):=\{(x,\,t)\in\bx\times\rr_+:A(x,\,t)\neq0\}\subset T(B)$;
\item[\rm{(ii)}] $\|A\|_{T_2^q(\bx)}:=\lf\|\mathcal{A}(A)\r\|_{L^q(\bx)}
\leq [\mu(B)]^\frac{1}{q}\|{\bf 1}_B\|^{-1}_{X(\bx)}$.
\end{enumerate}
Furthermore, if $A$ is a $(T_2^X,\,q)$-atom for any $q\in(1,\,\fz)$, then $A$ is called a $(T_2^X,\,\fz)$-atom.
\end{definition}

The following is the atomic decomposition of the elements of $T^{X}_2(\mathbb{X})$, which
is established by \cite{lyyy23}.
\begin{lemma}\label{lem:global-tent-atomic}
Assume that $X(\bx)$ is a BQBF space satisfying Assumptions \ref{a11} and \ref{a12} for some
$p\in(0,\,\fz)$.
Let $F\in T^X_2(\bx)$ and $s_0\in(0,\,p]$. Then there exists a sequence $\{\lambda_j\}_{j\in\nn}\subset[0,\,\fz)$ and a sequence
$\{A_j\}_{j\in\nn}$ of $(T_2^X,\,\fz)$-atoms associated, respectively, with the balls $\{B_j\}_{j\in\nn}$ such that, for almost every $(x,\,t)\in\bx\times\rr_+$, \begin{align}\label{atom1}
F(x,\,t)=\sum_{j\in\nn}\lambda_jA_j(x,\,t) \ \ \mathrm{in} \ \ T_2^X(\bx)
\end{align}
and
$$\lf\|\lf\{\sum_{j=1}^{\fz}\lf[\frac{\lambda_j{\bf 1}_{B_j}}{\lf\|{\bf 1}_{B_j}\r\|_{X(\bx)}}\r]^{s_0}
\r\}^{\frac{1}{s_0}}\r\|_{X(\bx)}
\lesssim \|F\|_{T_2^X(\bx)},$$
where the implicit positive constant is independent of $F$. Moreover, if $F\in
T_2^2(\bx)\cap T_2^X(\bx)$, then \eqref{atom1} holds true in $T^2_2(\bx)$.
\end{lemma}

Now we introduce the local tent space $\mathbf{T}_2^X(\bx)$ associated with $X(\bx)$.
Let $F$ be a measurable complex-valued function on $\mathbb{X}\times(0,1)$. For any $x\in\mathbb{X}$, we define
$$\mathcal{A}_{\mathrm{loc}}(F)(x):=\lf(\int_0^1\int_{B(x,\,t)}|F(y,\,t)|^2\frac{d\mu(y)dt}
{V_t(x)t}\r)^{\frac{1}{2}}.$$
For any given BQBF space $X(\bx)$, the local tent space $\mathbf{T}_2^X(\mathbb{X})$ is defined to be the set of all $F$ such that $\mathcal{A}_{\mathrm{loc}}(F)\in X(\mathbb{X})$, equipped with the quasi-norm
$$\|F\|_{\mathbf{T}_2^X(\mathbb{X})}:=
\lf\|\mathcal{A}_{\mathrm{loc}}(F)\r\|_{X(\mathbb{X})}.$$
\begin{definition}
Let $X(\bx)$ be a BQBF space and $q\in(1,\,\fz)$. A measurable function $A:\bx\times(0,1)\rightarrow\cc$
is said to be a $(\mathbf{T}_2^X,\,q)$-atom if there exists a ball $B\subset\bx$ such that
\begin{enumerate}
\item[\rm{(i)}]
$\supp (A):=\{(x,\,t)\in\bx\times(0,1):A(x,\,t)\neq0\}\subset T(B)\cap(\bx\times(0,\,1));$
\item[\rm{(ii)}] $\|A\|_{\mathbf{T}_2^q(\bx)}:=\lf\|\mathcal{A}_{\mathrm{loc}}(A)\r\|_{L^q(\bx)}
\leq [\mu(B)]^\frac{1}{q}\|{\bf 1}_B\|^{-1}_{X(\bx)}$.
\end{enumerate}
Furthermore, if $A$ is a $(\mathbf{T}_2^X,\,q)$-atom for any $q\in(1,\,\fz)$, then $A$ is called a $(\mathbf{T}_2^X,\,\fz)$-atom.
\end{definition}
\begin{lemma}
Assume that $X(\bx)$ is a BQBF space satisfying Assumptions \ref{a11} and \ref{a12} for some
$p\in(0,\,\fz)$.
Let $F\in\mathbf{T}^X_2(\bx)$ and $s_0\in(0,\,p]$. Then there exists a sequence $\{\lambda_j\}_{j\in\nn}\subset[0,\,\fz)$ and a sequence
$\{A_j\}_{j\in\nn}$ of $(\mathbf{T}_2^X,\,\fz)$-atoms associated, respectively, with the balls $\{B_j\}_{j\in\nn}$ such that, for almost every $(x,\,t)\in\bx\times(0,1)$, \begin{align}\label{atom}
F(x,\,t)=\sum_{j\in\nn}\lambda_jA_j(x,\,t) \ \ \mathrm{in} \ \ \mathbf{T}_2^X(\bx)
\end{align}
and
$$\lf\|\lf\{\sum_{j=1}^{\fz}\lf[\frac{\lambda_j{\bf 1}_{B_j}}{\lf\|{\bf 1}_{B_j}\r\|_{X(\bx)}}\r]^{s_0}
\r\}^{\frac{1}{s_0}}\r\|_{X(\bx)}
\lesssim \|F\|_{\mathbf{T}_2^X(\bx)},$$
where the implicit positive constant is independent of $F$ . Moreover, if $F\in
\mathbf{T}_2^2(\bx)\cap \mathbf{T}_2^X(\bx)$, then \eqref{atom} holds true in $\mathbf{T}^2_2(\bx)$.
\end{lemma}

\begin{proof}
For a measurable function $G$ on $\bx\times(0,\infty)$, put
\[
 \mathfrak P G(y,t):=G(y,t){\bf 1}_{(0,1)}(t),
 \qquad (y,t)\in\bx\times(0,\infty).
\]
Then, for every $x\in\bx$,
\[
 \mathcal A_{\rm loc}(\mathfrak P G)(x)
 =
 \left(
 \int_0^1\int_{B(x,t)}
 |G(y,t)|^2\,\frac{d\mu(y)\,dt}{V_t(x)t}
 \right)^{1/2}
 \leq
 \mathcal A(G)(x).
\]
Consequently, the truncation operator $\mathfrak P$ is bounded from
$T^X_2(\bx)$ to $\mathbf T^X_2(\bx)$, and also from $T^2_2(\bx)$ to
$\mathbf T^2_2(\bx)$.

Let $F\in\mathbf T^X_2(\bx)$. Define
\[
 \widetilde F(y,t):=
 \begin{cases}
 F(y,t), & 0<t<1,\\
 0,      & t\geq 1.
 \end{cases}
\]
Then
$
 \mathcal A(\widetilde F)(x)
 =
 \mathcal A_{\rm loc}(F)(x),
  x\in\bx,
$
and hence
$
 \|\widetilde F\|_{T^X_2(\bx)}
 =
 \|F\|_{\mathbf T^X_2(\bx)}.
$
Thus $\widetilde F\in T^X_2(\bx)$. By the atomic decomposition of
$T^X_2(\bx)$ in Lemma \ref{lem:global-tent-atomic}, there exist
$\{\lambda_j\}_{j\in\nn}\subset[0,\infty)$ and a sequence
$\{A_j\}_{j\in\nn}$ of $(T^X_2,\infty)$-atoms associated, respectively,
with balls $\{B_j\}_{j\in\nn}$ such that
\[
 \widetilde F
 =
 \sum_{j\in\nn}\lambda_jA_j
 \quad\text{in } T^X_2(\bx),
\]
and
\[
\left\|
\left\{
\sum_{j\in\nn}
\left[
\frac{\lambda_j{\bf 1}_{B_j}}
{\|{\bf 1}_{B_j}\|_{X(\bx)}}
\right]^{s_0}
\right\}^{1/s_0}
\right\|_{X(\bx)}
\lesssim
\lf\|\widetilde F\r\|_{T^X_2(\bx)}
=
\|F\|_{\mathbf T^X_2(\bx)}.
\]

For each $j\in\nn$, set
\[
 \widetilde A_j:=\mathfrak P A_j
 =
 A_j{\bf 1}_{\bx\times(0,1)}.
\]
Since $A_j$ is supported in $T(B_j)$, we have
\[
 \supp(\widetilde A_j)
 \subset
 T(B_j)\cap(\bx\times(0,1)).
\]
Moreover, for any $q\in(1,\infty)$,
\[
 \lf\|\widetilde A_j\r\|_{\mathbf T^q_2(\bx)}
 =
 \lf\|\mathcal A_{\rm loc}(\widetilde A_j)\r\|_{L^q(\bx)}
 \leq
 \lf\|\mathcal A(A_j)\r\|_{L^q(\bx)}
 =
 \lf\|A_j\r\|_{T^q_2(\bx)}
 \leq
 [\mu(B_j)]^{1/q}\lf\|{\bf 1}_{B_j}\r\|_{X(\bx)}^{-1}.
\]
Thus each $\widetilde A_j$ is a $(\mathbf T^X_2,\infty)$-atom
associated with $B_j$.

Now, for almost every $(x,t)\in\bx\times(0,1)$,
\[
 F(x,t)
 =
 \widetilde F(x,t)
 =
 \sum_{j\in\nn}\lambda_jA_j(x,t)
 =
 \sum_{j\in\nn}\lambda_j\widetilde A_j(x,t).
\]
It remains to verify the convergence in $\mathbf T^X_2(\bx)$. For every
$N\in\nn$,
\[
\begin{aligned}
\left\|
F-\sum_{j=1}^N\lambda_j\widetilde A_j
\right\|_{\mathbf T^X_2(\bx)}
&=
\left\|
\mathcal A_{\rm loc}
\left(
\mathfrak P\left[
\widetilde F-\sum_{j=1}^N\lambda_jA_j
\right]
\right)
\right\|_{X(\bx)}
\leq
\left\|
\mathcal A
\left(
\widetilde F-\sum_{j=1}^N\lambda_jA_j
\right)
\right\|_{X(\bx)}
\\
&=
\left\|
\widetilde F-\sum_{j=1}^N\lambda_jA_j
\right\|_{T^X_2(\bx)}
\to0
\end{aligned}
\]
as $N\to\infty$. Hence
\[
 F=\sum_{j\in\nn}\lambda_j\widetilde A_j
 \quad\text{in } \mathbf T^X_2(\bx).
\]
Renaming $\widetilde A_j$ as $A_j$, we obtain the desired atomic
decomposition and the required estimate of the coefficients.

Finally, assume in addition that
$F\in\mathbf T^2_2(\bx)\cap\mathbf T^X_2(\bx)$. Then
$\widetilde F\in T^2_2(\bx)\cap T^X_2(\bx)$. By the last assertion of
Lemma \ref{lem:global-tent-atomic}, we obtain
\[
 \widetilde F
 =
 \sum_{j\in\nn}\lambda_jA_j
 \quad\text{in } \ \ T^2_2(\bx).
\]
Using again the boundedness of $\mathfrak P$ from $T^2_2(\bx)$ to
$\mathbf T^2_2(\bx)$, we get
\[
\begin{aligned}
\left\|
F-\sum_{j=1}^N\lambda_j\widetilde A_j
\right\|_{\mathbf T^2_2(\bx)}
&\leq
\left\|
\widetilde F-\sum_{j=1}^N\lambda_jA_j
\right\|_{T^2_2(\bx)}
\to0.
\end{aligned}
\]
Therefore the decomposition also holds in $\mathbf T^2_2(\bx)$.
This finishes the proof.
\end{proof}

\section{The local Hardy space $h^{X}_{L}(\mathbb{X})$}\label{s3}
This section is devoted to the introduction of the local Hardy space \(h_L^X(\bx)\) associated with a ball quasi-Banach function space $X$ and a non-negative self-adjoint operator \(L\) on spaces of homogeneous type. We first formulate the standing assumptions used throughout this section.

\begin{assumption}\label{a2.5}
\rm Let $(\bx,\,d,\,\mu)$ be a space of homogeneous  type (see Section \ref{homogene} for definition). We assume that:
\begin{enumerate}
\item[\rm $(\mathbf{H}_1)$]
 $L$ is a non-negative self-adjoint densely--defined operator on $L^2(\bx)$;
\item[\rm $(\mathbf{H}_2)$] The kernel of $e^{-tL}$, denoted by $P_t(x,\,y)$, is a measurable function on
$\bx\times\bx$ and satisfies the Gaussian upper bound, that is, there are constants $C,c>0$ such that for any $t>0$ and $x$, $y\in\bx$,
$$\lf|P_t(x,\,y)\r|\leq \frac{C}{\mu(B(x,\sqrt{t}))}
\exp\lf(-c\frac{d(x, y)^2}{t}\r),$$
\end{enumerate}
where we denote by $B(x,\,r):=\{y\in\bx:d(x,\,y)<r\}$ the ball centered at $x\in \bx$ with radius $r>0$.
\end{assumption}

There exist many operators satisfying the above assumption, we give the following operators:
\begin{example}
Let $\mathbb X=\mathbb R^d$ equipped with the Euclidean metric and the
Lebesgue measure.
\begin{enumerate}
\item [{\rm(i)}]
Let
$
L=-\Delta_{\rr^d}
$
be the non-negative Laplacian on $L^2(\mathbb R^d)$. Then $L$ is a
non-negative self-adjoint operator. Moreover, the heat kernel of
$e^{-tL}=e^{t\Delta_{\rr^d}}$ is given by
\[
P_t(x,y)=\frac{1}{(4\pi t)^{d/2}}
\exp\left(-\frac{|x-y|^2}{4t}\right),
\qquad x,y\in\mathbb R^d,\ t>0.
\]
Since
$
\mu(B(x,\sqrt t))\sim t^{d/2},
$
we have
\[
|P_t(x,y)|
\lesssim
\frac{1}{\mu(B(x,\sqrt t))}
\exp\left(-c\frac{|x-y|^2}{t}\right).
\]
Thus $L=-\Delta_{\rr^d}$ satisfies Assumption \ref{a2.5}.
\item [{\rm(ii)}] Let $V$ be a non-negative locally integrable function
on $\mathbb R^d$, and consider the Schr\"odinger operator
$
L=-\Delta+V
$
defined via the closed quadratic form
\[
\mathfrak q(f,g)
:=
\int_{\mathbb R^d}\nabla f(x)\cdot \overline{\nabla g(x)}\,dx
+
\int_{\mathbb R^d}V(x)f(x)\overline{g(x)}\,dx.
\]
Then $L$ is a non-negative self-adjoint operator on $L^2(\mathbb R^d)$.
Moreover, by the domination property of the heat semigroup,
\[
0\leq e^{-tL}(x,y)
\leq
e^{t\Delta_{\rr^d}}(x,y)
=
\frac{1}{(4\pi t)^{d/2}}
\exp\left(-\frac{|x-y|^2}{4t}\right).
\]
Hence
\[
|P_t(x,y)|
\lesssim
\frac{1}{\mu(B(x,\sqrt t))}
\exp\left(-c\frac{|x-y|^2}{t}\right),
\qquad x,y\in\mathbb R^d,\ t>0.
\]
Therefore the Schr\"odinger operator $L=-\Delta_{\rr^d}+V$ with $V\geq0$
satisfies Assumption \ref{a2.5}.
\end{enumerate}
\end{example}

For any $x\in\bx$ and $f\in L^2(\bx)$, we consider the area square integral function $S_L(f)$ defined by
$$S_L(f)(x):=\lf(\int_0^{\fz}\int_{B(x,\,t)}\lf|t^2Le^{-t^2L}(f)(y)\r|^2
\frac{d\mu(y)dt}{V_t(x)t}\r)^{\frac{1}{2}}$$
and the local area square integral function $S_L^{\mathrm{loc}}$ defined by
$$S_L^{\mathrm{loc}}(f)(x):=\lf(\int_0^1\int_{B(x,\,t)}\lf|t^2Le^{-t^2L}(f)(y)\r|^2
\frac{d\mu(y)dt}{V_t(x)t}\r)^{\frac{1}{2}}$$

The Hardy space $H_L^X(\bx)$ associated with $L$ and $X(\bx)$ is defined by Lin et al. \cite{lyyy23}.
A function $f\in L^2(\bx)$ is in $\mathbf{H}_L^X(\bx)$
 when $S_L(f)\in X(\bx)$. We consider the quasi-norm $\|\cdot\|_{H_L^X(\bx)}$
 given by
$$\|f\|_{H_L^X(\bx)}:=\lf\|S_L(f)\r\|_{X(\bx)}, \ \ \ f\in\mathbf{H}_L^X(\bx).$$
Define the Hardy space $H_L^X(\bx)$ as the completion of $\mathbf{H}_L^X(\bx)$ with respect
to $\|\cdot\|_{H_L^X(\bx)}$.

Next we introduce the local Hardy space $h_L^X(\bx)$ associated with $L$ in the following way.
A function $f\in L^2(\bx)$ is in $\mathbf{h}_L^X(\bx)$
 when $S^{\mathrm{loc}}_L(f)\in X(\bx)$
 and $S_I(e^{-L}f)\in X(\bx)$. Here, $S_I$ denotes the area square integral function associated
with the identity operator $I$. We consider the quasi-norm $\|\cdot\|_{h_L^X(\bx)}$
 given by
$$\|f\|_{h_L^X(\bx)}:=\lf\|S^{\mathrm{loc}}_L(f)\r\|_{X(\bx)}+\lf\|S_I(e^{-L}f)
\r\|_{X(\bx)}, \ \ \ f\in\mathbf{h}_L^X(\bx).$$
Define the local Hardy space $h_L^X(\bx)$ as the completion of $\mathbf{h}_L^X(\bx)$ with respect
to $\|\cdot\|_{h_L^X(\bx)}$.

In what follows, for a ball $B=B(x_B,r_B)$ and $j\in\zz_+$, put
\[
S_0(B):=B,\qquad S_j(B):=2^jB\setminus 2^{j-1}B.
\]
To characterize  the local Hardy space $h_L^X(\bx)$ by using molecules,
we introduce two types of local molecules:
a measurable function
$m\in L^2(\bx)$ is a local $(X,\,2,\,M,\,\varepsilon)$-molecule associated with the ball $B=B(x_B,\,r_B)$
\begin{enumerate}
\item [(a)] of $(\mathrm{I})$-type: when $r_B\geq1$ and
$\|m\|_{L^2(S_j(B))}\leq2^{-j\varepsilon}\mu(2^jB)^{\frac{1}{2}}
\|{\bf 1}_{2^jB}\|_{X(\bx)}^{-1}$, for every
$j\in\zz_+$.
\item [(b)] of $(\mathrm{II})$-type: $r_B\in(0, 1)$ and there exists \(b\in D(L^M)\) such that
$m=L^M b$ and,
for every $k\in\{0,\,1,\ldots,M\}$,
$\lf\|L^k(b)\r\|_{L^2(S_j(B))}\leq2^{-j\varepsilon}r_B^{2(M-k)}\mu(2^jB)^{\frac{1}{2}}
\|{\bf 1}_{2^jB}\|_{X(\bx)}^{-1}$, $j\in\zz_+.$
\end{enumerate}
\begin{definition}
We say that $f\in L^2(\bx)$ is in the space $\mathbf{h}^X_{L,\mathrm{mol},M,\varepsilon}(\bx)$ if, for every $j\in\nn$,
there exist $\lambda_j>0$ and a local $(X,\,2,\,M,\,\varepsilon)$-molecule $m_j$ associated with the ball $B_j$
such that $$f=\sum_{j\in\nn}\lambda_jm_j \ \ \mathrm{in} \ \ L^2(\bx)$$ and
$$\|f\|_{\mathbf{h}^X_{L,\mathrm{mol},M,\varepsilon}(\bx)}:=
\inf\lf\|\lf\{\sum_{j=1}^{\fz}\lf[\frac{\lambda_j{\bf 1}_{B_j}}{\lf\|{\bf 1}_{B_j}\r\|_{X(\bx)}}\r]^{s_0}
\r\}^{\frac{1}{s_0}}\r\|_{X(\bx)}
<\fz.$$
The local molecular Hardy space $h^X_{L,\mathrm{mol},M,\varepsilon}(\bx)$ is the completion of $\mathbf{h}^X_{L,\mathrm{mol},M,\varepsilon}(\bx)$ with respect to the
quasi-norm $\|\cdot\|_{h^X_{L,\mathrm{mol},M,\varepsilon}(\bx)}$.
\end{definition}

The following are the main results of this paper.
\begin{theorem}\label{mol}
Let $L$ be a non-negative self-adjoint operator on $L^2(\bx)$ satisfying Assumption \ref{a2.5}, and
Let \(X(\bx)\) be a BQBF space satisfying both Assumptions
\(\ref{a11}\) and \(\ref{a12}\) for some $p\in(0,\fz)$, \(s_0\in(0,\min\{p,1\}]\) and \(q_0\in(s_0,2]\).
\begin{enumerate}
\item [{\rm(i)}]
If $\varepsilon>d(\frac{1}{s_0}-\frac{1}{q_0})$ and $M>\frac{d}{2}(\frac{2}{s_0}-\frac{1}{2}-\frac{1}{q_0})$, then we have
$$h^X_{L,\mathrm{mol},M,\varepsilon}(\bx)\subseteq
h^X_{L}(\bx);$$
\item [{\rm(ii)}]
If $\varepsilon>0$ and $M\in\nn$, then we have
$$h^X_{L}(\bx)\subseteq h^X_{L,\mathrm{mol},M,\varepsilon}(\bx).$$
\end{enumerate}
\end{theorem}

\begin{theorem}\label{eq0}
Let $L$ be a non-negative self-adjoint operator on $L^2(\bx)$ satisfying the Assumption \ref{a2.5}, and
$X(\bx)$ be a BQBF space satisfying both Assumptions \ref{a11} and \ref{a12}  for some
$p\in(0,\fz)$, $s_0\in(0,\,\min\{p,1\}]$ and $q_0\in(s_0,\,2]$. If $\inf\sigma(L)>0$, then we have
$$h^X_L(\bx)=H^X_{L}(\bx),$$
where \(\sigma(L)\) denotes the spectrum of \(L\) on \(L^2(\bx)\).
\end{theorem}

\begin{theorem}\label{eq1}
Let $L$ be a non-negative self-adjoint operator on $L^2(\bx)$ satisfying the Assumption \ref{a2.5}, and
$X(\bx)$ be a BQBF space satisfying Assumptions \ref{a11} and \ref{a12}  for some
$p\in(0,\fz)$, $s_0\in(0,\,\min\{p,1\}]$ and $q_0\in(s_0,\,2]$. Then we have
$$h^X_L(\bx)=H^X_{L+I}(\bx)$$
algebraically and topologically.
\end{theorem}

The proofs Theorems \ref{mol}, \ref{eq0} and \ref{eq1} is given in Section \ref{s5}.

\section{Global Hardy spaces associated with the
identity operator}\label{s4}
In this section, we investigate the auxiliary Hardy space $H_{I}^X(\bx)$ associated with the
identity operator $I$. This part will be used to treat the low-frequency term in
the definition of the local Hardy space $h_L^X(\bx)$.
\begin{definition}
Let $M\in \nn$ and $\varepsilon>0$.
 \begin{enumerate}
\item [{\rm(i)}]
A function
$a\in L^2(\bx)$ is called an $(X,2,M)_I$-atom associated with a ball
$B=B(x_B,r_B)$ if
$
 \operatorname{supp}a\subset B
$
and, for every $k\in\{0,1,\ldots,M\}$,
\begin{equation}\label{e:I-atom-size}
 \|a\|_{L^2(\bx)}
 \leq
 r_B^{2k}\mu(B)^{1/2}\|{\bf 1}_B\|_{X(\bx)}^{-1}.
\end{equation}

\item [{\rm(ii)}]
A function $m\in L^2(\bx)$ is called an
$(X,2,M,\varepsilon)_I$-molecule associated with $B$ if, for every
$i\in\zz_+$ and every $k\in\{0,1,\ldots,M\}$,
\begin{equation}\label{e:I-molecule-size}
 \|m\|_{L^2(S_i(B))}
 \leq
 2^{-i\varepsilon}r_B^{2k}
 \mu(2^iB)^{1/2}\lf\|{\bf 1}_{2^iB}\r\|_{X(\bx)}^{-1}.
\end{equation}
\end{enumerate}
Clearly, every $(X,2,M)_I$-atom is an
$(X,2,M,\varepsilon)_I$-molecule for every $\varepsilon>0$.
The atomic Hardy space $H_{I,\textrm{at},M}^{X}(\bx)$ and the quasi-norm $\|\cdot\|_{H_{I,\textrm{at},M}^{X}(\bx)}$, and
the molecular Hardy space $H_{I,\textrm{mol},M,\varepsilon}^{X}(\bx)$ and the quasi-norm $\|\cdot\|_{H_{I,\textrm{mol},M,\varepsilon}^{X}(\bx)}$ are defined in the usual way.
\end{definition}

%
\begin{remark}
Notice that for every $M\in \nn$, if $a$ is a $(X, 2,M)_I$-atom, then $a$ is also a $(X,2,M,\varepsilon)_I$-molecule, for each $\varepsilon >0$. Thus,
$H_{I,\textrm{at},M}^{X}(\bx)$ is continuously contained in $H_{I,\textrm{mol},M,\varepsilon}^{X}(\bx)$, for every $M\in \zz_+$ and $\varepsilon >0$.
\end{remark}

For $M\in\nn$, set
\[
 c_M^{-1}:=
 \int_0^\infty t^{2M+2}e^{-2t^2}\,\frac{dt}{t},
\]
and define, initially on $T_2^2(\bx)$,
\begin{equation}\label{e:PiM-I-def}
 \Pi_M(F)(x):=
 c_M\int_0^\infty t^{2M}e^{-t^2}F(x,t)\,\frac{dt}{t}.
\end{equation}


\begin{lemma}\label{PiM}
Let $M\in\nn$. Then $\Pi_M$ is bounded from $T_2^2(\bx)$ to
$L^2(\bx)$. Moreover, if $A$ is a $(T_2^X,\infty)$-atom associated with a
ball $B$, then $C\Pi_M(A)$ is an $(X,2,M)_I$-atom associated with $B$, where
$C>0$ is independent of $A$ and $B$. Consequently, $\Pi_M$ extends to a
bounded operator from $T_2^X(\bx)$ to $H_{I,\mathrm{at},M}^X(\bx)$ and hence
to $H_{I,\mathrm{mol},M,\varepsilon}^X(\bx)$ for every $\varepsilon>0$.
\end{lemma}

\begin{proof}
By Lemma \ref{l2.0}, for every $t>0$ and $x,y\in\bx$ with $d(x,y)<t$,
$V_t(x)\sim V_t(y)$. Hence
\[
 \|F\|_{T_2^2(\bx)}^2
 \sim
 \int_0^\infty\int_\bx |F(y,t)|^2\,d\mu(y)\frac{dt}{t}.
\]
Therefore, by H\"older's inequality, we obtain
\begin{align*}
	\|\Pi_M(F)\|_{L^2(\bx)}^2&\leq \int_{\bx} \left(\int_0^\infty t^{2M}e^{-t^2}|F(y,t)|\frac{dt}{t}\right)^2 d\mu(y)\\
	&\leq \int_{\bx} \int_0^\infty \lf(t^{2M}e^{-t^2}\r)^2 \frac{dt}{t}\int_0^\infty |F(y,t)|^2\frac{dt}{t} d\mu(y)\\
	&\lesssim \int_{\bx}  \int_0^\infty \int_{B(y,t)} |F(y,t)|^2d\mu(x)\frac{dtd\mu(y)}{V_t(x)t} \thicksim\|F\|_{T^2_2(\bx)}^2.
	\end{align*}
Thus $\Pi_M:T_2^2(\bx)\to L^2(\bx)$ is bounded.

Let now $A$ be a $(T_2^X,\infty)$-atom associated with
$B=B(x_B,r_B)$. Since $\operatorname{supp}A\subset T(B)$, we have
$\operatorname{supp}\Pi_M(A)\subset B$. Moreover, for every
$k\in\{0,1,\ldots,M\}$,
\[
 \left(\int_0^{r_B}t^{4M}e^{-2t^2}\,\frac{dt}{t}\right)^{1/2}
 \lesssim r_B^{2k}.
\]
Using again the equivalence between the $T_2^2$-norm and the
$L^2(\bx\times\rr_+,d\mu dt/t)$-norm, together with the atom condition, we
obtain
\[
\begin{aligned}
 \|\Pi_M(A)\|_{L^2(\bx)}
 &\lesssim
 \left(\int_0^{r_B}t^{4M}e^{-2t^2}\,\frac{dt}{t}\right)^{1/2}
 \left(\int_0^{r_B}\|A(\cdot,t)\|_{L^2(\bx)}^2\,\frac{dt}{t}\right)^{1/2} \\
 &\lesssim
 r_B^{2k}\|A\|_{T_2^2(\bx)}
 \lesssim
 r_B^{2k}\mu(B)^{1/2}\|{\bf 1}_B\|_{X(\bx)}^{-1}.
\end{aligned}
\]
This proves that $C\Pi_M(A)$ is an $(X,2,M)_I$-atom.

Finally, if $F\in T_2^X(\bx)\cap T_2^2(\bx)$, then the atomic decomposition of
$T_2^X(\bx)$ gives
\[
 F=\sum_{j\in\nn}\lambda_jA_j
 \ \ \hbox{in } \ \ T_2^2(\bx),
\]
where $A_j$ are $(T_2^X,\infty)$-atoms associated with balls $B_j$, and
\[
 \left\|
 \left\{
 \sum_{j\in\nn}
 \left[
 \frac{|\lambda_j|{\bf 1}_{B_j}}
      {\|{\bf 1}_{B_j}\|_{X(\bx)}}
 \right]^{s_0}
 \right\}^{1/s_0}
 \right\|_{X(\bx)}
 \lesssim
 \|F\|_{T_2^X(\bx)}.
\]
Since $\Pi_M$ is bounded from $T_2^2(\bx)$ to $L^2(\bx)$,
\[
 \Pi_M(F)=\sum_{j\in\nn}\lambda_j\Pi_M(A_j)
 \qquad\hbox{in }L^2(\bx).
\]
The previous atom estimate then implies
\[
 \|\Pi_M(F)\|_{H_{I,\mathrm{at},M}^X(\bx)}
 \lesssim
 \|F\|_{T_2^X(\bx)}.
\]
The extension follows by completion. Since every $(X,2,M)_I$-atom is an
$(X,2,M,\varepsilon)_I$-molecule, the asserted boundedness into the molecular
space is immediate.
\end{proof}

\begin{proposition}\label{prop:HI-atomic-inclusion}
Let $M\in\nn$ and $\varepsilon>0$. Then
\[
 H_I^X(\bx)
 \subset
 H_{I,\mathrm{at},M}^X(\bx)
 \subset
 H_{I,\mathrm{mol},M,\varepsilon}^X(\bx)
\]
algebraically and continuously.
\end{proposition}

\begin{proof}
Let $f\in\mathbf H_I^X(\bx)$ and put
$
 F_f(y,t):=t^2e^{-t^2}f(y),
(y,t)\in\bx\times\rr_+.
$
Then
\[
 \lf\|F_f\r\|_{T_2^X(\bx)}=\lf\|S_I(f)\r\|_{X(\bx)}=
 \|f\|_{H_I^X(\bx)}.
\]
Moreover, if $f\in L^2(\bx)$, then $F_f\in T_2^2(\bx)$ and, by the choice of
$c_M$, we have
\[
 \Pi_M(F_f)=f
\ \ \hbox{in } \ \ L^2(\bx).
\]
Lemma \ref{PiM} gives
\[
 \|f\|_{H_{I,\mathrm{at},M}^X(\bx)}
 \lesssim
 \lf\|F_f\r\|_{T_2^X(\bx)}
 =
 \|f\|_{H_I^X(\bx)}.
\]
Passing to the completion yields the first continuous embedding. The second
embedding follows from the fact that every $(X,2,M)_I$-atom is an
$(X,2,M,\varepsilon)_I$-molecule.
\end{proof}

To prove the reverse result, we need to establish a technical lemmas as follows.

\begin{lemma}\label{l2.8}
Let $X(\bx)$ be a BQBF space satisfying Assumption \ref{a12} for some
$s_0\in(0,\,1]$ and $q_0\in(s_0,\,2]$, and $T$ be a bounded sublinear operator on $L^2(\bx)$, and satisfy,
for all
\(f,g\in L^2(\bx)\) and all \(\alpha\in\mathbb C\),
\begin{align*}
 T(f+g)\leq T(f)+T(g)
 \quad\text{and}\quad
 T(\alpha f)=|\alpha|T(f).
\end{align*}
Assume that $f=\sum_{j\in\nn}\lambda_jm_j$ in $L^2(\bx)$ with $\lambda_j>0$, $m_j\in L^2(\bx)$, and $\{B_j\}_{j\in\nn}\subset\bx$ be a sequence of balls such that
$$\lf\|\lf\{\sum_{j\in\nn}\lf[\frac{|\lambda_{j}|
{\bf 1}_{B_j}}{\|{\bf 1}_{B_j}\|_{X(\bx)}}\r]
^{s_0}\r\}^{\frac{1}{s_0}}\r\|_{X(\bx)}<\fz.$$
If there exist $C>0$ and $\delta>d(\frac{1}{s_0}-\frac{1}{q_0})$ such that, for all
\(j\in\nn\) and \(i\in\zz_+\),
\begin{align*}
\lf\|T(m_j)\r\|_{L^2(S_i(B_j))}
\leq C 2^{-i\delta}\mu(2^iB_j)^{\frac{1}{2}}\lf\|{\bf 1}_{2^iB_j}\r\|_{X(\bx)}^{-1}.
\end{align*}
Then we have
$T(f)\in X(\bx)$ and
$$\lf\|T(f)\r\|_{X(\bx)}\leq C
\lf\|\lf\{\sum_{j\in\nn}\lf[\frac{|\lambda_{j}|
{\bf 1}_{B_j}}{\|{\bf 1}_{B_j}\|_{X(\bx)}}\r]
^{s_0}\r\}^{\frac{1}{s_0}}\r\|_{X(\bx)}<\fz.$$
\end{lemma}


\begin{proof}
Let
$
 Y:=X^{1/s_0}(\bx)$, $ r:=q_0/s_0>1.
$
By Assumption \ref{a12}, we know that \(Y\) is a ball Banach function space and
\(M_{\rm HL}\) is bounded on \((Y')^{1/r'}(\bx)\). We shall use the following
two standard consequences of this assumption. First, if
\(B\subset \widetilde B\), by Assumption \ref{a12}, we obtain
\begin{equation}\label{e:ind-comp-lem39-short}
 \frac{\|{\bf 1}_B\|_{X(\bx)}}{\|{\bf 1}_{\widetilde B}\|_{X(\bx)}}
 \lesssim
 \left[\frac{\mu(B)}{\mu(\widetilde B)}\right]^{1/q_0}
 \quad
\mathrm{and} \quad
 \frac{\|{\bf 1}_{\widetilde B}\|_{X(\bx)}}{\|{\bf 1}_B\|_{X(\bx)}}
 \lesssim
 \left[\frac{\mu(\widetilde B)}{\mu(B)}\right]^{1/s_0}.
\end{equation}
Second, if \(\operatorname{supp}a_j\subset B_j\) and
\[
 \lf\|a_j\r\|_{L^{q_0}(\bx)}
 \leq
 \mu(B_j)^{1/q_0}\lf\|{\bf 1}_{B_j}\r\|_{X(\bx)}^{-1},
\]
then
\begin{equation}\label{e:atom-est-lem39-short}
 \left\|
 \sum_j\lambda_j^{s_0}\lf|a_j\r|^{s_0}
 \right\|_Y
 \lesssim
 \left\|
 \sum_j
 \left[
 \frac{\lambda_j{\bf 1}_{B_j}}
      {\|{\bf 1}_{B_j}\|_{X(\bx)}}
 \right]^{s_0}
 \right\|_Y .
\end{equation}
Indeed, \eqref{e:atom-est-lem39-short} follows by testing against
\(0\leq g\in Y'\), \(\|g\|_{Y'}\leq1\), see [\cite{syy22}, Proposition 2.14] for more details.

Moreover, for every \(i\in\zz_+\),
\begin{equation}\label{e:dil-est-lem39-short}
 \left\|
 \sum_j
 \left[
 \frac{\lambda_j{\bf 1}_{2^iB_j}}
      {\|{\bf 1}_{2^iB_j}\|_{X(\bx)}}
 \right]^{s_0}
 \right\|_Y
 \lesssim
 2^{id(1-s_0/q_0)}
 \left\|
 \sum_j
 \left[
 \frac{\lambda_j{\bf 1}_{B_j}}
      {\|{\bf 1}_{B_j}\|_{X(\bx)}}
 \right]^{s_0}
 \right\|_Y .
\end{equation}
This is obtained from the doubling
property and \eqref{e:ind-comp-lem39-short}, see also [\cite{syy2022}, Proposition 4.8(i)].

For \(N\in\nn\), set
\[
 f_N:=\sum_{j=1}^N\lambda_jm_j .
\]
Since \(T\) is sublinear and bounded on \(L^2(\bx)\),
\[
 Tf_N\leq\sum_{j=1}^N\lambda_jT(m_j)
 \ \ \hbox{in } \ \ L^2(\bx).
\]
As \(0<s_0\leq1\),
\[
 \lf|Tf_N\r|^{s_0}
 \leq
 \sum_{j=1}^N\lambda_j^{s_0}|T(m_j)|^{s_0}
 =
 \sum_{i=0}^\infty
 \sum_{j=1}^N
 \lambda_j^{s_0}|T(m_j)|^{s_0}{\bf 1}_{S_i(B_j)} .
\]
Hence, using the Banach lattice property of \(Y\),
\begin{equation}\label{e:start-lem39-short}
 \|Tf_N\|_{X(\bx)}^{s_0}
 \leq
 \sum_{i=0}^\infty
 \left\|
 \sum_{j=1}^N
 \lambda_j^{s_0}|T(m_j)|^{s_0}{\bf 1}_{S_i(B_j)}
 \right\|_Y .
\end{equation}

Fix \(i\in\zz_+\). By the assumed annular estimate,
\[
 \lf\|T(m_j){\bf 1}_{S_i(B_j)}\r\|_{L^2(\bx)}
 \lesssim
 2^{-i\delta}
 \mu(2^iB_j)^{1/2}
 \lf\|{\bf 1}_{2^iB_j}\r\|_{X(\bx)}^{-1}.
\]
Since \(q_0\leq2\), H\"{o}lder's inequality gives
\[
 \lf\|T(m_j){\bf 1}_{S_i(B_j)}\r\|_{L^{q_0}(\bx)}
 \lesssim
 2^{-i\delta}
 \mu(2^iB_j)^{1/q_0}
 \lf\|{\bf 1}_{2^iB_j}\r\|_{X(\bx)}^{-1}.
\]
Applying \eqref{e:atom-est-lem39-short} with \(2^iB_j\) in place of \(B_j\),
and then \eqref{e:dil-est-lem39-short}, we obtain
\[
\begin{aligned}
 \left\|
 \sum_{j=1}^N
 \lf|\lambda_j\r|^{s_0}|T(m_j)|^{s_0}{\bf 1}_{S_i(B_j)}
 \right\|_Y
 \lesssim
 2^{-is_0\delta}
 \left\|
 \sum_{j=1}^N
 \left[
 \frac{\lf|\lambda_j\r|{\bf 1}_{2^iB_j}}
      {\lf\|{\bf 1}_{2^iB_j}\r\|_{X(\bx)}}
 \right]^{s_0}
 \right\|_Y
 \quad\lesssim
 2^{-is_0[\delta-d(\frac1{s_0}-\frac1{q_0})]}
 \left\|
 \sum_{j=1}^N
 \left[
 \frac{\lf|\lambda_j\r|{\bf 1}_{B_j}}
      {\lf\|{\bf 1}_{B_j}\r\|_{X(\bx)}}
 \right]^{s_0}
 \right\|_Y .
\end{aligned}
\]
Since
$
 \delta>d(\frac1{s_0}-\frac1{q_0}),
$
summing over \(i\) in \eqref{e:start-lem39-short} yields
\[
 \|Tf_N\|_{X(\bx)}^{s_0}
 \lesssim
 \left\|
 \sum_{j=1}^N
 \left[
 \frac{\lambda_j{\bf 1}_{B_j}}
      {\|{\bf 1}_{B_j}\|_{X(\bx)}}
 \right]^{s_0}
 \right\|_Y
 \leq
 \left\|
 \sum_{j=1}^\infty
 \left[
 \frac{\lambda_j{\bf 1}_{B_j}}
      {\|{\bf 1}_{B_j}\|_{X(\bx)}}
 \right]^{s_0}
 \right\|_Y .
\]

Finally, \(f_N\to f\) in \(L^2(\bx)\), and hence \(Tf_N\to Tf\) in
\(L^2(\bx)\). Passing to a subsequence if necessary, \(Tf_N\to Tf\) almost
everywhere. By the Fatou property of \(Y\),
\[
\begin{aligned}
 \|Tf\|_{X(\bx)}^{s_0}
 =
 \lf\||Tf|^{s_0}\r\|_Y
 &\leq
 \liminf_{N\to\infty}
 \lf\|\lf|Tf_N\r|^{s_0}\r\|_Y
 \lesssim
 \left\|
 \sum_{j=1}^\infty
 \left[
 \frac{\lambda_j{\bf 1}_{B_j}}
      {\lf\|{\bf 1}_{B_j}\r\|_{X(\bx)}}
 \right]^{s_0}
 \right\|_Y .
\end{aligned}
\]
Taking the power \(1/s_0\), we obtain
\[
 \|Tf\|_{X(\bx)}
 \lesssim
 \left\|
 \left\{
 \sum_{j=1}^\infty
 \left[
 \frac{\lambda_j{\bf 1}_{B_j}}
      {\lf\|{\bf 1}_{B_j}\r\|_{X(\bx)}}
 \right]^{s_0}
 \right\}^{1/s_0}
 \right\|_{X(\bx)} .
\]
This proves the desired estimate.
\end{proof}

\begin{proposition}\label{prop:HI-mol-to-HI}
Let \(X(\bx)\) be a ball quasi-Banach function space satisfying Assumption
\ref{a12} for some \(s_0\in(0,1]\) and \(q_0\in(s_0,2]\). Let
\(M\in\nn\) and \(\varepsilon>0\) satisfy
$
 \varepsilon>d(\frac1{s_0}-\frac1{q_0})
$
and
$
 M>\frac d2(\frac2{s_0}-\frac12-\frac1{q_0}).
$
Then there exists a positive constant \(C\) such that, if
\[
 f=\sum_{j\in\nn}\lambda_jm_j
 \ \ \mathrm{in} \ \ L^2(\bx),
\]
where, for every \(j\in\nn\), \(\lambda_j>0\) and \(m_j\) is an
\((X,2,M,\varepsilon)_I\)-molecule associated with a ball \(B_j\), and
\[
 \left\|
 \left\{
 \sum_{j\in\nn}
 \left[
 \frac{\lambda_j{\bf 1}_{B_j}}
      {\|{\bf 1}_{B_j}\|_{X(\bx)}}
 \right]^{s_0}
 \right\}^{1/s_0}
 \right\|_{X(\bx)}
 <\infty,
\]
then \(f\in H_I^X(\bx)\) and
\[
 \|f\|_{H_I^X(\bx)}
 \leq
 C
 \left\|
 \left\{
 \sum_{j\in\nn}
 \left[
 \frac{\lambda_j{\bf 1}_{B_j}}
      {\|{\bf 1}_{B_j}\|_{X(\bx)}}
 \right]^{s_0}
 \right\}^{1/s_0}
 \right\|_{X(\bx)}.
\]
Consequently,
\[
 H_{I,\mathrm{mol},M,\varepsilon}^X(\bx)\subset H_I^X(\bx)
\]
algebraically and topologically.
\end{proposition}

\begin{proof}
Notice that, although the semigroup \(\{e^{-tI}\}_{t>0}\) does not satisfy the reinforced off-diagonal estimates, it does satisfy the Davies--Gaffney estimates. Consequently, the argument developed in \cite{lyyy23} can be adapted to the present setting. For completeness, we provide the details below.

To prove this, we will use Lemma \ref{l2.8} in its sublinear form. Its proof only requires
the \(L^2\)-boundedness and sublinearity of the operator under consideration.
Thus it is enough to prove that there exists
$
 \eta>d(\frac1{s_0}-\frac1{q_0})
$
such that, for every \((X,2,M,\varepsilon)_I\)-molecule \(m\) associated with
a ball \(B=B(x_B,r_B)\) and every \(i\in\zz_+\),
\begin{equation}\label{e:SI-molecule-annular}
 \lf\|S_I(m)\r\|_{L^2(S_i(B))}
 \lesssim
 2^{-i\eta}
 \mu(2^iB)^{1/2}
 \lf\|{\bf 1}_{2^iB}\r\|_{X(\bx)}^{-1}.
\end{equation}

Set
\[
 A_i(B):=
 \mu(2^iB)^{1/2}
 \lf\|{\bf 1}_{2^iB}\r\|_{X(\bx)}^{-1},
 \qquad i\in\zz_+.
\]
Consequently, since \(q_0\leq2\), by \eqref{e:ind-comp-lem39-short}, we have
\begin{equation}\label{e:Ai-A0-upper}
 A_i(B)\lesssim A_0(B),\qquad i\in\zz_+,
\end{equation}
and, by the doubling condition,
\begin{equation}\label{e:A0-Ai-upper}
 A_0(B)\lesssim
 2^{id(1/s_0-1/2)}A_i(B),
 \qquad i\in\zz_+.
\end{equation}

We first note that \(S_I\) is bounded on \(L^2(\bx)\). Indeed, by Fubini's
theorem and the volume comparison \(V_t(x)\sim V_t(y)\) whenever
\(d(x,y)<t\),
\[
\begin{aligned}
 \lf\|S_I(g)\r\|_{L^2(\bx)}^2
 =
 \int_\bx\int_0^\infty\int_{B(x,t)}
 \lf|t^2e^{-t^2}g(y)\r|^2
 \frac{d\mu(y)\,dt}{V_t(x)t}\,d\mu(x)
 \lesssim
 \int_0^\infty t^3e^{-2t^2}\,dt\,\|g\|_{L^2(\bx)}^2
 \lesssim
 \|g\|_{L^2(\bx)}^2.
\end{aligned}
\]
Moreover, from the molecular size condition and \eqref{e:Ai-A0-upper}, for
each \(k\in\{0,1,\ldots,M\}\),
\[
\begin{aligned}
 \|m\|_{L^2(\bx)}^2
 =
 \sum_{\ell=0}^\infty \|m\|_{L^2(S_\ell(B))}^2
\lesssim
 r_B^{4k}
 \sum_{\ell=0}^\infty
 2^{-2\ell\varepsilon}A_\ell(B)^2
 \lesssim
 r_B^{4k}A_0(B)^2.
\end{aligned}
\]
Hence
\begin{equation}\label{e:molecule-L2-I}
 \|m\|_{L^2(\bx)}\lesssim r_B^{2k}A_0(B),
 \qquad k\in\{0,1,\ldots,M\}.
\end{equation}

For \(i=0,1\), the \(L^2\)-boundedness of \(S_I\), \eqref{e:molecule-L2-I}
with \(k=0\), and the finite-dilation comparison of \(A_i(B)\) imply
\[
 \|S_I(m)\|_{L^2(S_i(B))}
 \leq
 \|S_I(m)\|_{L^2(\bx)}
 \lesssim
 A_0(B)
 \lesssim
 A_i(B).
\]
Thus \eqref{e:SI-molecule-annular} holds for \(i=0,1\).

Let now \(i\geq2\). We split
\[
 \|S_I(m)\|_{L^2(S_i(B))}^2=\mathrm{I}_1+\mathrm{I}_2,
\]
where
\[
 \mathrm{I}_1:=
 \int_{S_i(B)}
 \int_0^{2^{i-2}r_B}\int_{B(x,t)}
 \lf|t^2e^{-t^2}m(y)\r|^2
 \frac{d\mu(y)\,dt}{V_t(x)t}\,d\mu(x)
\]
and
\[
\mathrm{ I}_2:=
 \int_{S_i(B)}
 \int_{2^{i-2}r_B}^{\infty}\int_{B(x,t)}
 \lf|t^2e^{-t^2}m(y)\r|^2
 \frac{d\mu(y)\,dt}{V_t(x)t}\,d\mu(x).
\]

For \(\mathrm{I}_1\), if \(x\in S_i(B)\), \(0<t<2^{i-2}r_B\), and \(y\in B(x,t)\),
then
$
 y\in 2^{i+1}B\setminus 2^{i-2}B.
$
Using Fubini's theorem and \(V_t(x)\sim V_t(y)\) when \(d(x,y)<t\), we get
\[
\begin{aligned}
 \mathrm{I}_1
 &\lesssim
 \int_0^{2^{i-2}r_B}
 t^3e^{-2t^2}
 \int_{2^{i+1}B\setminus 2^{i-2}B}|m(y)|^2\,d\mu(y)\,dt
 \lesssim
 \sum_{\ell=-1}^1
 \|m\|_{L^2(S_{i+\ell}(B))}^2 \\
 &\lesssim
 \sum_{\ell=-1}^1
 2^{-2(i+\ell)\varepsilon}A_{i+\ell}(B)^2
 \lesssim
 2^{-2i\varepsilon}A_i(B)^2.
\end{aligned}
\]
Thus
\begin{equation}\label{e:I1-estimate}
 (\mathrm{I}_1)^{1/2}
 \lesssim
 2^{-i\varepsilon}A_i(B).
\end{equation}

For the term \(\mathrm{I}_2\), using again the Fubini theorem and the volume comparison, we have
\[
 \mathrm{I}_2
 \lesssim
 \|m\|_{L^2(\bx)}^2
 \int_{2^{i-2}r_B}^{\infty}t^3e^{-2t^2}\,dt.
\]
Taking \(k=M\) in \eqref{e:molecule-L2-I}, and using
\[
 \int_R^\infty t^3e^{-2t^2}\,dt
 \lesssim
 (1+R^2)e^{-cR^2},
 \qquad R>0,
\]
we obtain
\[
\begin{aligned}
 \mathrm{I}_2
 \lesssim
 r_B^{4M}A_0(B)^2
 (1+(2^ir_B)^2)e^{-c(2^ir_B)^2}
 \lesssim
 2^{-4Mi}A_0(B)^2.
\end{aligned}
\]
By \eqref{e:A0-Ai-upper}, this yields
\[
 (\mathrm{I}_2)^{1/2}
 \lesssim
 2^{-2Mi}A_0(B)
 \lesssim
 2^{-i[2M-d(1/s_0-1/2)]}A_i(B).
\]
Equivalently,
\begin{equation}\label{e:I2-estimate}
 (\mathrm{I}_2)^{1/2}
 \lesssim
 2^{-i[2M+d(1/2-1/s_0)]}A_i(B).
\end{equation}
Combining \eqref{e:I1-estimate} and \eqref{e:I2-estimate}, we obtain
\[
 \|S_I(m)\|_{L^2(S_i(B))}
 \lesssim
 \left[
 2^{-i\varepsilon}
 +
 2^{-i[2M+d(1/2-1/s_0)]}
 \right]A_i(B).
\]
By the assumptions on \(\varepsilon\) and \(M\), choose
\[
 \eta\in
 \left(
 d\left[\frac1{s_0}-\frac1{q_0}\right],
 \min\left\{
 \varepsilon,\,
 2M+d\left[\frac12-\frac1{s_0}\right]
 \right\}
 \right).
\]
Then \eqref{e:SI-molecule-annular} follows.

Now let
\[
 f=\sum_{j\in\nn}\lambda_jm_j
\ \ \hbox{in }L^2(\bx)
\]
be as in the statement. Applying Lemma \ref{l2.8} to the sublinear
\(L^2\)-bounded operator \(S_I\), with the annular estimate
\eqref{e:SI-molecule-annular}, gives
\[
 \|S_I(f)\|_{X(\bx)}
 \lesssim
 \left\|
 \left\{
 \sum_{j\in\nn}
 \left[
 \frac{\lambda_j{\bf 1}_{B_j}}
      {\lf\|{\bf 1}_{B_j}\r\|_{X(\bx)}}
 \right]^{s_0}
 \right\}^{\frac{1}{s_0}}
 \right\|_{X(\bx)}.
\]
Thus \(f\in H_I^X(\bx)\), and
\[
 \|f\|_{H_I^X(\bx)}
 =
 \|S_I(f)\|_{X(\bx)}
 \lesssim
 \left\|
 \left\{
 \sum_{j\in\nn}
 \left[
 \frac{\lambda_j{\bf 1}_{B_j}}
      {\lf\|{\bf 1}_{B_j}\r\|_{X(\bx)}}
 \right]^{s_0}
 \right\}^{\frac{1}{s_0}}
 \right\|_{X(\bx)}.
\]
Taking the infimum over all molecular representations of \(f\) and passing to
the completions, we obtain
\[
 H_{I,\mathrm{mol},M,\varepsilon}^X(\bx)\subset H_I^X(\bx)
\]
algebraically and topologically.
\end{proof}

\begin{corollary}
Assume that \(X(\bx)\) satisfies Assumptions \ref{a11} and \ref{a12}
with \(s_0\in(0,1]\), \(q_0\in(s_0,2]\), and \(s_0<p\).
Let \(M\in\nn\) and \(\varepsilon>0\) satisfy
\[
\varepsilon>d\left(\frac1{s_0}-\frac1{q_0}\right)
\quad\text{and}\quad
M>\frac d2\left(\frac2{s_0}-\frac12-\frac1{q_0}\right).
\]
Then
\[
H_I^X(\bx)
=
H_{I,\mathrm{at},M}^X(\bx)
=
H_{I,\mathrm{mol},M,\varepsilon}^X(\bx)
\]
with equivalent quasi-norms.
\end{corollary}

\section{The proofs of Theorems \ref{mol}, \ref{eq0} and \ref{eq1}}\label{s5}
In this section, we will devote to the proofs of Theorems \ref{mol}, \ref{eq0} and \ref{eq1}.

Now, we prove Theorem \ref{mol} via applying Lemma \ref{l2.8} and Assumption \ref{a2.5}.
\begin{proof}[Proof of Theorem \ref{mol} ]
Firstly, we show ${\rm(i)}$.
Let
$
 \varepsilon>d(\frac1{s_0}-\frac1{q_0})
$
and \(M\in\nn\) such that
\[
 2M>d\left(\frac2{s_0}-\frac12-\frac1{q_0}\right).
\]
Take
\[
 \eta\in
 \left(
 d\left[\frac1{s_0}-\frac1{q_0}\right],
 \min\left\{\varepsilon,\,
 2M-d\left[\frac1{s_0}-\frac12\right]\right\}
 \right).
\]
Suppose that
\[
 f\in L^2(\bx)\cap h^X_{L,\mathrm{mol},M,\varepsilon}(\bx).
\]
We can write
\[
 f=\sum_{j\in\nn}\lambda_jm_j \ \ \mathrm{in} \ \ L^2(\bx),
\]
where, for every \(j\in\nn\), \(\lambda_j>0\) and
\(m_j\) is a local \((X,2,M,\varepsilon)_L\)-molecule associated with the
ball \(B_j\), satisfying that
\[
\left\|
\left\{
\sum_{j\in\nn}
\left[
\frac{\lambda_j{\bf 1}_{B_j}}
     {\|{\bf 1}_{B_j}\|_{X(\bx)}}
\right]^{s_0}
\right\}^{1/s_0}
\right\|_{X(\bx)}
<\infty.
\]
In order to prove \(f\in h_L^X(\bx)\), we have to show that
\(S_L^{\rm loc}(f)\) and \(S_I(e^{-L}f)\) belong to \(X(\bx)\).

Notice that the local square function \(S_L^{\rm loc}\), \(S_I\), and
\(e^{-L}\) are bounded from \(L^2(\bx)\) into itself. Therefore, according to
Lemma \ref{l2.8}, the proof will be finished when we show that there exist
\(C>0\) and
\[
 \eta>d\left(\frac1{s_0}-\frac1{q_0}\right)
\]
such that
\begin{equation}\label{objective}
\begin{aligned}
\lf\|S_L^{\rm loc}(m_j)\r\|_{L^2(S_i(B_j))}
+
\lf\|S_I(e^{-L}m_j)\r\|_{L^2(S_i(B_j))}
\lesssim2^{-i\eta}
\mu(2^iB_j)^{1/2}
\lf\|{\bf 1}_{2^iB_j}\r\|_{X(\bx)}^{-1},
\end{aligned}
\end{equation}
for each \(i\in\zz_+\) and \(j\in\nn\).

Let \(m\) be a local \((X,2,M,\varepsilon)_L\)-molecule associated to the
ball \(B=B(x_B,r_B)\), with \(x_B\in\bx\) and \(r_B>0\). By proceeding as in
the standard annular summation and using Assumption \ref{a12}, we can see
that
\begin{equation}\label{m}
 \|m\|_{L^2(\bx)}
 \lesssim\mu(B)^{1/2}\|{\bf 1}_B\|_{X(\bx)}^{-1}.
\end{equation}
Indeed, if we set
\[
 A_i(B):=
 \mu(2^iB)^{1/2}
 \lf\|{\bf 1}_{2^iB}\r\|_{X(\bx)}^{-1},
\]
then, by Assumption \ref{a12} and \(q_0\leq2\), we obtain
\[
 A_i(B)\lesssim A_0(B),\qquad i\in\zz_+.
\]
Moreover, by Assumption \ref{a12}, we obtain, for all \(i,\nu\in\zz_+\),
\begin{equation}\label{Acomparison}
 A_\nu(B)
 \lesssim
 \begin{cases}
 2^{(i-\nu)d/s_0}A_i(B), & \nu<i,\\[3pt]
 2^{(\nu-i)d(1/2-1/q_0)}A_i(B)\lesssim A_i(B), & \nu\geq i.
 \end{cases}
\end{equation}
Therefore, the molecule condition gives \eqref{m}. Then, by the
\(L^2\)-boundedness of \(S_L^{\rm loc}\), \(S_I\), and \(e^{-L}\), we can
write, for every \(i\in\zz_+\),
\begin{equation}\label{gen}
\begin{aligned}
\lf\|S_L^{\rm loc}(m)\r\|_{L^2(S_i(B))}
+
\lf\|S_I(e^{-L}m)\r\|_{L^2(S_i(B))}
\lesssim\|m\|_{L^2(\bx)}
\lesssim\mu(B)^{1/2}\|{\bf 1}_B\|_{X(\bx)}^{-1}.
\end{aligned}
\end{equation}

From Assumption \ref{a2.5}, we know that, the kernel of integral operator $(t^2L)^ke^{-at^2L}$ and $L^ke^{-L}$ satisfy the following Gaussian upper bound estimates:
For every \(N>0\), \(a>0\), \(k\in\zz_+\), \(t>0\),
$$
\left|K_{k,a,t}(x,y)\right|
\leq
\frac{C_{a,k}}{\mu(B(x,t))}
\exp\left(
-c_{a,k}\frac{d(x,y)^2}{t^2}
\right).
$$
and
$$\left|K_k(x,y)\right|
\leq
\frac{C_k}{\mu(B(x,1))}
\exp\left(-c_k d(x,y)^2\right),$$
respectively.
This further implies that, for every \(N>0\), \(a>0\), \(k\in\zz_+\), \(t>0\),
and all closed sets \(E,F\subset\bx\),
\begin{equation}\label{OD1}
\left\|
{\bf 1}_E(t^2L)^ke^{-at^2L}{\bf 1}_F
\right\|_{L^2(\bx)\to L^2(\bx)}
\leq
C_N
\left(1+\frac{\dist(E,F)}{t}\right)^{-N}
\end{equation}
and
\begin{equation}\label{OD2}
\left\|
{\bf 1}_EL^ke^{-L}{\bf 1}_F
\right\|_{L^2(\bx)\to L^2(\bx)}
\leq
C_N\frac{1}{(1+\dist(E,F))^{N}}.
\end{equation}
Also, for every measurable set \(E\subset\bx\), if
\[
 E^{[t]}:=\{y\in\bx:\dist(y,E)<t\},
\]
then the doubling condition and Fubini's theorem imply that
\begin{equation}\label{cone1}
\lf\|S_L^{\rm loc}(g)\r\|_{L^2(E)}^2
\lesssim\int_0^1
\lf\|{\bf 1}_{E^{[t]}}t^2Le^{-t^2L}g\r\|_{L^2(\bx)}^2
\frac{dt}{t},
\end{equation}
and
\begin{equation}\label{cone2}
\lf\|S_I(e^{-L}g)\r\|_{L^2(E)}^2
\lesssim\int_0^\infty
t^4e^{-2t^2}
\lf\|{\bf 1}_{E^{[t]}}e^{-L}g\r\|_{L^2(\bx)}^2
\frac{dt}{t}.
\end{equation}

Now we assume that \(r_B\geq1\), that is, \(m\) is a local
\((X,2,M,\varepsilon)_L\)-molecule of \(({\rm I})\)-type. Let
\(i\in\nn\), \(i\geq3\). We decompose
\[
 m=m_1+m_2,
 \qquad
 m_1=m{\bf 1}_{2^{i+2}B\setminus2^{i-3}B}.
\]
We can write
\[
\begin{aligned}
\lf\|S_L^{\rm loc}(m)\r\|_{L^2(S_i(B))}^2
\leq
C\lf[
\lf\|S_L^{\rm loc}(m_1)\r\|_{L^2(S_i(B))}^2
+
\lf\|S_L^{\rm loc}(m_2)\r\|_{L^2(S_i(B))}^2
\r]
=:\mathrm{J}_1+\mathrm{J}_2.
\end{aligned}
\]
Since \(S_L^{\rm loc}\) is bounded from \(L^2(\bx)\) into itself, it follows
that
\begin{align}\label{ref1}
\mathrm{J}_1
\lesssim\|m_1\|_{L^2(\bx)}^2
\lesssim\sum_{\ell=-2}^{2}
\|m\|_{L^2(S_{i+\ell}(B))}^2
&\lesssim2^{-2i\varepsilon}
\mu(2^iB)
\lf\|{\bf 1}_{2^iB}\r\|_{X(\bx)}^{-2}.
\end{align}
On the other hand, since \(r_B\geq1\), for \(x\in S_i(B)\), \(t\in(0,1)\),
\(y\in S_i(B)^{[t]}\), and
\(z\in(2^{i+2}B\setminus2^{i-3}B)^\complement\), one has
\[
 d(y,z)\gtrsim 2^ir_B.
\]
Then, by \eqref{cone1} and \eqref{OD1}, choosing \(N>0\) large enough, we
obtain
\begin{align}\label{ref2}
\mathrm{J}_2
\lesssim\int_0^1
\left(\frac{t}{2^ir_B}\right)^{2N}
\|m\|_{L^2(\bx)}^2
\frac{dt}{t}
\lesssim(2^ir_B)^{-2N}
\mu(B)\|{\bf 1}_B\|_{X(\bx)}^{-2}.
\end{align}
By Assumption \ref{a12}, we obtain
\[
 \mu(B)^{1/2}\|{\bf 1}_B\|_{X(\bx)}^{-1}
 \lesssim2^{id(1/s_0-1/2)}
 \mu(2^iB)^{1/2}
 \lf\|{\bf 1}_{2^iB}\r\|_{X(\bx)}^{-1}.
\]
Since \(r_B\geq1\), by choosing \(N\) large enough, \eqref{ref2} gives
\[
\mathrm{J}_2
\lesssim2^{-2i\varepsilon}
\mu(2^iB)
\lf\|{\bf 1}_{2^iB}\r\|_{X(\bx)}^{-2}.
\]
By combining the above estimates, we obtain
\begin{equation}\label{acot1}
\lf\|S_L^{\rm loc}(m)\r\|_{L^2(S_i(B))}
\lesssim2^{-i\varepsilon}
\mu(2^iB)^{1/2}
\lf\|{\bf 1}_{2^iB}\r\|_{X(\bx)}^{-1},
\quad i\in\nn,\ i\geq3.
\end{equation}

Now we deal with \(S_I(e^{-L}m)\). By \eqref{cone2} and the above
decomposition of \(m\), we write
\[
\lf\|S_I(e^{-L}m)\r\|_{L^2(S_i(B))}^2
\leq
C(\mathrm{K}_1+\mathrm{K}_2),
\]
where
\[
\mathrm{K}_1:=
\int_0^\infty
t^4e^{-2t^2}
\lf\|e^{-L}m_1\r\|_{L^2(\bx)}^2
\frac{dt}{t}
\]
and
\[
\mathrm{K}_2:=
\int_0^\infty
t^4e^{-2t^2}
\lf\|{\bf 1}_{S_i(B)^{[t]}}e^{-L}m_2\r\|_{L^2(\bx)}^2
\frac{dt}{t}.
\]
The \(L^2\)-boundedness of \(e^{-L}\) and the molecule condition imply
\[
\mathrm{K}_1
\lesssim\|m_1\|_{L^2(\bx)}^2
\lesssim2^{-2i\varepsilon}
\mu(2^iB)
\lf\|{\bf 1}_{2^iB}\r\|_{X(\bx)}^{-2}.
\]
To estimate \(\mathrm{K}_2\), we split the integral into
\[
 (0,c2^ir_B)\quad\hbox{and}\quad [c2^ir_B,\infty),
\]
where \(c>0\) is sufficiently small. For \(t<c2^ir_B\), the sets
\(S_i(B)^{[t]}\) and \(\supp m_2\) are separated by a distance comparable to
\(2^ir_B\). Hence, by \eqref{OD2}, we have
\[
\int_0^{c2^ir_B}
t^4e^{-2t^2}
\lf\|{\bf 1}_{S_i(B)^{[t]}}e^{-L}m_2\r\|_{L^2(\bx)}^2
\frac{dt}{t}
\leq
C_N(2^ir_B)^{-2N}\|m\|_{L^2(\bx)}^2.
\]
For \(t\geq c2^ir_B\), we use only the \(L^2\)-boundedness of \(e^{-L}\)
and the exponential factor \(e^{-2t^2}\), obtaining
\[
\int_{c2^ir_B}^{\infty}
t^4e^{-2t^2}
\lf\|{\bf 1}_{S_i(B)^{[t]}}e^{-L}m_2\r\|_{L^2(\bx)}^2
\frac{dt}{t}
\lesssim e^{-c(2^ir_B)^2}\|m\|_{L^2(\bx)}^2.
\]
Using \eqref{m}, Assumption \ref{a12}, and \(r_B\geq1\), we take \(N\)
sufficiently large to get
\[
\mathrm{K}_2
\lesssim2^{-2i\varepsilon}
\mu(2^iB)
\lf\|{\bf 1}_{2^iB}\r\|_{X(\bx)}^{-2}.
\]
Thus
\begin{equation}\label{acot2}
\lf\|S_I(e^{-L}m)\r\|_{L^2(S_i(B))}
\lesssim2^{-i\varepsilon}
\mu(2^iB)^{1/2}
\lf\|{\bf 1}_{2^iB}\r\|_{X(\bx)}^{-1},
\quad i\in\nn,\ i\geq3.
\end{equation}
It is clear also, from \eqref{gen}, that, for \(i=0,1,2\),
\[
\lf\|S_L^{\rm loc}(m)\r\|_{L^2(S_i(B))}
+
\lf\|S_I(e^{-L}m)\r\|_{L^2(S_i(B))}
\lesssim2^{-i\eta}
\mu(2^iB)^{1/2}
\lf\|{\bf 1}_{2^iB}\r\|_{X(\bx)}^{-1}.
\]
This, jointly with \eqref{acot1} and \eqref{acot2}, leads to the estimate
\eqref{objective} for molecules of \(({\rm I})\)-type.

Assume next that \(m\) is a molecule of \(({\rm II})\)-type. Then
\(r_B\in(0,1)\) and there exists \(b\in L^2(\bx)\) such that \(m=L^Mb\) and,
for every \(k\in\{0,1,\ldots,M\}\),
\[
\lf\|L^kb\r\|_{L^2(S_i(B))}
\leq
2^{-i\varepsilon}r_B^{2(M-k)}
\mu(2^iB)^{1/2}
\lf\|{\bf 1}_{2^iB}\r\|_{X(\bx)}^{-1},
\qquad i\in\zz_+.
\]
By proceeding as in the proof of \eqref{m}, we can see that
\begin{equation}\label{B6}
\|b\|_{L^2(\bx)}
\lesssim r_B^{2M}\mu(B)^{1/2}
\|{\bf 1}_B\|_{X(\bx)}^{-1}.
\end{equation}
Put
\[
 \alpha:=d\left(\frac1{s_0}-\frac1{q_0}\right)
 \quad\text{and}\quad
 \beta:=d\left(\frac1{s_0}-\frac12\right).
\]
Since
\[
 2M>
 d\left(\frac2{s_0}-\frac12-\frac1{q_0}\right)
 =
 \alpha+\beta
\]
and \(\varepsilon>\alpha\), we may fix
\begin{equation}\label{choice-eta-II}
 \eta\in
 \left(
 \alpha,\,
 \min\{\varepsilon,2M-\beta\}
 \right).
\end{equation}
Moreover, by Assumption \(\ref{a12}\), for all
\(i,\nu\in\mathbb Z_+\), we obtain
\begin{equation}\label{Acomparison-sharp}
 A_\nu(B)
 \lesssim
 \begin{cases}
 2^{(i-\nu)\beta}A_i(B),&\nu<i,\\[3pt]
 A_i(B),&\nu\geq i.
 \end{cases}
\end{equation}
In particular,
\begin{equation}\label{A0-Ai-II}
 A_0(B)\lesssim 2^{i\beta}A_i(B).
\end{equation}

Now we prove the desired estimate for
\(S_L^{\rm loc}(m)\). The estimate for a fixed finite number of indices
\(i\) follows from \(\eqref{gen}\), \(\eqref{Acomparison-sharp}\), and
the finite-dilation property. Thus, it suffices to consider all
sufficiently large \(i\). Put
\[
 E_i:=S_i(B), \quad R_i:=2^ir_B, \quad
\mathrm{and} \quad
 \Omega_i:=2^{i+2}B\setminus 2^{i-3}B.
\]
We may choose \(i_0\in\mathbb N\), depending only on the geometric
constants of \((\mathbb X,d,\mu)\), such that, for every \(i\geq i_0\),
there exists a constant \(c\in(0,1)\), independent of \(i\) and \(B\),
satisfying
\begin{equation}\label{separation-II}
 \operatorname{dist}
 \left(E_i^{[t]},\Omega_i^{\complement}\right)
 \gtrsim R_i,
 \qquad 0<t<cR_i,
\end{equation}
and \(cR_i\geq r_B\).

By \(\eqref{cone1}\), the identity \(m=L^Mb\), and
\[
 t^2Le^{-t^2L}L^Mb
 =
 t^{-2M}(t^2L)^{M+1}e^{-t^2L}b,
\]
we have
\begin{equation}\label{II-time-decomposition}
\begin{aligned}
 \left\|S_L^{\rm loc}(m)\right\|_{L^2(S_i(B))}^2
 \lesssim{}&
 \int_0^{r_B}
 \left\|
 {\bf1}_{E_i^{[t]}}
 t^2Le^{-t^2L}m
 \right\|_{L^2(\bx)}^2\frac{dt}{t}
+
 \int_{r_B}^1
 t^{-4M}
 \left\|
 {\bf1}_{E_i^{[t]}}
 (t^2L)^{M+1}e^{-t^2L}b
 \right\|_{L^2(\bx)}^2\frac{dt}{t}
\\
 =:{}&\mathcal J_<+\mathcal J_>.
\end{aligned}
\end{equation}

We first estimate \(\mathcal J_<\). Decompose
\[
 m=m_1+m_2,
 \qquad
 m_1:=m{\bf1}_{\Omega_i}.
\]
The \(L^2\)-boundedness of the vertical square function associated with
\(L\), together with the molecular size condition, gives
\begin{align}
 \mathcal J_{<,1}
 \lesssim
 \left\|m_1\right\|_{L^2(\bx)}^2
 \lesssim
 \sum_{\nu=i-2}^{i+2}
 2^{-2\nu\varepsilon}A_\nu(B)^2
 \lesssim
 2^{-2i\varepsilon}A_i(B)^2.
 \label{II-short-near}
\end{align}
For $m_2$, \(\eqref{separation-II}\) and
\(\eqref{OD1}\) imply that, for every \(N>0\),
\begin{align}
 \mathcal J_{<,2}
 \lesssim
 \int_0^{r_B}
 \left(\frac{t}{R_i}\right)^{2N}
 \left\|m_2\right\|_{L^2(\mathbb X)}^2
 \frac{dt}{t}
 \lesssim
 2^{-2iN}\|m\|_{L^2(\bx)}^2
 \lesssim
 2^{-2iN}A_0(B)^2
 \lesssim
 2^{-2i(N-\beta)}A_i(B)^2.
 \label{II-short-far}
\end{align}
Here we used \(\eqref{m}\) and \(\eqref{A0-Ai-II}\).

We next estimate \(\mathcal J_>\). Write
\[
 b=b_1+b_2,
 \qquad
 b_1:=b{\bf1}_{\Omega_i}.
\]
By the spectral theorem, we obtain
\begin{equation}\label{vertical-QM}
 \int_0^\infty
 \left\|
 (t^2L)^{M+1}e^{-t^2L}h
 \right\|_{L^2(\mathbb X)}^2\frac{dt}{t}
 \lesssim
 \|h\|_{L^2(\mathbb X)}^2,
 \qquad h\in L^2(\mathbb X).
\end{equation}
Consequently,
\begin{align}
 \mathcal J_{>,1}
 &\leq
 r_B^{-4M}
 \int_0^\infty
 \left\|
 (t^2L)^{M+1}e^{-t^2L}b_1
 \right\|_{L^2(\mathbb X)}^2\frac{dt}{t} \notag
 \lesssim
 r_B^{-4M}\left\|b_1\right\|_{L^2(\mathbb X)}^2 \notag\\
 &\lesssim
 r_B^{-4M}
 \sum_{\nu=i-2}^{i+2}
 2^{-2\nu\varepsilon}r_B^{4M}A_\nu(B)^2
 \lesssim
 2^{-2i\varepsilon}A_i(B)^2.
 \label{II-long-near}
\end{align}

For the term $\mathcal J_{>,2}$, fix \(N>2M\). By
\(\eqref{separation-II}\), \(\eqref{OD1}\), and the
\(L^2\)-boundedness of
\((t^2L)^{M+1}e^{-t^2L}\), we obtain
\begin{align}
 \mathcal J_{>,2}
 \lesssim{}&
 \left[
 \int_{r_B}^{\min\{1,cR_i\}}
 t^{-4M}
 \left(\frac{t}{R_i}\right)^{2N}
 \frac{dt}{t}
 +
 \int_{\max\{r_B,cR_i\}}^1
 t^{-4M}\frac{dt}{t}
 \right]
 \|b_2\|_{L^2(\mathbb X)}^2
 \lesssim
 R_i^{-4M}\|b\|_{L^2(\mathbb X)}^2.
 \label{II-long-far-first}
\end{align}
The integrals over empty intervals are understood to be zero. Indeed,
since \(N>2M\),
\[
 R_i^{-2N}
 \int_{r_B}^{\min\{1,cR_i\}}
 t^{2N-4M}\frac{dt}{t}
 \lesssim R_i^{-4M} \quad
\mathrm{and} \quad
 \int_{\max\{r_B,cR_i\}}^1t^{-4M}\frac{dt}{t}
 \lesssim R_i^{-4M}.
\]
Using \(\eqref{B6}\) and \(\eqref{A0-Ai-II}\), we further obtain
\begin{align}
 \mathcal J_{>,2}
 \lesssim
 (2^ir_B)^{-4M}r_B^{4M}A_0(B)^2
 \lesssim
 2^{-2i(2M-\beta)}A_i(B)^2.
 \label{II-long-far}
\end{align}

Combining
\(\eqref{II-short-near}\)--\(\eqref{II-long-far}\), and choosing
\(N>\beta+\eta\), we conclude that
\[
\begin{aligned}
 \left\|S_L^{\rm loc}(m)\right\|_{L^2(S_i(B))}
 \lesssim
 \left[
 2^{-i\varepsilon}
 +
 2^{-i(N-\beta)}
 +
 2^{-i(2M-\beta)}
 \right]A_i(B)
 \lesssim
 2^{-i\eta}A_i(B),
\end{aligned}
\]
where the last inequality follows from
\(\eqref{choice-eta-II}\). Therefore,
\begin{equation}\label{acot3}
 \left\|S_L^{\rm loc}(m)\right\|_{L^2(S_i(B))}
 \lesssim
 2^{-i\eta}
 \mu(2^iB)^{1/2}
 \left\|{\bf1}_{2^iB}\right\|_{X(\mathbb X)}^{-1},
 \qquad i\in\mathbb Z_+.
\end{equation}


Let us now show that this estimate is also satisfied by \(S_I(e^{-L}m)\).
We observe that \(L^Me^{-L}\) is a bounded operator from \(L^2(\bx)\) into
itself and, by \eqref{OD2}, satisfies the off-diagonal estimate
\[
\left\|
{\bf 1}_EL^Me^{-L}{\bf 1}_F
\right\|_{L^2(\bx)\to L^2(\bx)}
\leq
C_N(1+\dist(E,F))^{-N}.
\]
Let \(i\in\nn\), \(i\geq3\). If \(2^ir_B\leq C_0\), then, by the
\(L^2\)-boundedness of \(S_I\), \(L^Me^{-L}\), and \eqref{B6}, we obtain
\begin{align*}
\lf\|S_I(e^{-L}m)\r\|_{L^2(S_i(B))}
&=\lf\|S_I(L^Me^{-L}b)\r\|_{L^2(S_i(B))}
\lesssim\|b\|_{L^2(\bx)}       \\
&\lesssim2^{-i[2M-d(1/s_0-1/2)]}
\mu(2^iB)^{1/2}
\lf\|{\bf 1}_{2^iB}\r\|_{X(\bx)}^{-1}  \\
&\lesssim2^{-i\eta}
\mu(2^iB)^{1/2}
\lf\|{\bf 1}_{2^iB}\r\|_{X(\bx)}^{-1}.
\end{align*}
If \(2^ir_B>C_0\), we decompose \(b=b_1+b_2\), where
\[
 b_1=b{\bf 1}_{2^{i+3}B\setminus2^{i-3}B}.
\]
By the \(L^2\)-boundedness of \(S_I\) and \(L^Me^{-L}\),
\[
\lf\|S_I(L^Me^{-L}b_1)\r\|_{L^2(S_i(B))}
\leq
C\|b_1\|_{L^2(\bx)}
\leq
C2^{-i\varepsilon}r_B^{2M}
\mu(2^iB)^{1/2}
\lf\|{\bf 1}_{2^iB}\r\|_{X(\bx)}^{-1}.
\]
Since \(r_B<1\), this is bounded by the required right-hand side. For the
term involving \(b_2\), we use \eqref{cone2}, the above off-diagonal estimate
for \(L^Me^{-L}\), and the exponential decay of \(e^{-t^2}\). Splitting the
integral in \(t\) into
\[
 (0,c2^ir_B)\quad\hbox{and}\quad [c2^ir_B,\infty),
\]
and using \eqref{B6}, we obtain
\[
\lf\|S_I(L^Me^{-L}b_2)\r\|_{L^2(S_i(B))}
\lesssim2^{-i\eta}
\mu(2^iB)^{1/2}
\lf\|{\bf 1}_{2^iB}\r\|_{X(\bx)}^{-1}.
\]
Hence
\begin{equation}\label{acot4}
\lf\|S_I(e^{-L}m)\r\|_{L^2(S_i(B))}
\lesssim2^{-i\eta}
\mu(2^iB)^{1/2}
\lf\|{\bf 1}_{2^iB}\r\|_{X(\bx)}^{-1},
\quad i\in\zz_+.
\end{equation}
These estimations, jointly with \eqref{acot3}, lead to \eqref{objective} for
molecules of \(({\rm II})\)-type.

We conclude that \(f\in h_L^X(\bx)\). Moreover, by Lemma \ref{l2.8},
\[
\begin{aligned}
\|f\|_{h_L^X(\bx)}
=
\lf\|S_L^{\rm loc}(f)\r\|_{X(\bx)}
+
\lf\|S_I(e^{-L}f)\r\|_{X(\bx)}
\lesssim
\left\|
\left\{
\sum_{j\in\nn}
\left[
\frac{\lambda_j{\bf 1}_{B_j}}
     {\|{\bf 1}_{B_j}\|_{X(\bx)}}
\right]^{s_0}
\right\}^{1/s_0}
\right\|_{X(\bx)}.
\end{aligned}
\]
Taking the infimum over all molecular representations of \(f\), we obtain
\[
\|f\|_{h_L^X(\bx)}
\lesssim\|f\|_{h^X_{L,\mathrm{mol},M,\varepsilon}(\bx)}.
\]
By using a density argument, we obtain that
\[
h^X_{L,\mathrm{mol},M,\varepsilon}(\bx)\subset h_L^X(\bx)
\]
algebraically and topologically.

Next we prove \({\rm(ii)}\). Let
$
 f\in L^2(\bx)\cap h_L^X(\bx).
$
Define
\[
 F(y,t):=
 t^2Le^{-t^2L}f(y){\bf 1}_{(0,1)}(t)
\ \
\mathrm{and} \ \
 G(y,t):=
 t^2e^{-t^2}e^{-L}f(y){\bf 1}_{[1,\infty)}(t).
\]
Then
\[
 F\in \mathbf T_2^X(\bx)\cap \mathbf T_2^2(\bx),
 \qquad
 G\in T_2^X(\bx)\cap T_2^2(\bx),
\]
and
\[
 \|F\|_{\mathbf T_2^X(\bx)}
 =
 \lf\|S_L^{\mathrm{loc}}(f)\r\|_{X(\bx)},
 \qquad
 \|G\|_{T_2^X(\bx)}
 \leq
 \lf\|S_I(e^{-L}f)\r\|_{X(\bx)}.
\]
From the atomic decomposition of the local tent space, we deduce that there exist local
\((\mathbf T_2^X,\infty)\)-atoms \(A_j\), coefficients
\(\lambda_j\geq0\), and balls \(B_j\), such that
\[
 F=\sum_{j\in\nn}\lambda_jA_j
 \ \ \hbox{in } \ \ \mathbf T_2^2(\bx),
\]
and
\[
 \left\|
 \left\{
 \sum_{j\in\nn}
 \left[
 \frac{\lambda_j{\bf 1}_{B_j}}
      {\|{\bf 1}_{B_j}\|_{X(\bx)}}
 \right]^{s_0}
 \right\}^{1/s_0}
 \right\|_{X(\bx)}
 \lesssim
 \|F\|_{\mathbf T_2^X(\bx)}.
\]
Similarly, by the atomic decomposition of the global tent space, there exist
 global \((T_2^X,\infty)\)-atoms \(\widetilde A_j\), coefficients
\(\mu_j\geq0\), and balls \(\widetilde B_j\), such that
\[
 G=\sum_{j\in\nn}\mu_j\widetilde A_j
\ \ \hbox{in } \ \ T_2^2(\bx),
\]
and
\[
 \left\|
 \left\{
 \sum_{j\in\nn}
 \left[
 \frac{\mu_j{\bf 1}_{\widetilde B_j}}
      {\|{\bf 1}_{\widetilde B_j}\|_{X(\bx)}}
 \right]^{s_0}
 \right\}^{1/s_0}
 \right\|_{X(\bx)}
 \lesssim
 \|G\|_{T_2^X(\bx)}.
\]
Since \(G\) is supported in \(\bx\times[1,\infty)\), every nonzero atom
\(\widetilde A_j\) may be taken to be associated with a ball
\(\widetilde B_j\) satisfying \(r_{\widetilde B_j}\geq1\).

Let
$
 c_M:=\frac{2^{M+2}}{M!},
$
and define, for \(H\in \mathbf T_{2,c}^2(\bx)\),
\[
 \pi_{M,L}^{\mathrm{loc}}(H)(x)
 :=
 c_M\int_0^1
 \lf(t^2L\r)^Me^{-t^2L}H(\cdot,t)(x)\,\frac{dt}{t}.
\]

We first record the following fact, whose proof is included here in order to
keep the argument self-contained. The operator
\(\pi_{M,L}^{\mathrm{loc}}\) extends to a bounded operator from
\(\mathbf T_2^2(\bx)\) to \(L^2(\bx)\). Indeed, for
\(H\in\mathbf T_{2,c}^2(\bx)\) and \(g\in D(L)\), by the self-adjointness of
\(L\), the H\"older inequality, and the spectral theorem, we have
\[
\begin{aligned}
\left|
\int_{\bx}\pi_{M,L}^{\mathrm{loc}}(H)(x)g(x)\,d\mu(x)
\right|
&\lesssim
\left(
\int_0^1\int_{\bx}|H(x,t)|^2\,d\mu(x)\frac{dt}{t}
\right)^{1/2}
\left(
\int_0^1
\lf\|\lf(t^2L\r)^Me^{-t^2L}g\r\|_{L^2(\bx)}^2
\frac{dt}{t}
\right)^{1/2}                                      \\
&\lesssim
\|H\|_{\mathbf T_2^2(\bx)}\|g\|_{L^2(\bx)}.
\end{aligned}
\]
Here we used \(M\in\nn\), so that
\[
 \sup_{\lambda\geq0}
 \int_0^1
 \lf(t^2\lambda\r)^{2M}e^{-2t^2\lambda}\frac{dt}{t}
 <\infty .
\]
By duality and the density of \(D(L)\) in \(L^2(\bx)\),
\(\pi_{M,L}^{\mathrm{loc}}\) extends boundedly from
\(\mathbf T_2^2(\bx)\) into \(L^2(\bx)\).

Next we show that, if \(A\) is a \((\mathbf T_2^X,2)\)-atom associated with a
ball \(B=B(x_B,r_B)\), then \(\pi_{M,L}^{\mathrm{loc}}(A)\) is, up to a
harmless multiplicative constant, a local
\((X,2,M,\varepsilon)_L\)-molecule associated with \(B\).

Assume first that \(0<r_B<1\). Define
\[
 b_A
 :=
 c_M\int_0^{r_B}
 t^{2M}e^{-t^2L}A(\cdot,t)\,\frac{dt}{t},
\]
where the integral is understood in \(L^2(\bx)\). Then, by the closedness of
\(L\),
\[
 \pi_{M,L}^{\mathrm{loc}}(A)=L^M b_A,
\]
and, for every \(k\in\{0,1,\ldots,M\}\),
\[
 L^kb_A
 =
 c_M\int_0^{r_B}
 t^{2(M-k)}\lf(t^2L\r)^ke^{-t^2L}A(\cdot,t)\,\frac{dt}{t}
\]
in \(L^2(\bx)\). Notice that \(A(\cdot,t)=0\) for \(t\ge r_B\), since
\(\supp A\subset T(B)\cap(\bx\times(0,1))\).

Let \(i\in\zz_+\). If \(i\ge2\), then
\(\dist(S_i(B),B)\sim 2^ir_B\). Hence, by the off-diagonal estimate
\eqref{OD1}, for every \(N>0\),
\[
\begin{aligned}
 \lf\|L^kb_A\r\|_{L^2(S_i(B))}
&\lesssim
 \int_0^{r_B}
 t^{2(M-k)}
 \left(1+\frac{\dist(S_i(B),B)}{t}\right)^{-N}
 \|A(\cdot,t)\|_{L^2(\bx)}\,\frac{dt}{t}     \\
&\lesssim
 2^{-iN}r_B^{2(M-k)}
 \left(
 \int_0^{r_B}\|A(\cdot,t)\|_{L^2(\bx)}^2
 \frac{dt}{t}
 \right)^{1/2}                                \\
&\lesssim
 2^{-iN}r_B^{2(M-k)}
 \mu(B)^{1/2}\|{\bf 1}_B\|_{X(\bx)}^{-1}.
\end{aligned}
\]
Using the comparison estimate \eqref{Acomparison}, we further obtain
\[
 \lf\|L^kb_A\r\|_{L^2(S_i(B))}
 \lesssim
 2^{-i[N-d/s_0]}
 r_B^{2(M-k)}
 A_i(B).
\]
Choosing \(N>d/s_0+\varepsilon\), we get
\[
 \lf\|L^kb_A\r\|_{L^2(S_i(B))}
 \lesssim
 2^{-i\varepsilon}
 r_B^{2(M-k)}
 A_i(B),
 \qquad i\ge2.
\]
For \(i=0,1\), the same estimate follows from the \(L^2\)-boundedness of the
operators involved and the finite-dilation comparison of \(A_i(B)\). Thus,
for every \(i\in\zz_+\) and every \(k\in\{0,1,\ldots,M\}\),
\[
 \lf\|L^kb_A\r\|_{L^2(S_i(B))}
 \lesssim
 2^{-i\varepsilon}
 r_B^{2(M-k)}
 A_i(B).
\]
Therefore \(C\pi_{M,L}^{\mathrm{loc}}(A)\) is a local
\((X,2,M,\varepsilon)_L\)-molecule of type \({\rm(II)}\).

Assume next that \(r_B\ge1\). Since \(A(\cdot,t)\) is supported in \(B\) for
\(0<t<1\), by \eqref{OD1}, for \(i\ge2\) and every \(N>0\), we have
\[
\begin{aligned}
 \lf\|\pi_{M,L}^{\mathrm{loc}}(A)\r\|_{L^2(S_i(B))}
&\lesssim
 \int_0^1
 \left(1+\frac{\dist(S_i(B),B)}{t}\right)^{-N}
 \|A(\cdot,t)\|_{L^2(\bx)}\,\frac{dt}{t}       \\
&\lesssim
 2^{-iN}r_B^{-N}
 \left(
 \int_0^1 t^{2N}\frac{dt}{t}
 \right)^{1/2}
 \left(
 \int_0^1\|A(\cdot,t)\|_{L^2(\bx)}^2
 \frac{dt}{t}
 \right)^{1/2}                                 \\
&\lesssim
 2^{-iN}
 \mu(B)^{1/2}\|{\bf 1}_B\|_{X(\bx)}^{-1}.
\end{aligned}
\]
Since \(r_B\ge1\), by \eqref{Acomparison}, we have
\[
 \mu(B)^{1/2}\|{\bf 1}_B\|_{X(\bx)}^{-1}
 =
 A_0(B)
 \lesssim
 2^{id/s_0}A_i(B).
\]
Choosing \(N>d/s_0+\varepsilon\), we obtain
\[
 \lf\|\pi_{M,L}^{\mathrm{loc}}(A)\r\|_{L^2(S_i(B))}
 \lesssim
 2^{-i\varepsilon}A_i(B),
 \qquad i\ge2.
\]
For \(i=0,1\), the estimate follows from the \(L^2\)-boundedness of
\(\pi_{M,L}^{\mathrm{loc}}\), the atom condition, and the finite-dilation
comparison. Consequently,
\[
 \lf\|\pi_{M,L}^{\mathrm{loc}}(A)\r\|_{L^2(S_i(B))}
 \lesssim
 2^{-i\varepsilon}A_i(B),
 \qquad i\in\zz_+.
\]
Thus \(C\pi_{M,L}^{\mathrm{loc}}(A)\) is a local
\((X,2,M,\varepsilon)_L\)-molecule of type \({\rm(I)}\).

Therefore, since
\[
 F=\sum_{j\in\nn}\lambda_jA_j
 \ \ \hbox{in } \ \ \mathbf T_2^2(\bx),
\]
the boundedness of \(\pi_{M,L}^{\mathrm{loc}}\) from \(\mathbf T_2^2(\bx)\) to
\(L^2(\bx)\) gives
\[
 \pi_{M,L}^{\mathrm{loc}}(F)
 =
 \sum_{j\in\nn}\lambda_j\pi_{M,L}^{\mathrm{loc}}(A_j)
\ \ \hbox{in } \ \ L^2(\bx).
\]
Moreover, each \(\pi_{M,L}^{\mathrm{loc}}(A_j)\), after multiplication by a
fixed harmless constant, is a local
\((X,2,M,\varepsilon)_L\)-molecule associated with \(B_j\).

Now define
\[
 c_I^{-1}:=
 \int_1^\infty t^4e^{-2t^2}\,\frac{dt}{t},
\]
and, for a global tent atom \(A\), put
\[
 \pi_I(A):=
 c_I\int_1^\infty
 t^2e^{-t^2}A(\cdot,t)\,\frac{dt}{t}.
\]
Then, for every global atom \(A\) associated with a ball \(B\) satisfying
\(r_B\geq1\), and for every \(\ell\in\{0,1,\ldots,M\}\), we know that
$
 C L^\ell e^{-L}\pi_I(A)
$
is a local molecule of type \({\rm(I)}\) associated with \(B\). Indeed,
\[
 \lf\|\pi_I(A)\r\|_{L^2(\bx)}
 \lesssim
 \left(\int_1^\infty t^4e^{-2t^2}\,\frac{dt}{t}\right)^{1/2}
 \|A\|_{T_2^2(\bx)}
 \lesssim
 A_0(B).
\]
By \eqref{OD2}, we have
\[
\begin{aligned}
 \lf\|L^\ell e^{-L}\pi_I(A)\r\|_{L^2(S_i(B))}
 \lesssim
 [1+\dist(S_i(B),B)]^{-N}\|\pi_I(A)\|_{L^2(\bx)}
 \lesssim
 2^{-iN}A_0(B).
\end{aligned}
\]
Since \(r_B\geq1\), choosing \(N>d/s_0+\varepsilon\) gives
\[
 \lf\|L^\ell e^{-L}\pi_I(A)\r\|_{L^2(S_i(B))}
 \lesssim
 2^{-i\varepsilon}A_i(B).
\]
Hence \(CL^\ell e^{-L}\pi_I(A)\) is of type \({\rm(I)}\).

For every \(\lambda\geq0\), one has
\[
 1
 =
 c_M\int_0^1
 (t^2\lambda)^{M+1}e^{-2t^2\lambda}\,\frac{dt}{t}
 +
 e^{-2\lambda}
 \sum_{\ell=0}^M
 \frac{(2\lambda)^\ell}{\ell!}.
\]
Therefore, by the spectral theorem, we obtain
\[
 f
 =
 \pi_{M,L}^{\mathrm{loc}}(F)
 +
 \sum_{\ell=0}^M
 \frac{2^\ell}{\ell!}
 L^\ell e^{-L}\pi_I(G)
 \ \ \hbox{in } \ \ L^2(\bx).
\]
Substituting the two atomic decompositions of \(F\) and \(G\), we obtain an
\(L^2\)-convergent molecular representation of \(f\). The coefficient
functional is bounded by
\begin{align*}
 &\left\|
 \left\{
 \sum_{j\in\nn}
 \left[
 \frac{\lambda_j{\bf 1}_{B_j}}
      {\|{\bf 1}_{B_j}\|_{X(\bx)}}
 \right]^{s_0}
 \right\}^{1/s_0}
 \right\|_{X(\bx)}
+
 \sum_{\ell=0}^M
 \left\|
 \left\{
 \sum_{j\in\nn}
 \left[
 \frac{\mu_j{\bf 1}_{\widetilde B_j}}
      {\|{\bf 1}_{\widetilde B_j}\|_{X(\bx)}}
 \right]^{s_0}
 \right\}^{1/s_0}
 \right\|_{X(\bx)}                                      \\
 &\hs\lesssim
 \|F\|_{\mathbf T_2^X(\bx)} +\|G\|_{T_2^X(\bx)}
\lesssim
 \lf\|S_L^{\mathrm{loc}}(f)\r\|_{X(\bx)}
 +\lf\|S_I(e^{-L}f)\r\|_{X(\bx)}
 =\|f\|_{h_L^X(\bx)}.
\end{align*}
Hence
\[
 \|f\|_{h^X_{L,\mathrm{mol},M,\varepsilon}(\bx)}
 \lesssim
 \|f\|_{h_L^X(\bx)},
 \qquad
 f\in L^2(\bx)\cap h_L^X(\bx).
\]

Finally, passing to the completions gives
\[
 h_L^X(\bx)
 \subset
 h^X_{L,\mathrm{mol},M,\varepsilon}(\bx).
\]
This finishes the proof of Theorem \ref{mol}.
\end{proof}

Next, we prove Theorem \ref{eq0} by using Theorem \ref{mol} and Lemma \ref{l2.8}.
\begin{proof}[Proof of Theorem \ref{eq0}]
We have to prove that
\[
 h^X_L(\bx)=H^X_L(\bx).
\]
Suppose that \(f\in h^X_L(\bx)\cap L^2(\bx)\) and consider
\(M\in\nn\), such that
\[
 2M>d\left(\frac2{s_0}-\frac12-\frac1{q_0}\right)
\ \
\mathrm{and} \ \
 \varepsilon>d\left(\frac1{s_0}-\frac1{q_0}\right).
\]
Then there exists
$
 \eta>d(\frac1{s_0}-\frac1{q_0})
$
such that
\[
 \eta<\varepsilon
 \quad\hbox{and}\quad
 \eta<2M-d\left(\frac1{s_0}-\frac12\right).
\]
According to Theorem \ref{mol}, there exist, for every \(j\in\nn\),
\(\lambda_j>0\) and a local \((X,2,M,\varepsilon)\)-molecule \(m_j\)
associated with the ball \(B_j\) such that
\[
 f=\sum_{j\in\nn}\lambda_jm_j \ \  \mathrm{in} \ \ L^2(\bx) \ \ \mathrm{and} \ \ h^X_L(\bx),
\]
 and
\[
\left\|
\left\{
\sum_{j\in\nn}
\left[
\lambda_j
\frac{{\bf 1}_{B_j}}{\|{\bf 1}_{B_j}\|_{X(\bx)}}
\right]^{s_0}
\right\}^{1/s_0}
\right\|_{X(\bx)}
\sim
\|f\|_{h^X_L(\bx)}.
\]
We will show that \(f\in H^X_L(\bx)\) and
\[
\|f\|_{H^X_L(\bx)}
\leq
C
\left\|
\left\{
\sum_{j\in\nn}
\left[
\lambda_j
\frac{{\bf 1}_{B_j}}{\|{\bf 1}_{B_j}\|_{X(\bx)}}
\right]^{s_0}
\right\}^{1/s_0}
\right\|_{X(\bx)}.
\]
Notice that
\[
 S_L(f)\leq S_L^{\rm loc}(f)+S_L^\infty(f),
\]
where
\[
 S_L^\infty(f)(x)
 :=
 \left(
 \int_1^\infty
 \int_{B(x,t)}
 \lf|t^2Le^{-t^2L}(f)(y)\r|^2
 \frac{d\mu(y)dt}{V_t(x)t}
 \right)^{1/2},
 \quad x\in\bx.
\]
Since \(f\in L^2(\bx)\cap h^X_L(\bx)\),
\(S_L^{\rm loc}(f)\in X(\bx)\) and
\[
\lf\|S_L^{\rm loc}(f)\r\|_{X(\bx)}
\leq
\|f\|_{h^X_L(\bx)}
\lesssim
\left\|
\left\{
\sum_{j\in\nn}
\left[
\lambda_j
\frac{{\bf 1}_{B_j}}{\|{\bf 1}_{B_j}\|_{X(\bx)}}
\right]^{s_0}
\right\}^{1/s_0}
\right\|_{X(\bx)}.
\]
By Lemma \ref{l2.8}, in order to prove that
\(S_L^\infty(f)\in X(\bx)\) and that
\[
\lf\|S_L^\infty(f)\r\|_{X(\bx)}
\lesssim
\left\|
\left\{
\sum_{j\in\nn}
\left[
\lambda_j
\frac{{\bf 1}_{B_j}}{\|{\bf 1}_{B_j}\|_{X(\bx)}}
\right]^{s_0}
\right\}^{1/s_0}
\right\|_{X(\bx)},
\]
it is sufficient to see that there exist \(C>0\) and
$
 \eta>d(\frac1{s_0}-\frac1{q_0})
$
such that
\begin{equation}\label{C1}
\lf\|S_L^\infty(m_j)\r\|_{L^2(S_i(B_j))}
\lesssim2^{-i\eta}
\mu(2^iB_j)^{1/2}
\lf\|{\bf 1}_{2^iB_j}\r\|_{X(\bx)}^{-1},
\quad i\in\zz_+,\ j\in\nn.
\end{equation}
Assume that \(m\) is a local \((X,2,M,\varepsilon)\)-molecule associated
with the ball \(B=B(x_B,r_B)\), being \(x_B\in\bx\) and \(r_B>0\).
We need the following estimate. Since \(\inf\sigma(L)>0\), for every
\(N>0\), every \(k\in\zz_+\), every \(t>0\), and all closed sets
\(E,F\subset\bx\), one has
\begin{equation}\label{gapOD}
\left\|
{\bf 1}_E(t^2L)^ke^{-t^2L}{\bf 1}_F
\right\|_{L^2(\bx)\to L^2(\bx)}
\leq
C_N e^{-ct^2}
\left(1+\frac{\dist(E,F)}{t}\right)^{-N}.
\end{equation}
Moreover, if
\[
 E^{[t]}:=\{y\in\bx:\dist(y,E)<t\},
\]
then, by the doubling condition and Fubini's theorem, we have
\begin{equation}\label{coneestimate}
 \lf\|S_L^\infty(g)\r\|_{L^2(E)}^2
 \lesssim\int_1^\infty
 \lf\|{\bf 1}_{E^{[t]}}t^2Le^{-t^2L}g\r\|_{L^2(\bx)}^2
 \frac{dt}{t}.
\end{equation}

Suppose first that \(m\) is of \(({\rm I})\)-type. Then \(r_B\geq1\).
We observe that
\[
 \|m\|_{L^2(\bx)}
 \lesssim\mu(B)^{1/2}\|{\bf 1}_B\|_{X(\bx)}^{-1}.
\]
Consider \(i\in\nn\), \(i\geq3\). As in the proof of Theorem \ref{mol}
\({\rm(i)}\), we decompose
\[
 m=m_1+m_2,
 \qquad
 m_1=m{\bf 1}_{2^{i+2}B\setminus 2^{i-3}B}.
\]
Then, by the \(L^2\)-boundedness of \(S_L^\infty\),
\[
\begin{aligned}
 \lf\|S_L^\infty(m)\r\|_{L^2(S_i(B))}^2
 \leq C\lf(\lf\|S_L^\infty(m_1)\r\|_{L^2(S_i(B))}^2
 +\lf\|S_L^\infty(m_2)\r\|_{L^2(S_i(B))}^2\r)
 =:\mathbb I_1+\mathbb I_2.
\end{aligned}
\]
As in \eqref{ref1}, \eqref{e:ind-comp-lem39-short} and the \(L^2\)-boundedness of \(S_L^\infty\) gives
\[
\begin{aligned}
 \mathbb I_1
 \lesssim\|m_1\|_{L^2(\bx)}^2
 \lesssim\sum_{\ell=-3}^{2}
 \|m\|_{L^2(S_{i+\ell}(B))}^2
 \lesssim2^{-2i\varepsilon}
 \mu(2^iB)
 \lf\|{\bf 1}_{2^iB}\r\|_{X(\bx)}^{-2}.
\end{aligned}
\]
On the other hand, by \eqref{coneestimate} and \eqref{gapOD}, we have
\begin{align}\label{acotI2}
 \mathbb I_2
 \leq
 C\|m\|_{L^2(\bx)}^2
 \left(
 \int_1^{c2^ir_B}
 e^{-ct^2}
 \left(\frac{t}{2^ir_B}\right)^{2N}
 \frac{dt}{t}
 +
 \int_{c2^ir_B}^{\infty}
 e^{-ct^2}
 \frac{dt}{t}
 \right)
 =:\mathbb I_{2,1}+\mathbb I_{2,2}.
\end{align}
By choosing \(N>0\) large enough and using \eqref{e:ind-comp-lem39-short},
we have
\[
\begin{aligned}
 \mathbb I_{2,1}
 &\lesssim\|m\|_{L^2(\bx)}^2(2^ir_B)^{-2N}
 \lesssim\mu(B)\|{\bf 1}_B\|_{X(\bx)}^{-2}
 2^{-2iN}r_B^{-2N}                                         \\
 &\lesssim2^{-2i[N-d(1/s_0-1/2)]}
 \mu(2^iB)\lf\|{\bf 1}_{2^iB}\r\|_{X(\bx)}^{-2}                  \\
 &\lesssim2^{-2i\eta}
 \mu(2^iB)\lf\|{\bf 1}_{2^iB}\r\|_{X(\bx)}^{-2}.
\end{aligned}
\]
Here we used \(r_B\geq1\) in the last estimate. Also, in an analogous way,
\[
\begin{aligned}
 \mathbb I_{2,2}
 \leq
 C\|m\|_{L^2(\bx)}^2e^{-c(2^ir_B)^2}
\lesssim2^{-2i\eta}
 \mu(2^iB)\lf\|{\bf 1}_{2^iB}\r\|_{X(\bx)}^{-2}.
\end{aligned}
\]
By putting together the above estimates, we conclude that, for \(i\geq3\),
\[
 \lf\|S_L^\infty(m)\r\|_{L^2(S_i(B))}
 \lesssim2^{-i\eta}
 \mu(2^iB)^{1/2}
 \lf\|{\bf 1}_{2^iB}\r\|_{X(\bx)}^{-1}.
\]
By taking into account the \(L^2\)-boundedness of \(S_L^\infty\) and the
molecule condition, this estimation is also true for \(i=0,1,2\). Hence
\eqref{C1} is established for molecules of \(({\rm I})\)-type.

We now consider that \(m\) is a molecule of \(({\rm II})\)-type. Then
\(r_B\in(0,1)\) and there exists \(b\in L^2(\bx)\) such that
\(m=L^Mb\) and, for every \(k\in\{0,\ldots,M\}\),
\begin{equation}\label{acotb}
\lf\|L^kb\r\|_{L^2(S_i(B))}
\leq
2^{-i\varepsilon}r_B^{2(M-k)}
\mu(2^iB)^{1/2}
\lf\|{\bf 1}_{2^iB}\r\|_{X(\bx)}^{-1},
\quad i\in\zz_+.
\end{equation}
In particular,
\[
 \|b\|_{L^2(\bx)}
 \lesssim r_B^{2M}\mu(B)^{1/2}\|{\bf 1}_B\|_{X(\bx)}^{-1}.
\]
We have that
\[
\begin{aligned}
 \lf\|S_L^\infty(m)\r\|_{L^2(S_i(B))}^2
 &\lesssim\int_1^\infty
 \left\|
 {\bf 1}_{S_i(B)^{[t]}}
 t^2Le^{-t^2L}L^Mb
 \right\|_{L^2(\bx)}^2
 \frac{dt}{t}                                                \\
 &\thicksim
 \int_1^\infty
 t^{-4M}
 \left\|
 {\bf 1}_{S_i(B)^{[t]}}
 (t^2L)^{M+1}e^{-t^2L}b
 \right\|_{L^2(\bx)}^2
 \frac{dt}{t}.
\end{aligned}
\]

Let \(i\in\nn\), \(i\geq3\). We write
\[
 b=b_1+b_2,
 \qquad
 b_1=b{\bf 1}_{2^{i+2}B\setminus 2^{i-3}B}.
\]
Then we get
\[
\begin{aligned}
 \lf\|S_L^\infty(m)\r\|_{L^2(S_i(B))}^2
 \leq C\left(\lf\|S_L^\infty(L^Mb_1)\r\|_{L^2(S_i(B))}^2
 +\lf\|S_L^\infty(L^Mb_2)\r\|_{L^2(S_i(B))}^2\right)
 =:\mathbb J_1+\mathbb J_2.
\end{aligned}
\]
By using \eqref{acotb} and the \(L^2\)-boundedness of the Littlewood--Paley
function
\[
G_{L,M}(g)(x)
:=
\left(
\int_0^\infty
\lf|(t^2L)^{M+1}e^{-t^2L}g(x)\r|^2
\frac{dt}{t}
\right)^{1/2},
\]
we obtain
\[
\begin{aligned}
 \mathbb J_1
 &\lesssim\|b_1\|_{L^2(\bx)}^2
 \lesssim2^{-2i\varepsilon}r_B^{4M}\mu(2^iB)\lf\|{\bf 1}_{2^iB}\r\|_{X(\bx)}^{-2}                            \\
 &\lesssim2^{-2i\eta}
 \mu(2^iB)
 \lf\|{\bf 1}_{2^iB}\r\|_{X(\bx)}^{-2}.
\end{aligned}
\]
On the other hand, by \eqref{coneestimate} and \eqref{gapOD}, if
\(2^ir_B>C_0\), where \(C_0\) is a sufficiently large positive constant, then
\[
\begin{aligned}
 \mathbb J_2
 &\lesssim\|b\|_{L^2(\bx)}^2
 \left(\int_1^{c2^ir_B}t^{-4M}e^{-ct^2}\left(\frac{t}{2^ir_B}\right)^{4M}\frac{dt}{t}
 +\int_{c2^ir_B}^{\infty}t^{-4M}e^{-ct^2}\frac{dt}{t}\right)                                                    \\
 &\lesssim\lf(2^ir_B\r)^{-4M}\|b\|_{L^2(\bx)}^2
 \lesssim2^{-4Mi}\mu(B)\|{\bf 1}_B\|_{X(\bx)}^{-2}                           \\
 &\lesssim2^{-2i[2M-d(1/s_0-1/2)]}
 \mu(2^iB)\lf\|{\bf 1}_{2^iB}\r\|_{X(\bx)}^{-2}                   \\
 &\lesssim2^{-2i\eta}
 \mu(2^iB)\lf\|{\bf 1}_{2^iB}\r\|_{X(\bx)}^{-2}.
\end{aligned}
\]
If \(2^ir_B\leq C_0\), then, by the spectral theorem and the above estimate
for \(\|b\|_{L^2(\bx)}\), we have
\begin{align*}
\lf\|S_L^\infty(m)\r\|_{L^2(S_i(B))}
&\lesssim\|b\|_{L^2(\bx)}
\lesssim r_B^{2M}\mu(B)^{1/2}\|{\bf 1}_B\|_{X(\bx)}^{-1}    \\
&\lesssim2^{-2Mi}\mu(B)^{1/2}\|{\bf 1}_B\|_{X(\bx)}^{-1}   \\
&\lesssim2^{-i[2M-d(1/s_0-1/2)]}\mu(2^iB)^{1/2}\lf\|{\bf 1}_{2^iB}\r\|_{X(\bx)}^{-1} \\
&\lesssim2^{-i\eta}\mu(2^iB)^{1/2}\lf\|{\bf 1}_{2^iB}\r\|_{X(\bx)}^{-1}.
\end{align*}
Combining  with the above estimates, we have
\[
\lf\|S_L^\infty(m)\r\|_{L^2(S_i(B))}
\lesssim2^{-i\eta}
\mu(2^iB)^{1/2}
\lf\|{\bf 1}_{2^iB}\r\|_{X(\bx)}^{-1},
\quad i\in\nn,\ i\geq3,
\]
with
\[
 \eta>d\left(\frac1{s_0}-\frac1{q_0}\right).
\]
When \(i=0,1,2\), as in the previous case, we also get this estimation.
Thus we have proved that \eqref{C1} holds for local molecules of
\(({\rm II})\)-type.

Therefore, by Lemma \ref{l2.8}, we obtain
\[
\lf\|S_L^\infty(f)\r\|_{X(\bx)}
\lesssim
\left\|
\left\{
\sum_{j\in\nn}
\left[
\lambda_j
\frac{{\bf 1}_{B_j}}{\|{\bf 1}_{B_j}\|_{X(\bx)}}
\right]^{s_0}
\right\}^{1/s_0}
\right\|_{X(\bx)}.
\]
Consequently,
\[
\|f\|_{H^X_L(\bx)}
=
\lf\|S_L(f)\r\|_{X(\bx)}
\lesssim
\left\|
\left\{
\sum_{j\in\nn}
\left[
\lambda_j
\frac{{\bf 1}_{B_j}}{\|{\bf 1}_{B_j}\|_{X(\bx)}}
\right]^{s_0}
\right\}^{1/s_0}
\right\|_{X(\bx)}.
\]
This proves
\[
 h^X_L(\bx)\subset H^X_L(\bx).
\]

Conversely, suppose that \(f\in H^X_L(\bx)\cap L^2(\bx)\). By the molecular
characterization of \(H^X_L(\bx)\) (see \cite{lyyy23}), we know that there exist, for every \(j\in\nn\),
\(\lambda_j>0\) and an \((X,2,M,\varepsilon)_L\)-molecule \(m_j\) associated
with the ball \(B_j\) such that
\[
 f=\sum_{j\in\nn}\lambda_jm_j
 \ \ \hbox{in } \ \ L^2(\bx),
\]
and
\[
\left\|
\left\{
\sum_{j\in\nn}
\left[
\lambda_j
\frac{{\bf 1}_{B_j}}{\|{\bf 1}_{B_j}\|_{X(\bx)}}
\right]^{s_0}
\right\}^{1/s_0}
\right\|_{X(\bx)}
\lesssim\|f\|_{H^X_L(\bx)}.
\]
Moreover, we observe that every global \((X,2,M,\varepsilon)_L\)-molecule is a local
\((X,2,M,\varepsilon)\)-molecule. Indeed, if \(r_{B_j}<1\), it is a molecule
of \(({\rm II})\)-type; if \(r_{B_j}\geq1\), then the estimate corresponding
to \(k=M\) gives
\[
 \lf\|m_j\r\|_{L^2(S_i(B_j))}
 \lesssim2^{-i\varepsilon}
 \mu(2^iB_j)^{1/2}
 \lf\|{\bf 1}_{2^iB_j}\r\|_{X(\bx)}^{-1},
 \quad i\in\zz_+,
\]
and hence \(m_j\) is a molecule of \(({\rm I})\)-type. Therefore, by
Theorem \ref{mol}, we obtain
\[
 \|f\|_{h^X_L(\bx)}
 \lesssim\|f\|_{H^X_L(\bx)}.
\]
The reverse proof can be completed by using a density argument. Thus
\[
 h^X_L(\bx)=H^X_L(\bx)
\]
algebraically and topologically. This completes the proof of Theorem \ref{eq0}.
\end{proof}

Finally, we use Theorem \ref{mol} to prove Theorem \ref{eq1} with the help of
some ideas from \cite[Theorem 1.6]{abdr20}.
\begin{proof}[Proof of Theorem \ref{eq1}]
Firstly, we prove that
\[
 H^X_{L+I}(\bx)\subset h^X_L(\bx).
\]
Let \(f\in H^X_{L+I}(\bx)\cap L^2(\bx)\). Consider
\(M\in\nn\) such that
$
 2M>d(\frac2{s_0}-\frac12-\frac1{q_0})+\frac12
$
and
$
 \varepsilon>d(\frac1{s_0}-\frac1{q_0}).
$
According to the molecular characterization of \(H^X_{L+I}(\bx)\), for every
\(j\in\nn\), there exist \(\lambda_j>0\) and an
\((X,2,M,\varepsilon)_{L+I}\)-molecule \(m_j\) associated to the ball
\(B_j=B(x_{B_j},r_{B_j})\), with \(x_{B_j}\in\bx\) and \(r_{B_j}>0\), such
that
\[
 f=\sum_{j\in\nn}\lambda_jm_j \ \ \mathrm{in} \ \ L^2(\bx) \ \mathrm{and} \ \ H^X_{L+I}(\bx),
\]
 and
\[
\left\|
\left\{
\sum_{j\in\nn}
\left[
\lambda_j
\frac{{\bf 1}_{B_j}}{\|{\bf 1}_{B_j}\|_{X(\bx)}}
\right]^{s_0}
\right\}^{1/s_0}
\right\|_{X(\bx)}
\sim
\|f\|_{H^X_{L+I}(\bx)}.
\]
We only need to show that \(S_L^{\rm loc}(f)\in X(\bx)\),
\(e^{-L}f\in H^X_I(\bx)\), and
\[
\lf\|S_L^{\rm loc}(f)\r\|_{X(\bx)}
+
\lf\|S_I(e^{-L}f)\r\|_{X(\bx)}
\lesssim
\left\|
\left\{
\sum_{j\in\nn}
\left[
\lambda_j
\frac{{\bf 1}_{B_j}}{\|{\bf 1}_{B_j}\|_{X(\bx)}}
\right]^{s_0}
\right\}^{1/s_0}
\right\|_{X(\bx)}.
\]
By applying Lemma \ref{l2.8}, it is sufficient to show that there exist
\(C>0\) and
$
 \eta>d(\frac1{s_0}-\frac1{q_0})
$
such that, for every \(j\in\nn\),
\begin{align}\label{D1}
\lf\|S_L^{\rm loc}(m_j)\r\|_{L^2(S_i(B_j))}
+
\lf\|S_I(e^{-L}m_j)\r\|_{L^2(S_i(B_j))}
\lesssim2^{-i\eta}
\mu(2^iB_j)^{1/2}
\lf\|{\bf 1}_{2^iB_j}\r\|_{X(\bx)}^{-1},
\qquad i\in\zz_+.
\end{align}

Let \(m\) be an \((X,2,M,\varepsilon)_{L+I}\)-molecule associated to
\(B=B(x_B,r_B)\), with \(x_B\in\bx\) and \(r_B>0\). Then there exists
\(b\in L^2(\bx)\) such that \(m=(L+I)^Mb\) and, for every
\(k\in\{0,\ldots,M\}\),
\begin{equation}\label{Lkb2}
\lf\|(L+I)^kb\r\|_{L^2(S_i(B))}
\leq
2^{-i\varepsilon}r_B^{2(M-k)}
\mu(2^iB)^{1/2}
\lf\|{\bf 1}_{2^iB}\r\|_{X(\bx)}^{-1},
\qquad i\in\zz_+.
\end{equation}
We note that if \(r_B\geq1\), then \(m\) is also a local
\((X,2,M,\varepsilon)_L\)-molecule of \(({\rm I})\)-type. Then, according to
the estimates established in the proof of Theorem \ref{mol}, \eqref{D1}
holds for \(m\), provided that \(r_B\geq1\).
Assume now that \(r_B\in(0,1)\). We can write
\begin{equation}\label{D2}
Le^{-t^2L}-(L+I)e^{-t^2(L+I)}
=
(1-e^{-t^2})Le^{-t^2L}
-
e^{-t^2}e^{-t^2L},
\qquad t>0.
\end{equation}
Thus we have that, for every \(g\in L^2(\bx)\),
\[
 S_L^{\rm loc}(g)
 \leq
 S_{L+I}^{\rm loc}(g)+{\mathcal H}(g)+{\mathcal Q}(g),
\]
where
\[
{\mathcal H}(g)(x)
:=
\left(
\int_0^1\int_{B(x,t)}
\lf|t^2(1-e^{-t^2})Le^{-t^2L}g(y)\r|^2
\frac{d\mu(y)dt}{V_t(x)t}
\right)^{1/2},
\qquad x\in\bx,
\]
and
\[
{\mathcal Q}(g)(x)
:=
\left(
\int_0^1\int_{B(x,t)}
\lf|t^2e^{-t^2}e^{-t^2L}g(y)\r|^2
\frac{d\mu(y)dt}{V_t(x)t}
\right)^{1/2},
\qquad x\in\bx.
\]

According to the molecular estimates for \(H^X_{L+I}(\bx)\), we can find
\(C>0\) and
$
 \eta>d(\frac1{s_0}-\frac1{q_0})
$
such that
\[
\lf\|S_{L+I}^{\rm loc}(m)\r\|_{L^2(S_i(B))}
\leq
\lf\|S_{L+I}(m)\r\|_{L^2(S_i(B))}
\lesssim2^{-i\eta}
\mu(2^iB)^{1/2}
\lf\|{\bf 1}_{2^iB}\r\|_{X(\bx)}^{-1},
\ \  i\in\zz_+.
\]

For the operator \({\mathcal H}\), we can write
\begin{align*}
\|{\mathcal H}(m)\|_{L^2(S_i(B))}^2
&\leq
C\int_{S_i(B)}
\left(
\int_0^{2^{i/2}r_B}
+
\int_{2^{i/2}r_B}^1
\right)
\int_{B(x,t)}
\lf|t^4Le^{-t^2L}m(y)\r|^2
\frac{d\mu(y)dt}{V_t(x)t}
d\mu(x)                                      \\
&=:\mathrm{L}_1+\mathrm{L}_2,
\end{align*}
provided that \(2^{i/2}r_B<1\), and, in the other case, \(\mathrm{L}_2=0\).

Let \(i\in\nn\), \(i\geq3\). We observe that if \(x\in S_i(B)\),
\(t\in(0,2^{i/2}r_B)\), \(y\in B(x,t)\), and
\(z\in(2^{i+2}B\setminus2^{i-3}B)^\complement\), then
$
 d(y,z)\geq c2^ir_B.
$
From the decomposition
\[
 m=m_1+m_2,\qquad
 m_1=m{\bf 1}_{2^{i+2}B\setminus2^{i-3}B},
\]
and the \(L^2\)-boundedness of \(S_L\), the off-diagonal estimates
deduced from Assumption \ref{a2.5}, and \eqref{Lkb2} with \(k=M\), we conclude
that, for every sufficiently large \(\gamma>0\),
\begin{align*}
 \mathrm{L}_1
&\lesssim
\lf\|S_L(m_1)\r\|_{L^2(\bx)}^2
+
\|m\|_{L^2(\bx)}^2
\int_0^{2^{i/2}r_B}
t^4
\left(\frac{t}{2^ir_B}\right)^\gamma
\frac{dt}{t}                                                  \\
&\lesssim
2^{-2i\varepsilon}
\mu(2^iB)
\lf\|{\bf 1}_{2^iB}\r\|_{X(\bx)}^{-2}
+
\|m\|_{L^2(\bx)}^2
\frac{(2^{i/2}r_B)^{\gamma+4}}{(2^ir_B)^\gamma}
.
\end{align*}
By \eqref{Lkb2} and Assumption \ref{a12}, we have
\[
 \|m\|_{L^2(\bx)}
 \lesssim\mu(B)^{1/2}\|{\bf 1}_B\|_{X(\bx)}^{-1}
 \lesssim2^{id(1/s_0-1/2)}
 \mu(2^iB)^{1/2}\lf\|{\bf 1}_{2^iB}\r\|_{X(\bx)}^{-1}.
\]
Consequently,
\begin{align*}
 \mathrm{L}_1
&\lesssim
\mu(2^iB)
\lf\|{\bf 1}_{2^iB}\r\|_{X(\bx)}^{-2}
\left[
2^{-2i\varepsilon}
+
2^{2id(1/s_0-1/2)}
2^{-i(\gamma/2-2)}r_B^4
\right]                                                   \\
&\lesssim2^{-2i\eta}
\mu(2^iB)
\lf\|{\bf 1}_{2^iB}\r\|_{X(\bx)}^{-2},
\end{align*}
where \(\gamma\) is chosen large enough and
$
 \eta>d(\frac1{s_0}-\frac1{q_0}).
$

On the other hand, by the spectral theorem, it follows that the family of
operators
$
\{
 e^{-t^2L/2}(L+I)^Mt^{2M}e^{-t^2L/2}
\}_{t\in(0,1)}
$
is bounded in \({\mathcal L}(L^2(\bx))\). Then we have that
\begin{align*}
\mathrm{L}_2
&=
\int_{S_i(B)}
\int_{2^{i/2}r_B}^1
\int_{B(x,t)}
\lf|t^4Le^{-t^2L}(L+I)^Mb(y)\r|^2
\frac{d\mu(y)dt}{V_t(x)t}
d\mu(x)                                                  \\
&\lesssim
\int_{2^{i/2}r_B}^1
\int_{\bx}
\lf|t^2Le^{-t^2L/2}(L+I)^Mt^{2M}e^{-t^2L/2}b(y)\r|^2
\frac{d\mu(y)dt}{t^{4M-3}}                               \\
&\lesssim
\frac{\|b\|_{L^2(\bx)}^2}{(2^{i/2}r_B)^{4M-2}}.
\end{align*}
From \eqref{Lkb2} with \(k=0\), it follows that
\[
 \|b\|_{L^2(\bx)}
 \lesssim r_B^{2M}\mu(B)^{1/2}\|{\bf 1}_B\|_{X(\bx)}^{-1}.
\]
Thus, by Assumption \ref{a12} and the condition
$
 M>d(\frac2{s_0}-\frac12-\frac1{q_0})+\frac12
$, we obtain
\begin{align*}
\mathrm{L}_2
&\lesssim
\frac{r_B^{4M}\mu(B)\|{\bf 1}_B\|_{X(\bx)}^{-2}}
     {(2^{i/2}r_B)^{4M-2}}
\lesssim
2^{-i(2M-1)}
\mu(B)\|{\bf 1}_B\|_{X(\bx)}^{-2}                          \\
&\lesssim
2^{-2i[M-1/2-d(1/s_0-1/2)]}
\mu(2^iB)\lf\|{\bf 1}_{2^iB}\r\|_{X(\bx)}^{-2}
\lesssim2^{-2i\eta}
\mu(2^iB)\lf\|{\bf 1}_{2^iB}\r\|_{X(\bx)}^{-2}.
\end{align*}
Combining with the estimates of \(\mathrm{L}_1\) and \(\mathrm{L}_2\), we conclude that
\[
\|{\mathcal H}(m)\|_{L^2(S_i(B))}
\lesssim2^{-i\eta}
\mu(2^iB)^{1/2}
\lf\|{\bf 1}_{2^iB}\r\|_{X(\bx)}^{-1},
\qquad i\in\nn,\ i\geq3.
\]
Also, the \(L^2\)-boundedness of \(S_L\) and \eqref{Lkb2} give this
estimation also for \(i=0,1,2\), because, in these cases,
\[
\|{\mathcal H}(m)\|_{L^2(S_i(B))}
\leq
\|{\mathcal H}(m)\|_{L^2(\rr^d)}
\lesssim\|S_L(m)\|_{L^2(\bx)}
\lesssim\|m\|_{L^2(\bx)}
\lesssim
\mu(B)^{1/2}
\|{\bf 1}_B\|_{X(\bx)}^{-1}.
\]

In order to estimate \({\mathcal Q}\), we proceed in a similar way.
We observe that the operator \({\mathcal Q}\) is bounded from \(L^2(\bx)\)
into itself. Indeed, for every \(g\in L^2(\bx)\), by the Fubini theorem, the
doubling condition, and the \(L^2\)-boundedness of \(e^{-t^2L}\), we obtain
\begin{align*}
\lf\|{\mathcal Q}(g)\r\|_{L^2(\bx)}^2
&=
\int_{\bx}
\int_0^1
\int_{B(x,t)}
\lf|t^2e^{-t^2}e^{-t^2L}g(y)\r|^2
\frac{d\mu(y)dt}{V_t(x)t}
d\mu(x)                                                     \\
&\lesssim\int_0^1
t^3e^{-2t^2}
\lf\|e^{-t^2L}g\r\|_{L^2(\bx)}^2
dt
\lesssim\|g\|_{L^2(\bx)}^2.
\end{align*}
This fact, the \({\mathcal L}(L^2(\bx))\)-boundedness of the family of
operators
\[
\left\{
e^{-t^2L/2}(L+I)^Mt^{2M}e^{-t^2L/2}
\right\}_{t\in(0,1)}
\]
and the Gaussian upper estimate for the heat kernel associated with \(L\)
allow us to show, in the same way as above, that
\[
\lf\|{\mathcal Q}(m)\r\|_{L^2(S_i(B))}
\lesssim2^{-i\eta}
\mu(2^iB)^{1/2}
\lf\|{\bf 1}_{2^iB}\r\|_{X(\bx)}^{-1},
\qquad i\in\zz_+,
\]
for some
$
 \eta>d(\frac1{s_0}-\frac1{q_0}),
$
which, together with  the above estimates, implies that
\[
\lf\|S_L^{\rm loc}(m)\r\|_{L^2(S_i(B))}
\lesssim2^{-i\eta}
\mu(2^iB)^{1/2}
\lf\|{\bf 1}_{2^iB}\r\|_{X(\bx)}^{-1},
\qquad i\in\zz_+.
\]
Moreover, we observe that the operator \((L+I)^Me^{-L}\) is
bounded from \(L^2(\bx)\) into itself.
Let \(K_M(x,y)\) denote the integral kernel of
$
(L+I)^M e^{-L}.
$
Since
\[
(L+I)^M e^{-L}
=
\sum_{j=0}^M\binom{M}{j}L^j e^{-L},
\]
the Gaussian kernel estimates for \(L^j e^{-L}\),
\(j\in\{0,\ldots,M\}\), imply that there exist constants
\(C_M,c_M>0\) such that
\[
\left|K_M(x,y)\right|
\leq
\frac{C_M}{\mu(B(x,1))}
\exp\left(-c_M d(x,y)^2\right)
\]
for all \(x,y\in\mathbb X\).
This further implies that, for every \(N>0\), the operator
\((L+I)^Me^{-L}\) satisfies the following off-diagonal estimate: for all closed sets \(E,F\subset\bx\),
\[
\left\|
{\bf 1}_E(L+I)^Me^{-L}{\bf 1}_F
\right\|_{L^2(\bx)\to L^2(\bx)}
\leq
C_N
(1+\dist(E,F))^{-N}.
\]
 By proceeding as in the proof of
Theorem \ref{mol}, we deduce that, for certain
$
 \eta>d(\frac1{s_0}-\frac1{q_0}),
$
\[
\lf\|S_I(e^{-L}m)\r\|_{L^2(S_i(B))}
\lesssim2^{-i\eta}
\mu(2^iB)^{1/2}
\lf\|{\bf 1}_{2^iB}\r\|_{X(\bx)}^{-1},
\qquad i\in\zz_+.
\]
Thus \eqref{D1} is established. A density argument allows us to conclude that
\[
 H^X_{L+I}(\bx)\subset h^X_L(\bx)
\]
algebraically and topologically.

Finally, we prove that
\[
 h^X_L(\bx)\subset H^X_{L+I}(\bx).
\]
Let us consider \(M\in\nn\) such that
\[
 2(M-1)>d\left(\frac2{s_0}-\frac12-\frac1{q_0}\right)
\ \
\mathrm{and}
\ \
 \varepsilon>d\left(\frac1{s_0}-\frac1{q_0}\right).
\]
Let \(f\in h^X_L(\bx)\cap L^2(\bx)\). According to Theorem \ref{mol}, for
every \(j\in\nn\), there exist \(\lambda_j>0\) and a local
\((X,2,M,\varepsilon)_L\)-molecule \(m_j\) associated to the ball \(B_j\)
such that
\[
 f=\sum_{j\in\nn}\lambda_jm_j \ \ \mathrm{in} \ \ L^2(\bx) \ \ \mathrm{and} \ \ h^X_L(\bx),
\]
 and
\[
\left\|
\left\{
\sum_{j\in\nn}
\left[
\lambda_j
\frac{{\bf 1}_{B_j}}{\|{\bf 1}_{B_j}\|_{X(\bx)}}
\right]^{s_0}
\right\}^{1/s_0}
\right\|_{X(\bx)}
\lesssim\|f\|_{h^X_L(\bx)}.
\]
Our objective is to show that \(f\in H^X_{L+I}(\bx)\) and
\[
\lf\|S_{L+I}(f)\r\|_{X(\bx)}
\leq
C
\left\|
\left\{
\sum_{j\in\nn}
\left[
\lambda_j
\frac{{\bf 1}_{B_j}}{\lf\|{\bf 1}_{B_j}\r\|_{X(\bx)}}
\right]^{s_0}
\right\}^{1/s_0}
\right\|_{X(\bx)}.
\]
According to Lemma \ref{l2.8}, we only need to show that,
for every \(j\in\nn\),
\[
\lf\|S_{L+I}(m_j)\r\|_{L^2(S_i(B_j))}
\lesssim2^{-i\eta}
\mu(2^iB_j)^{1/2}
\lf\|{\bf 1}_{2^iB_j}\r\|_{X(\bx)}^{-1},
\qquad i\in\zz_+,
\]
for some
$
 \eta>d(\frac1{s_0}-\frac1{q_0}).
$

Assume that \(m\) is a local \((X,2,M,\varepsilon)_L\)-molecule associated
to the ball \(B=B(x_B,r_B)\), with \(x_B\in\bx\) and \(r_B>0\). We can write
\[
 S_{L+I}(m)
 \leq
 S_{L+I}^{\rm loc}(m)+S_{L+I}^{\infty}(m),
\]
where
\[
S_{L+I}^{\infty}(m)(x)
:=
\left(
\int_1^\infty
\int_{B(x,t)}
\lf|t^2(L+I)e^{-t^2(L+I)}m(y)\r|^2
\frac{d\mu(y)dt}{V_t(x)t}
\right)^{1/2},
\qquad x\in\bx.
\]

Suppose that \(r_B\geq1\). Then \(m\) is also a local
\((X,2,M,\varepsilon)_{L+I}\)-molecule. Hence, by the estimates established
above with \(L+I\) instead of \(L\),
\[
\lf\|S_{L+I}^{\rm loc}(m)\r\|_{L^2(S_i(B))}
\lesssim2^{-i\eta}
\mu(2^iB)^{1/2}
\lf\|{\bf 1}_{2^iB}\r\|_{X(\bx)}^{-1},
\qquad i\in\zz_+,
\]
for some
$
 \eta>d(\frac1{s_0}-\frac1{q_0}).
$
On the other hand, since
\[
 e^{-t(L+I)}=e^{-t}e^{-tL},
 \qquad t>0,
\]
the heat kernel associated with \(L+I\) satisfies the same Gaussian upper
estimate with the additional factor \(e^{-t}\). Moreover, by the spectral
theorem, for every \(N>0\),
\[
\left\|
{\bf 1}_E(t^2(L+I))e^{-t^2(L+I)}{\bf 1}_F
\right\|_{L^2(\bx)\to L^2(\bx)}
\leq
C_Ne^{-ct^2}
\left(1+\frac{\dist(E,F)}{t}\right)^{-N}.
\]
Since \(S_{L+I}\) is a bounded sublinear operator from \(L^2(\bx)\) into
itself, we can proceed as in the proof of \eqref{D1} for molecules of
\(({\rm I})\)-type to obtain that
\[
\lf\|S_{L+I}^{\infty}(m)\r\|_{L^2(S_i(B))}
\lesssim2^{-i\eta}
\mu(2^iB)^{1/2}
\lf\|{\bf 1}_{2^iB}\r\|_{X(\bx)}^{-1},
\qquad i\in\zz_+,
\]
for some
$
 \eta>d(\frac1{s_0}-\frac1{q_0}).
$

Suppose now that \(r_B\in(0,1)\). There exists \(b\in L^2(\bx)\) such that
\(m=L^Mb\) and, for every \(k\in\{0,\ldots,M\}\),
\[
\lf\|L^kb\r\|_{L^2(S_i(B))}
\leq
2^{-i\varepsilon}r_B^{2(M-k)}
\mu(2^iB)^{1/2}
\lf\|{\bf 1}_{2^iB}\r\|_{X(\bx)}^{-1},
\qquad i\in\zz_+.
\]
We define
$
 \mathfrak b:=L^{M-1}b.
$
Since \(r_B\in(0,1)\), \(\mathfrak b\) is a local
\((X,2,M-1,\varepsilon)_L\)-molecule associated with the ball \(B\) of
\(({\rm II})\)-type. It is clear that
$m=L\mathfrak b.$
We can write
\begin{align*}
\lf[S_{L+I}^{\rm loc}(m)(x)\r]^2
&=
\int_0^1
\int_{B(x,t)}
\lf|t^2(L+I)e^{-t^2(L+I)}m(y)\r|^2
\frac{d\mu(y)dt}{V_t(x)t}                                  \\
&\lesssim\lf[
\int_0^1
\int_{B(x,t)}
\lf|e^{-t^2}t^2Le^{-t^2L}m(y)\r|^2
\frac{d\mu(y)dt}{V_t(x)t}
+
\int_0^1
\int_{B(x,t)}
\lf|e^{-t^2}t^2Le^{-t^2L}\mathfrak b(y)\r|^2
\frac{d\mu(y)dt}{V_t(x)t}\r]                                                      \\
&\lesssim\left[
\lf(S_L^{\rm loc}(m)(x)\r)^2+
\lf(S_L^{\rm loc}(\mathfrak b)(x)\r)^2
\right],
\qquad x\in\bx.
\end{align*}
According to the estimates established in the proof of Theorem \ref{mol}, we
conclude that
\begin{align*}
\lf\|S_{L+I}^{\rm loc}(m)\r\|_{L^2(S_i(B))}
&\lesssim\left(
\lf\|S_L^{\rm loc}(m)\r\|_{L^2(S_i(B))}
+
\lf\|S_L^{\rm loc}(\mathfrak b)\r\|_{L^2(S_i(B))}
\right)                                                    \\
&\lesssim2^{-i\eta}
\mu(2^iB)^{1/2}
\lf\|{\bf 1}_{2^iB}\r\|_{X(\bx)}^{-1},
\qquad i\in\zz_+,
\end{align*}
for some
$
 \eta>d(\frac1{s_0}-\frac1{q_0}).
$
Moreover, for every \(i\in\zz_+\), we have
\begin{align*}
\lf\|S_{L+I}^{\infty}(m)\r\|_{L^2(S_i(B))}^2
&\lesssim\int_{S_i(B)}
\int_1^\infty
\int_{B(x,t)}
\lf|t^2(L+I)e^{-t^2(L+I)}L^Mb(y)\r|^2
\frac{d\mu(y)dt}{V_t(x)t}
d\mu(x)                                                     \\
&\lesssim\int_{S_i(B)}
\int_1^\infty
\int_{B(x,t)}
\lf|e^{-t^2}(t^2L)^{M+1}e^{-t^2L}b(y)\r|^2
\frac{d\mu(y)dt}{V_t(x)t^{4M+1}}
d\mu(x)                                                     \\
&\quad+
\int_{S_i(B)}
\int_1^\infty
\int_{B(x,t)}
\lf|e^{-t^2}(t^2L)^Me^{-t^2L}b(y)\r|^2
\frac{d\mu(y)dt}{V_t(x)t^{4M-3}}
d\mu(x).
\end{align*}
From the proof of \eqref{D1} for molecules of
\(({\rm II})\)-type, we deduce that
\[
\lf\|S_{L+I}^{\infty}(m)\r\|_{L^2(S_i(B))}
\lesssim2^{-i\eta}
\mu(2^iB)^{1/2}
\lf\|{\bf 1}_{2^iB}\r\|_{X(\bx)}^{-1},
\qquad i\in\zz_+.
\]
Combining with the above estimates, we have
\[
\lf\|S_{L+I}(m)\r\|_{L^2(S_i(B))}
\lesssim2^{-i\eta}
\mu(2^iB)^{1/2}
\lf\|{\bf 1}_{2^iB}\r\|_{X(\bx)}^{-1},
\qquad i\in\zz_+.
\]

Therefore, we can deduce that \(f\in H^X_{L+I}(\bx)\) and
\[
\|f\|_{H^X_{L+I}(\bx)}
\lesssim
\left\|
\left\{
\sum_{j\in\nn}
\left[
\lambda_j
\frac{{\bf 1}_{B_j}}{\|{\bf 1}_{B_j}\|_{X(\bx)}}
\right]^{s_0}
\right\}^{1/s_0}
\right\|_{X(\bx)}.
\]
By using a standard density arguments, we obtain that
$h^X_L(\bx)\subset H^X_{L+I}(\bx)
$
and the inclusion is continuous. This finishes the proof of Theorem
\ref{eq1}.
\end{proof}

\section{Applications}\label{s6x}
In this section, we will establish some applications. The key tools are the molecular characterization of the local Hardy space \(h_L^X(\bx)\) and the relation between \(h_L^X(\bx)\) and the global Hardy spaces associated with the shifted operator \(L+I\). We apply these characterizations to establish the boundedness of local Riesz transforms and H\"ormander-type spectral multipliers on local Hardy spaces.

\subsection{Boundedness of local Riesz transforms}\label{s61}
In this subsection, we restrict ourselves to the Euclidean setting
\(\mathbb X=\mathbb R^d\) and to second-order divergence-form elliptic
operators
\[
L=-\operatorname{div}(A\nabla),
\]
where \(A\) is a real symmetric matrix satisfying the uniform ellipticity
condition.
 Let
\(A\in L^\infty(\rr^d;\mathbb R^{d\times d})\) be real symmetric and
uniformly elliptic; namely, there exist constants
\(0<\lambda\leq\Lambda<\infty\) such that, for a.e. \(x\in\rr^d\) and all
\(\xi\in\mathbb C^d\),
\[
 \lambda|\xi|^2
 \leq
 A(x)\xi\cdot\overline{\xi}
 \leq
 \Lambda|\xi|^2 .
\]
Let \(L:=-\operatorname{div}(A\nabla)\) be the non-negative self-adjoint
operator on \(L^2(\rr^d)\) associated with the closed quadratic form
\[
 \mathfrak a(f,g)
 :=
 \int_{\rr^d}A(x)\nabla f(x)\cdot\overline{\nabla g(x)}\,dx,
 \qquad
 f,g\in W^{1,2}(\rr^d).
\]
Put
\[
        T_{\kappa}:=L+\kappa I.
\]
Then we say that $\mathcal R_{L,\kappa}^{\rm loc}:=\nabla(L+\kappa I)^{-1/2}:=\nabla T_{\kappa}^{-1/2}$
is the local, or inhomogeneous, Riesz transform associated with $L$.

In what follows,
for a vector-valued function \(F\), we write
\[
X(\mathbb R^d;\mathbb C^d)
:=
\left\{
F:\mathbb R^d\to\mathbb C^d:
|F|_{\mathbb C^d}\in X(\mathbb R^d)
\right\},
\]
equipped with
\[
\|F\|_{X(\mathbb R^d;\mathbb C^d)}
:=
\big\||F|_{\mathbb C^d}\big\|_{X(\mathbb R^d)}.
\]

\begin{theorem}\label{thm:local-Riesz}
Let \(X(\mathbb R^d)\) be a ball quasi-Banach function space satisfying
Assumptions \ref{a11} and \ref{a12} for some
$
p\in(0,\infty),
s_0\in(0,\min\{p,1\}],
q_0\in(s_0,2].
$
Let
$
L=-\operatorname{div}(A\nabla)
$
on \(\mathbb R^d\), where \(A\) is real symmetric, bounded, and uniformly
elliptic. For every fixed \(\kappa>0\), put
\[
T_\kappa:=L+\kappa I
\]
and define
\[
\mathcal R_{L,\kappa}^{\rm loc}
:=
\nabla T_\kappa^{-1/2}.
\]
Then \(\mathcal R_{L,\kappa}^{\rm loc}\) extends uniquely to a bounded
operator from \(h_L^X(\mathbb R^d)\) to
\(X(\mathbb R^d;\mathbb C^d)\). More precisely, there exists
\(C_\kappa>0\) such that, for every
\(f\in\mathbf h_L^X(\mathbb R^d)\),
\[
\left\|
\mathcal R_{L,\kappa}^{\rm loc}f
\right\|_{X(\mathbb R^d;\mathbb C^d)}
\leq
C_\kappa
\|f\|_{h_L^X(\mathbb R^d)}.
\]
\end{theorem}

The following lemma is inspired by
\cite[Lemmas 2.3 and 2.4]{hm09}.

\begin{lemma}\label{grad}
Let
$
L=-\operatorname{div}(A\nabla)
$
be the non-negative self-adjoint operator on \(L^2(\mathbb R^d)\)
associated with a real symmetric, bounded, and uniformly elliptic
matrix \(A\). For every \(\kappa>0\), put
$
T_\kappa:=L+\kappa I.
$
Then, for every \(k\in\mathbb Z_+\) and every \(N>0\), there exists
a constant \(C_{N,k}>0\), independent of
\(\kappa\), \(E\), \(F\), and \(t\), such that, for all Borel sets
\(E,F\subset\mathbb R^d\) and all \(t>0\),
\[
\left\|
{\bf 1}_E\,t\nabla(t^2L)^k
e^{-t^2T_\kappa}{\bf 1}_F
\right\|_{L^2(\mathbb R^d)
\to L^2(\mathbb R^d;\mathbb C^d)}
\leq
C_{N,k}e^{-\kappa t^2}
\left(
1+\frac{\operatorname{dist}(E,F)}{t}
\right)^{-N}.
\]
Moreover,
\[
\left\|
{\bf 1}_E\,t\nabla(t^2T_\kappa)^k
e^{-t^2T_\kappa}{\bf 1}_F
\right\|_{L^2(\mathbb R^d)
\to L^2(\mathbb R^d;\mathbb C^d)}
\leq
C_{N,k}e^{-\kappa t^2/2}
\left(
1+\frac{\operatorname{dist}(E,F)}{t}
\right)^{-N}.
\]
\end{lemma}

\begin{proof}
Firstly, we establish the higher-order gradient Gaffney estimate
\begin{equation}\label{eq:higher-gradient-Gaffney}
\left\|
{\bf 1}_E\,t\nabla(t^2L)^\ell
e^{-t^2L}{\bf 1}_F
\right\|_{L^2(\mathbb R^d)
\to L^2(\mathbb R^d;\mathbb C^d)}
\leq
C_\ell
\exp\left(
-c_\ell\frac{[\operatorname{dist}(E,F)]^2}{t^2}
\right)
\end{equation}
for every \(\ell\in\mathbb Z_+\).
Indeed, by \cite{hm09}, the standard \(L^2\)-Gaffney estimates for divergence-form
elliptic operators state that, for all Borel sets
\(G,H\subset\mathbb R^d\) and all \(s>0\),
\begin{align}
&
\left\|
{\bf 1}_G\,\sqrt{s}\nabla e^{-sL}{\bf 1}_H
\right\|_{L^2(\mathbb R^d)
\to L^2(\mathbb R^d;\mathbb C^d)}
\lesssim
\exp\left(
-c\frac{\operatorname{dist}(G,H)^2}{s}
\right),
\label{eq:gradient-Gaffney}
\\
&
\left\|
{\bf 1}_G\,sLe^{-sL}{\bf 1}_H
\right\|_{L^2(\mathbb R^d)
\to L^2(\mathbb R^d)}
\lesssim
\exp\left(
-c\frac{\operatorname{dist}(G,H)^2}{s}
\right).
\label{eq:L-Gaffney}
\end{align}
Let us recall briefly the composition property used below. Suppose that
two \(L^2\)-bounded operators \(S_t\) and \(R_t\) satisfy Gaussian
off-diagonal estimates at scale \(t\). If \(\operatorname{dist}(E,F)>0\), set
\[
\Omega
:=
\left\{
x\in\mathbb R^d:
\operatorname{dist}(x,E)<\frac {\operatorname{dist}(E,F)}{2}
\right\}.
\]
Then
\[
\operatorname{dist}(E,\Omega^c)\geq\frac {\operatorname{dist}(E,F)}{2},
\qquad
\operatorname{dist}(\Omega,F)\geq\frac {\operatorname{dist}(E,F)}{2},
\]
and
\[
{\bf 1}_ES_tR_t{\bf 1}_F
=
{\bf 1}_ES_t{\bf 1}_\Omega R_t{\bf 1}_F
+
{\bf 1}_ES_t{\bf 1}_{\Omega^c}R_t{\bf 1}_F.
\]
For the first term, since
\[
\operatorname{dist}(\Omega,F)
\geq
\frac{\operatorname{dist}(E,F)}{2},
\]
we apply the off-diagonal estimate to
$
{\bf 1}_{\Omega}R_t{\bf 1}_F
$
and use the global \(L^2\)-boundedness of \(S_t\). For the second term, since
$
\operatorname{dist}(E,\Omega^c)
\geq
\frac{\operatorname{dist}(E,F)}{2},
$
we apply the off-diagonal estimate to
$
{\bf 1}_ES_t{\bf 1}_{\Omega^c}
$
and use the global \(L^2\)-boundedness of \(R_t\).
 Consequently,
\[
\left\|
{\bf 1}_ES_tR_t{\bf 1}_F
\right\|_{L^2(\rr^d)\to L^2(\rr^d)}
\lesssim
\exp\left(
-c'\frac{[\operatorname{dist}(E,F)]^2}{t^2}
\right).
\]
If \(\operatorname{dist}(E,F)=0\), the estimate follows from the global
\(L^2\)-boundedness. Repeating the above argument gives the same Gaussian
off-diagonal estimate for any finite composition of such operators, including
the case where the leftmost operator is vector-valued.

For \(\ell=0\), the estimate \eqref{eq:higher-gradient-Gaffney} follows
directly from \eqref{eq:gradient-Gaffney} with \(s=t^2\). Let now
\(\ell\geq1\). Since \(L\) commutes with its semigroup, the spectral
calculus gives
\[
\begin{aligned}
t\nabla\lf(t^2L\r)^\ell e^{-t^2L}
&=
(2\ell)^\ell
\left[
t\nabla e^{-\frac{t^2}{2}L}
\right]
\left[
\frac{t^2}{2\ell}L
e^{-\frac{t^2}{2\ell}L}
\right]^\ell
=(2\ell)^\ell
\sqrt{2}
\left[
\left(\frac{t^2}{2}\right)^{1/2}
\nabla e^{-\frac{t^2}{2}L}
\right]
\left[
\frac{t^2}{2\ell}L
e^{-\frac{t^2}{2\ell}L}
\right]^\ell.
\end{aligned}
\]
The first factor satisfies a gradient Gaffney estimate at scale
\(t\) from \eqref{eq:gradient-Gaffney}. 
The second factor satisfies an \(L^2\) Gaffney estimate
at scale \(t/\sqrt{2\ell}\), and hence at a scale comparable to
\(t\) from \eqref{eq:L-Gaffney}. Therefore, the composition property implies \eqref{eq:higher-gradient-Gaffney}. The constants
\(C_\ell,c_\ell>0\) may depend on \(\ell\), but are independent of
\(E\), \(F\), and \(t\).

Now we prove the first estimate. Since
\[
e^{-t^2T_\kappa}
=
e^{-\kappa t^2}e^{-t^2L},
\]
we have
\[
t\nabla(t^2L)^ke^{-t^2T_\kappa}
=
e^{-\kappa t^2}
t\nabla(t^2L)^ke^{-t^2L}.
\]
It follows from \eqref{eq:higher-gradient-Gaffney} that
\[
\begin{aligned}
\left\|
{\bf 1}_E\,t\nabla(t^2L)^k
e^{-t^2T_\kappa}{\bf 1}_F
\right\|_{L^2(\mathbb R^d)
\to L^2(\mathbb R^d;\mathbb C^d)}
\leq
C_k e^{-\kappa t^2}
\exp\left(
-c_k\frac{[\operatorname{dist}(E,F)]^2}{t^2}
\right).
\end{aligned}
\]
Since, for every \(N>0\),
\[
e^{-c_ku^2}
\leq
C_{N,k}(1+u)^{-N},
\qquad u\geq0,
\]
we obtain
\[
\left\|
{\bf 1}_E\,t\nabla(t^2L)^k
e^{-t^2T_\kappa}{\bf 1}_F
\right\|_{L^2(\mathbb R^d)
\to L^2(\mathbb R^d;\mathbb C^d)}
\leq
C_{N,k}e^{-\kappa t^2}
\left(
1+\frac {\operatorname{dist}(E,F)}{t
}\right)^{-N}.
\]

It remains to prove the second estimate. By the binomial formula and
the identity
\(e^{-t^2T_\kappa}=e^{-\kappa t^2}e^{-t^2L}\), we obtain
\[
\begin{aligned}
t\nabla(t^2T_\kappa)^ke^{-t^2T_\kappa}
=
e^{-\kappa t^2}
t\nabla
\left(t^2L+\kappa t^2I\right)^k
e^{-t^2L}
=
e^{-\kappa t^2}
\sum_{\ell=0}^k
\binom{k}{\ell}
(\kappa t^2)^{k-\ell}
t\nabla(t^2L)^\ell e^{-t^2L}.
\end{aligned}
\]
Let
$
c_k^*:=\min_{0\leq\ell\leq k}c_\ell>0.
$
Using \eqref{eq:higher-gradient-Gaffney}, we obtain
\[
\begin{aligned}
\left\|
{\bf 1}_E\,t\nabla(t^2T_\kappa)^k
e^{-t^2T_\kappa}{\bf 1}_F
\right\|_{L^2(\mathbb R^d)
\to L^2(\mathbb R^d;\mathbb C^d)}
&\leq
e^{-\kappa t^2}
\sum_{\ell=0}^k
\binom{k}{\ell}
(\kappa t^2)^{k-\ell}
C_\ell
\exp\left(
-c_\ell\frac{[\operatorname{dist}(E,F)]^2}{t^2}
\right)
\\
&\leq
C_k e^{-\kappa t^2}
\left[
\sum_{j=0}^k(\kappa t^2)^j
\right]
\exp\left(
-c_k^*\frac{[\operatorname{dist}(E,F)]^2}{t^2}
\right).
\end{aligned}
\]
For every \(j\in\{0,\ldots,k\}\),
\[
x^je^{-x}
\leq
C_j e^{-x/2},
\qquad x\geq0.
\]
Consequently,
\[
e^{-\kappa t^2}
\sum_{j=0}^k\lf(\kappa t^2\r)^j
\leq
C_k e^{-\kappa t^2/2},
\]
Combining this with
\[
\exp\left(
-c_k^*\frac{[\operatorname{dist}(E,F)]^2}{t^2}
\right)
\leq
C_{N,k}
\left(
1+\frac {\operatorname{dist}(E,F)}{t
}\right)^{-N},
\]
we conclude that
\[
\left\|
{\bf 1}_E\,t\nabla(t^2T_\kappa)^k
e^{-t^2T_\kappa}{\bf 1}_F
\right\|_{L^2(\mathbb R^d)
\to L^2(\mathbb R^d;\mathbb C^d)}
\leq
C_{N,k}e^{-\kappa t^2/2}
\left(
1+\frac {\operatorname{dist}(E,F)}{t}
\right)^{-N}.
\]
All the constants above depend at most on \(d\), the ellipticity
constants of \(A\), \(N\), and \(k\), and are independent of
\(\kappa\), \(E\), \(F\), and \(t\). This finishes the proof.
\end{proof}

\begin{proof}[Proof of Theorem \ref{thm:local-Riesz}]
Firstly, by ellipticity and the form characterization of \(T_\kappa\), we obtain that, for any $g\in L^2(\mathbb R^d)$,
\[
\left\|
\nabla T_\kappa^{-1/2}g
\right\|_{L^2(\mathbb R^d;\mathbb C^d)}
\lesssim
\|g\|_{L^2(\mathbb R^d)}.
\]
Moreover, by the first estimate in Lemma \ref{grad}, we know that, for every
\(k\in\mathbb Z_+\), \(N>0\),
\(E,F\subset\mathbb R^d\), and \(t>0\),
\begin{equation}\label{eq:Riesz-gradient-OD}
\left\|
{\bf1}_E\,t\nabla(t^2L)^k
e^{-t^2T_\kappa}{\bf1}_F
\right\|_{L^2(\rr^d)\to L^2(\rr^d;\mathbb{C}^d)}
\lesssim
e^{-\kappa t^2}
\left(
1+\frac{\operatorname{dist}(E,F)}t
\right)^{-N}.
\end{equation}

We also need the fixed-scale off-diagonal estimate
\begin{equation}\label{eq:Riesz-fixed-OD}
\left\|
{\bf1}_E\mathcal R_{L,\kappa}^{\rm loc}{\bf1}_F
\right\|_{L^2(\mathbb R^d)
\longrightarrow
L^2(\mathbb R^d;\mathbb C^d)}
\lesssim_{N,\kappa}
\left[
1+\operatorname{dist}(E,F)
\right]^{-N}.
\end{equation}
To justify it rigorously, use
\[
\mathcal R_{L,\kappa}^{\rm loc}
=
\frac2{\sqrt\pi}
\lim_{a\downarrow0,\,b\uparrow\infty}
\int_a^b
\nabla e^{-t^2T_\kappa}\,dt
\]
in the strong \(L^2\)-sense. Indeed, the truncated integral equals
\[
\mathcal R_{L,\kappa}^{\rm loc}
\psi_{a,b}(T_\kappa),
\qquad
\psi_{a,b}(\lambda)
:=
\lambda^{1/2}
\int_a^b e^{-t^2\lambda}\,dt,
\]
and the multipliers \(\psi_{a,b}\) are uniformly bounded and converge
pointwise to \(\sqrt\pi/2\). Thus, when
\(D:=\operatorname{dist}(E,F)>1\), \eqref{eq:Riesz-gradient-OD} gives
\[
\left\|
{\bf1}_E\mathcal R_{L,\kappa}^{\rm loc}{\bf1}_F
\right\|_{L^2(\mathbb R^d)
\longrightarrow
L^2(\mathbb R^d;\mathbb C^d)}
\lesssim
\int_0^\infty
t^{-1}e^{-\kappa t^2}
\left(1+\frac Dt\right)^{-K}\,dt
\lesssim_{N,\kappa}
D^{-N},
\]
provided \(K>N\) is sufficiently large. When \(D\leq1\), the result
follows from the \(L^2\)-boundedness.

Put
\[
\delta_0
:=
d\left(\frac1{s_0}-\frac1{q_0}\right),
\qquad
a_0
:=
d\left(\frac1{s_0}-\frac12\right).
\]
Choose \(M\in\mathbb N\) and \(\varepsilon>0\) such that
\[
2M>a_0+\delta_0
=
d\left(
\frac2{s_0}-\frac12-\frac1{q_0}
\right)
\ \
\mathrm{and}
 \ \
\varepsilon>\delta_0.
\]
Then choose
$
\eta\in
\left(
\delta_0,\,
\min\{\varepsilon,2M-a_0\}
\right).
$
Since the heat semigroup generated by \(L\) satisfies
Assumption \ref{a2.5}. Hence Theorem \ref{mol} is applicable.
Let \(f\in\mathbf h_L^X(\mathbb R^d)\). By Theorem \ref{mol}, there
exist local \((X,2,M,\varepsilon)\)-molecules
\(\{m_j\}_{j\in\mathbb N}\), associated respectively with balls
\(\{B_j\}_{j\in\mathbb N}\), and coefficients
\(\{\lambda_j\}_{j\in\mathbb N}\subset[0,\infty)\) such that
\[
f=\sum_{j\in\mathbb N}\lambda_jm_j
\ \ \text{in } \ \ L^2(\mathbb R^d)
\]
and
\begin{equation}\label{eq:Riesz-local-coefficients}
\left\|
\left\{
\sum_{j\in\mathbb N}
\left[
\frac{\lambda_j{\bf1}_{B_j}}
{\|{\bf1}_{B_j}\|_{X(\mathbb R^d)}}
\right]^{s_0}
\right\}^{1/s_0}
\right\|_{X(\mathbb R^d)}
\lesssim
\|f\|_{h_L^X(\mathbb R^d)}.
\end{equation}
For a ball \(B\), set
\[
\alpha_i(B)
:=
\lf|2^iB\r|^{1/2}
\lf\|{\bf1}_{2^iB}\r\|_{X(\mathbb R^d)}^{-1},
\qquad i\in\mathbb Z_+.
\]
By \eqref{e:ind-comp-lem39-short}, we know that, for all
\(i,\nu\in\mathbb Z_+\),
\begin{equation}\label{eq:Riesz-alpha-local}
\alpha_\nu(B)
\lesssim
\begin{cases}
2^{(i-\nu)a_0}\alpha_i(B),&\nu<i,\\
\alpha_i(B),&\nu\geq i.
\end{cases}
\end{equation}

Now we claim that every local \((X,2,M,\varepsilon)\)-molecule \(m\),
associated with \(B\), satisfies
\begin{equation}\label{eq:Riesz-local-molecule}
\left\|
\mathcal R_{L,\kappa}^{\rm loc}m
\right\|_{L^2(S_i(B);\mathbb C^d)}
\lesssim
2^{-i\eta}\alpha_i(B),
\qquad i\in\mathbb Z_+.
\end{equation}

Assume first that \(m\) is of type \({\rm(I)}\). Then \(r_B\geq1\).
Writing
\[
m_\nu:={\bf1}_{S_\nu(B)}m
\]
and using \eqref{eq:Riesz-fixed-OD}, we obtain
\[
\begin{aligned}
\left\|
{\bf1}_{S_i(B)}
\mathcal R_{L,\kappa}^{\rm loc}m
\right\|_{L^2(\rr^d)}
&\lesssim
\sum_{|\nu-i|\leq1}\|m_\nu\|_{L^2(\rr^d)}
+
2^{-iN}\sum_{\nu=0}^{i-2}\|m_\nu\|_{L^2(\rr^d)}
+
\sum_{\nu=i+2}^\infty
2^{-\nu N}\|m_\nu\|_{L^2(\rr^d)}.
\end{aligned}
\]
Here and below, a sum is understood to be zero whenever its upper limit is smaller than its lower limit. By the molecular estimates and
\eqref{eq:Riesz-alpha-local}, we obtain
\[
\begin{aligned}
\left\|
{\bf1}_{S_i(B)}
\mathcal R_{L,\kappa}^{\rm loc}m
\right\|_{L^2(\rr^d)}
&\lesssim
\left[
2^{-i\varepsilon}
+
2^{-i(N-a_0)}
+
2^{-i(N+\varepsilon)}
\right]\alpha_i(B)\lesssim
2^{-i\eta}\alpha_i(B)
\end{aligned}
\]
when \(N>a_0+\eta\).

Next we assume that \(m\) is of type \({\rm(II)}\). Then
\(r:=r_B\in(0,1)\), and there exists \(b\in D(L^M)\) such that
$
m=L^Mb.
$
Define
\[
\phi_r(\lambda)
:=
\lambda^{1/2}\int_0^r e^{-t^2\lambda}\,dt
 \ \
\mathrm{and}
 \ \
U_r
:=
\frac2{\sqrt\pi}
\mathcal R_{L,\kappa}^{\rm loc}\phi_r(T_\kappa).
\]
Then, for every \(g\in L^2(\mathbb R^d)\),
\[
U_rg
=
\frac2{\sqrt\pi}
\lim_{\rho\downarrow0}
\int_\rho^r
\nabla e^{-t^2T_\kappa}g\,dt
\ \ \text{in } \ \
L^2(\mathbb R^d;\mathbb C^d).
\]
For every \(g\in L^2(\mathbb R^d)\), define
\[
W_rg
:=
\frac2{\sqrt\pi}
\int_r^\infty
t^{-2M-1}
\left[
t\nabla(t^2L)^Me^{-t^2T_\kappa}
\right]g\,dt.
\]
Moreover,
\[
\begin{aligned}
\int_r^\infty
t^{-2M-1}
\left\|
t\nabla(t^2L)^Me^{-t^2T_\kappa}g
\right\|_{L^2(\mathbb R^d;\mathbb C^d)}\,dt\lesssim
\|g\|_{L^2(\mathbb R^d)}
\int_r^\infty
t^{-2M-1}e^{-\kappa t^2}\,dt
<\infty.
\end{aligned}
\]
Consequently,
\[
\mathcal R_{L,\kappa}^{\rm loc}m
=
U_rm+W_rb.
\]
It follows from \eqref{eq:Riesz-gradient-OD} that
\[
\lf\|{\bf1}_EU_r{\bf1}_F\r\|_{L^2(\mathbb R^d)
\longrightarrow
L^2(\mathbb R^d;\mathbb C^d)}
\lesssim_N
\left(
1+\frac{\operatorname{dist}(E,F)}r
\right)^{-N}
\]
and, for \(N>2M\),
\[
\|{\bf1}_EW_r{\bf1}_F\|_{L^2(\mathbb R^d)
\longrightarrow
L^2(\mathbb R^d;\mathbb C^d)}
\lesssim
\left[
r+\operatorname{dist}(E,F)
\right]^{-2M}.
\]
Consequently, the same annular summation gives
\[
\left\|
{\bf1}_{S_i(B)}U_rm
\right\|_{L^2(\rr^d)}
\lesssim
\left[
2^{-i\varepsilon}
+
2^{-i(N-a_0)}
+
2^{-i(N+\varepsilon)}
\right]\alpha_i(B).
\]

Writing \(b_\nu:={\bf1}_{S_\nu(B)}b\), the molecular condition gives
\[
\|b_\nu\|_{L^2(\rr^d)}
\lesssim
2^{-\nu\varepsilon}r^{2M}\alpha_\nu(B).
\]
Therefore,
\[
\begin{aligned}
\left\|
{\bf1}_{S_i(B)}W_rb
\right\|_{L^2(\rr^d)}
&\lesssim
\sum_{|\nu-i|\leq1}
r^{-2M}\|b_\nu\|_{L^2(\rr^d)}
+
\lf(2^ir\r)^{-2M}
\sum_{\nu=0}^{i-2}\|b_\nu\|_{L^2(\rr^d)}
+
\sum_{\nu=i+2}^\infty
\lf(2^\nu r\r)^{-2M}\|b_\nu\|_{L^2(\rr^d)}\\
&\lesssim
\left[
2^{-i\varepsilon}
+
2^{-i(2M-a_0)}
+
2^{-i(2M+\varepsilon)}
\right]\alpha_i(B)\\
&\lesssim
2^{-i\eta}\alpha_i(B).
\end{aligned}
\]
This proves \eqref{eq:Riesz-local-molecule} for both types of molecules.

Finally, we define the non-negative sublinear operator
\[
\mathscr R_\kappa g(x)
:=
\left|
\mathcal R_{L,\kappa}^{\rm loc}g(x)
\right|_{\mathbb C^d}.
\]
It is bounded on \(L^2(\mathbb R^d)\), and
\eqref{eq:Riesz-local-molecule} is exactly the annular estimate required
by Lemma \ref{l2.8}, since
$
\eta>
d(\frac1{s_0}-\frac1{q_0}).
$
Applying Lemma \ref{l2.8} and
\eqref{eq:Riesz-local-coefficients}, we obtain
\[
\left\|
\mathcal R_{L,\kappa}^{\rm loc}f
\right\|_{X(\mathbb R^d;\mathbb C^d)}
\lesssim_\kappa
\|f\|_{h_L^X(\mathbb R^d)}.
\]
The desired extension follows from the density of
\(\mathbf h_L^X(\mathbb R^d)\) in \(h_L^X(\mathbb R^d)\) and the
completeness of \(X(\mathbb R^d;\mathbb C^d)\).
\end{proof}

\subsection{Boundedness of H\"ormander-type spectral multipliers}\label{s62}

In this subsection, we establish a H\"ormander-type spectral multiplier
theorem on the local Hardy space \(h_L^X(\bx)\). Put
\[
        T:=L+I .
\]
Then \(T\) is a non-negative self-adjoint operator on \(L^2(\bx)\), and
\[
        e^{-tT}=e^{-t}e^{-tL},\qquad t>0.
\]
Thus \(T\) also satisfies Assumption \(\ref{a2.5}\). Moreover,
\[
        \sigma(T)\subset[1,\infty).
\]

Let \(\psi\in C_{\mathrm c}^{\infty}((1/2,2))\) be a non-trivial cut-off
function such that
\[
        \sum_{j\in\mathbb Z}\psi(2^{-j}\lambda)=1,
        \qquad \lambda\in(0,\infty).
\]
For \(\beta>0\) and a bounded Borel function \(F\) on \([0,\infty)\), define
\[
 \|F\|_{\mathcal H^\beta}
 :=
 \sup_{r>0}
 \left\|
 \psi(\cdot)F(r\cdot)
 \right\|_{W^{\beta,\infty}(\mathbb R)}.
\]
Here \(\psi(\cdot)F(r\cdot)\), initially defined on \((0,\infty)\), is
extended to \(\mathbb R\) by zero.

We shall also use the following global H\"ormander multiplier theorem for
\(H_T^X(\bx)\). By the H\"ormander-type spectral multiplier theorem for Hardy spaces
associated with ball quasi-Banach function spaces, see
Lin et al. \cite[Theorem 4.1]{lyyy23}, applied to
\(T=L+I\), we obtain that \(F(T)\) is bounded on \(H_T^X(\mathbb X)\).

\begin{lemma}\label{lem:global-Hormander}
Let \(T\) be a non-negative self-adjoint operator on \(L^2(\bx)\) satisfying
Assumption \(\ref{a2.5}\). Let \(X(\bx)\) be a BQBF space satisfying
Assumption \(\ref{a11}\) for some \(p\in(0,\infty)\) and Assumption
\(\ref{a12}\) for some
\[
        s_0\in(0,\min\{p,1\}]
        \quad\text{and}\quad
        q_0\in(s_0,2].
\]
Then there exists \(\beta_0>0\) such that, whenever
\(\beta>\beta_0\), every bounded Borel function \(F\) on \([0,\infty)\)
satisfying
\[
        \|F\|_{L^\infty([0,\infty))}+\|F\|_{\mathcal H^\beta}<\infty
\]
gives rise to a bounded operator \(F(T)\) on \(H_T^X(\bx)\). Moreover,
\[
 \|F(T)g\|_{H_T^X(\bx)}
 \leq
 C
 \left[
 \|F\|_{L^\infty([0,\infty))}+\|F\|_{\mathcal H^\beta}
 \right]
 \|g\|_{H_T^X(\bx)}
\]
for all \(g\in H_T^X(\bx)\).
\end{lemma}

\begin{theorem}
\label{thm:Hormander-local}
Let \(L\) be a non-negative self-adjoint operator on \(L^2(\bx)\) satisfying
Assumption \(\ref{a2.5}\). Let \(X(\bx)\) be a BQBF space satisfying
Assumption \(\ref{a11}\) for some \(p\in(0,\infty)\) and Assumption
\(\ref{a12}\) for some
\[
        s_0\in(0,\min\{p,1\}]
        \quad\text{and}\quad
        q_0\in(s_0,2].
\]
Put \(T:=L+I\). Then there exists \(\beta_0>0\) such that, for every
\(\beta>\beta_0\) and every bounded Borel function \(F\) on
\([0,\infty)\) satisfying
\[
        \|F\|_{L^\infty([0,\infty))}+\|F\|_{\mathcal H^\beta}<\infty,
\]
the operator \(F(L+I)\), initially defined on \(L^2(\bx)\) by the spectral
theorem, extends uniquely to a bounded operator on \(h_L^X(\bx)\). Moreover,
\[
 \|F(L+I)f\|_{h_L^X(\bx)}
 \leq
 C
 \left[
 \|F\|_{L^\infty([0,\infty))}+\|F\|_{\mathcal H^\beta}
 \right]
 \|f\|_{h_L^X(\bx)}
\]
for all \(f\in h_L^X(\bx)\).
\end{theorem}

\begin{proof}[Proof of Theorem \ref{thm:Hormander-local}]
Let
$
        T:=L+I.
$
As observed above, \(T\) is a non-negative self-adjoint operator satisfying
Assumption \(\ref{a2.5}\), and
$
        \sigma(T)\subset[1,\infty).
$
By Lemma \(\ref{lem:global-Hormander}\), the H\"{o}rmander condition on \(F\)
implies
\[
 \|F(T)g\|_{H_T^X(\bx)}
 \leq
 C
 \left[
 \|F\|_{L^\infty([0,\infty))}+\|F\|_{\mathcal H^\beta}
 \right]
 \|g\|_{H_T^X(\bx)}
\]
for all \(g\in H_T^X(\bx)\).

On the other hand, Theorem \(\ref{eq1}\) gives the algebraic and topological
identity
\[
        h_L^X(\bx)=H_{L+I}^X(\bx)=H_T^X(\bx).
\]
Hence, for every \(f\in\mathbf h_L^X(\bx)\),
\[
\begin{aligned}
 \|F(L+I)f\|_{h_L^X(\bx)}
 &\sim
 \|F(T)f\|_{H_T^X(\bx)}
 \lesssim
 \left[\|F\|_{L^\infty([0,\infty))}+\|F\|_{\mathcal H^\beta}\right]\|f\|_{H_T^X(\bx)}\\
 &\sim
 \left[
 \|F\|_{L^\infty([0,\infty))}+\|F\|_{\mathcal H^\beta}
 \right]
 \|f\|_{h_L^X(\bx)}.
\end{aligned}
\]
Since \(\mathbf h_L^X(\bx)\) is dense in \(h_L^X(\bx)\) by definition, the
operator \(F(L+I)\) admits a unique bounded extension to \(h_L^X(\bx)\).
This proves the theorem.
\end{proof}

\section{Applications to specific Hardy-type function spaces}\label{s6}
In this section, we recall several concrete examples of ball quasi-Banach
function spaces and describe the corresponding local Hardy spaces. These
examples include local Orlicz-Hardy spaces (see Subsection \ref{OHS}), local
variable Hardy spaces (see Subsection \ref{VHS}), and local mixed-norm Hardy
spaces (see Subsection \ref{MHS}). They illustrate the applicability of the
abstract results obtained in the preceding sections to various concrete
function space settings.

\subsection{Local Orlicz-Hardy spaces}\label{OHS}

Let us recall the definitions of Orlicz functions and Orlicz spaces. A
function
\[
 \Phi:[0,\infty)\longrightarrow [0,\infty)
\]
is called an Orlicz function if \(\Phi\) is non-decreasing,
\(\Phi(0)=0\), \(\Phi(t)>0\) for all \(t\in(0,\infty)\), and
\[
 \lim_{t\to\infty}\Phi(t)=\infty.
\]
Moreover, an Orlicz function \(\Phi\) is said to be of lower type \(r\)
for some \(r\in\mathbb R\) if there exists a positive constant \(C_r\)
such that, for all \(t\in[0,\infty)\) and \(s\in(0,1)\),
\[
 \Phi(st)\leq C_rs^r\Phi(t).
\]
Similarly, \(\Phi\) is said to be of upper type \(r\) if there exists a
positive constant \(C_r\) such that, for all \(t\in[0,\infty)\) and
\(s\in[1,\infty)\),
\[
 \Phi(st)\leq C_rs^r\Phi(t).
\]
Let \(\Phi\) be an Orlicz function with positive lower type
\(r_\Phi^{-}\) and positive upper type \(r_\Phi^{+}\). The Orlicz space
\(L^\Phi(\bx)\) is defined to be the collection of all \(\mu\)-measurable
functions \(f\) on \(\bx\) such that
\[
 \|f\|_{L^\Phi(\bx)}
 :=
 \inf\left\{
 \lambda\in(0,\infty):
 \int_{\bx}
 \Phi\left(\frac{|f(x)|}{\lambda}\right)\,d\mu(x)
 \leq1
 \right\}
 <\infty.
\]
Then \(L^\Phi(\bx)\) is a quasi-Banach function space and, in particular, a
ball quasi-Banach function space; see, for instance, \cite[Subsection
8.3]{syy22}.

In what follows, taking
\[
 X(\bx):=L^\Phi(\bx),
\]
we write
\[
 h_L^\Phi(\bx):=h_L^X(\bx),\quad
 h^\Phi_{L,\mathrm{mol},M,\varepsilon}(\bx)
 :=
 h^X_{L,\mathrm{mol},M,\varepsilon}(\bx),
\quad {\rm{and}} \quad
 H^\Phi_{L+I}(\bx):=H^X_{L+I}(\bx).
\]
Applying Theorems \ref{mol} and \ref{eq1} with
\(X(\bx)=L^\Phi(\bx)\), we obtain the following conclusions.

\begin{theorem}\label{mol-OHS}
Let \(L\) be a non-negative self-adjoint operator on \(L^2(\bx)\) satisfying
Assumption \ref{a2.5}. Let \(\Phi\) be an Orlicz function with positive lower
type \(r_\Phi^{-}\) and positive upper type \(r_\Phi^{+}\).
Assume that
 $0<r_\Phi^{-}\leq 1$ and $r_\Phi^{-}\leq r_\Phi^{+}< 2$.
Then the following assertions hold.

\begin{enumerate}
\item[{\rm(i)}]
If $\varepsilon>d(\frac1{r_\Phi^{-}}-\frac1{2}) $ and $2M>d(\frac2{r_\Phi^{-}}-1),$
then
\[
 h^\Phi_{L,\mathrm{mol},M,\varepsilon}(\bx)
 \subseteq
 h^\Phi_L(\bx).
\]

\item[{\rm(ii)}]
If \(\varepsilon>0\) and \(M\in\nn\), then
\[
 h^\Phi_L(\bx)
 \subseteq
 h^\Phi_{L,\mathrm{mol},M,\varepsilon}(\bx).
\]
\end{enumerate}
\end{theorem}

\begin{proof}[Proof of Theorem \ref{mol-OHS}]

To show Theorem \ref{mol-OHS}, we only need  to show that the Orlicz space $L^{\Phi}$ satisfies Assumptions \ref{a11} and \ref{a12}.
Take $p:=r_{\Phi}^-$, $q_0:=2$, and $s_0\in(0,\min\{1,r^-_{\Phi}\}]$ such that $2M>
 d(
 \frac2{s_0}-\frac12-\frac1{q_0}
 )$ and $\varepsilon>
 d(\frac1{s_0}-\frac1{q_0})$.
By \cite{shyy17}, for the above $p$, we deduce that the BQBF space $X(\bx)=L^{\Phi}(\bx)$ satisfies Assumption \ref{a11}.

By \cite{shyy17}, we obtain that
$$\lf[\lf(\lf[L^{\Phi}(\bx)\r]^{\frac{1}{s_0}}\r)^{'}\r]^{\frac{1}{(2/s_0)^{'}}}=L^{\Psi}(\bx),$$
where, for any $t\in(0,\fz)$,
$$\Psi(t):=\sup_{u\in(0,\fz)}\lf(t^{\frac{1}{(2/s_0)^{'}}}u-\Phi(u^{\frac{1}{s_0}})\r).$$
Moreover, by \cite{shyy17}, we find that $\Psi$ is an Orlicz function with positive lower type $r^-_{\Psi}:=\frac{(r_{\Phi}^+/s_0)^{'}}{(2/s_0)^{'}}\in(1,\fz)$.
By the above facts and \cite{shyy17}, we conclude that the BQBF space $X(\bx)=L^{\Phi}(\bx)$ satisfies Assumption \ref{a12}.
This completes the proof Theorem \ref{mol-OHS}.
\end{proof}

\begin{theorem}\label{eq0-OHS}
Let \(L\) be a non-negative self-adjoint operator on \(L^2(\bx)\) satisfying
Assumption \ref{a2.5}. Let \(\Phi\) be an Orlicz function with positive lower
type \(r_\Phi^{-}\) and positive upper type \(r_\Phi^{+}\).
Assume that $0<r_\Phi^{-}\leq 1$ and $ r_\Phi^{-}\leq r_\Phi^{+}< 2.$
If $\inf\sigma(L)>0,$ then
$$h^{\Phi}_L(\bx)=H^{\Phi}_L(\bx)$$
algebraically and topologically, where \(\sigma(L)\) denotes the spectrum of
\(L\) on \(L^2(\bx)\).
\end{theorem}

\begin{theorem}\label{eq1-OHS}
Let \(L\) be a non-negative self-adjoint operator on \(L^2(\bx)\) satisfying
Assumption \ref{a2.5}. Let \(\Phi\) be an Orlicz function with positive lower
type \(r_\Phi^{-}\) and positive upper type \(r_\Phi^{+}\). Assume that
$0<r_\Phi^{-}\leq 1$ and $ r_\Phi^{-}\leq r_\Phi^{+}< 2.$ Then
$$h^\Phi_L(\bx)=H^\Phi_{L+I}(\bx)$$
algebraically and topologically.
\end{theorem}

\begin{theorem}\label{Riesz-OHS}
Let \(\Phi\) be an Orlicz function with positive lower
type \(r_\Phi^{-}\) and positive upper type \(r_\Phi^{+}\). Assume that
$0<r_\Phi^{-}\leq 1$ and $ r_\Phi^{-}\leq r_\Phi^{+}< 2.$
Let
\(L:=-\operatorname{div}(A\nabla)\) be as above, and let \(T:=L+\kappa I\) with $\kappa>0$.
Then the local Riesz transform
$
        \mathcal R_{L,\kappa}^{\rm loc}
$
admits a bounded extension from \(h_L^{\Phi}(\rr^d)\) to
\(L^{\Phi}(\rr^d;\mathbb C^d)\). More precisely, there exists a positive constant
\(C\) such that, for all \(f\in h_L^{\Phi}(\rr^d)\cap L^2(\rr^d)\),
\[
 \left\|
 \mathcal R_{L,\kappa}^{\rm loc}f
 \right\|_{L^{\Phi}(\rr^d;\mathbb C^d)}
 \leq
 C\|f\|_{h_L^{\Phi}(\rr^d)}.
\]
\end{theorem}

\begin{theorem}\label{hor-OHS}Let \(L\) be a non-negative self-adjoint operator on \(L^2(\rr^d)\) satisfying
Assumption \(\ref{a2.5}\), and
 \(\Phi\) be an Orlicz function with positive lower
type \(r_\Phi^{-}\) and positive upper type \(r_\Phi^{+}\). Assume that
$0<r_\Phi^{-}\leq 1$ and $ r_\Phi^{-}\leq r_\Phi^{+}< 2.$
Put \(T:=L+I\). Then there exists \(\beta_0>0\) such that, for every
\(\beta>\beta_0\) and every bounded Borel function \(F\) on
\([0,\infty)\) satisfying
\[
        \|F\|_{L^\infty([0,\infty))}+\|F\|_{\mathcal H^\beta}<\infty,
\]
the operator \(F(L+I)\), initially defined on \(L^2(\rr^d)\) by the spectral
theorem, extends uniquely to a bounded operator on \(h_L^{\Phi}(\rr^d)\). Moreover,
\[
 \|F(L+I)f\|_{h_L^{\Phi}(\rr^d)}
 \leq
 C
 \left[
 \|F\|_{L^\infty([0,\infty))}+\|F\|_{\mathcal H^\beta}
 \right]
 \|f\|_{h_L^{\Phi}(\rr^d)}
\]
for all \(f\in h_L^{\Phi}(\rr^d)\).
\end{theorem}
The proofs of Theorems \ref{eq0-OHS}, \ref{eq1-OHS}, \ref{Riesz-OHS}, and \ref{hor-OHS} are similar to that of Theorem \ref{mol-OHS}; therefore, we omit the details.

\subsection{Local variable Hardy spaces}\label{VHS}

Let \(\mathcal P(\bx)\) denote the set of all measurable functions
\(p(\cdot):\bx\to(0,\infty)\) such that
$0<p_- \leq p_+<\infty,$
where
\[
 p_-:=\mathop{\mathrm{ess\,inf}}_{x\in\bx}p(x)
 \quad\hbox{and}\quad
 p_+:=\mathop{\mathrm{ess\,sup}}_{x\in\bx}p(x).
\]
Let \(\mathcal B(\bx)\) be the collection of all
\(p(\cdot)\in\mathcal P(\bx)\) for which the Hardy--Littlewood maximal
operator \(M_{\mathrm{HL}}\) is bounded on \(L^{p(\cdot)}(\bx)\). It is well
known that, if \(p(\cdot)\in\mathcal P(\bx)\) satisfies
\(1<p_-\leq p_+<\infty\) and the global log-H\"older continuity condition,
then \(p(\cdot)\in\mathcal B(\bx)\).

For any measurable function \(p(\cdot):\bx\to(0,\infty)\), the variable
Lebesgue space \(L^{p(\cdot)}(\bx)\) is defined to be the set of all
\(\mu\)-measurable functions \(f\) on \(\bx\) such that
\[
\|f\|_{L^{p(\cdot)}(\bx)}
:=
\inf\left\{
\lambda\in(0,\infty):
\int_{\bx}
\left(\frac{|f(x)|}{\lambda}\right)^{p(x)}
\,d\mu(x)
\leq1
\right\}
<\infty.
\]
By \cite[Remark 2.7]{yhyy23}, \(L^{p(\cdot)}(\bx)\) is a ball
quasi-Banach function space whenever
$0<p_-\leq p_+<\infty,$
and it is a ball Banach function space whenever
$1<p_-\leq p_+<\infty.$

We recall the log-H\"older continuity conditions. The variable exponent
\(p(\cdot)\) is said to be locally log-H\"older continuous if there exists a
positive constant \(C_{\log}(p)\) such that, for all \(x,y\in\bx\),
\[
 |p(x)-p(y)|
 \leq
 \frac{C_{\log}(p)}
      {\log\left(e+\frac1{d(x,y)}\right)}.
\]
Moreover, \(p(\cdot)\) is said to satisfy the log-H\"older decay condition
with respect to a base point \(x_\infty\in\bx\) if there exist
\(p_\infty\in\mathbb R\) and a positive constant \(C_\infty\) such that, for
all \(x\in\bx\),
\[
 |p(x)-p_\infty|
 \leq
 \frac{C_\infty}
      {\log\left(e+d(x,x_\infty)\right)}.
\]
The exponent \(p(\cdot)\) is called globally log-H\"older continuous if it
satisfies both the local log-H\"older continuity condition and the
log-H\"older decay condition.

Let \(p(\cdot):\bx\to(0,\infty)\) be a globally log-H\"older continuous
function, and let
\[
 X(\bx):=L^{p(\cdot)}(\bx).
\]
In this case, we write
$$h_L^{p(\cdot)}(\bx):=h_L^X(\bx),\quad
 h^{p(\cdot)}_{L,\mathrm{mol},M,\varepsilon}(\bx):=h^X_{L,\mathrm{mol},M,\varepsilon}(\bx),
\quad {\rm{and}} \quad
H^{p(\cdot)}_{L+I}(\bx):=H^X_{L+I}(\bx).$$

Applying Theorems \ref{mol}, \ref{eq0}, and \ref{eq1} with
\(X(\bx)=L^{p(\cdot)}(\bx)\), we obtain the following conclusions.

\begin{theorem}\label{mol-VHS}
Let \(L\) be a non-negative self-adjoint operator on \(L^2(\bx)\) satisfying
Assumption \ref{a2.5}. Let \(p(\cdot):\bx\to(0,\infty)\) be a globally
log-H\"older continuous function such that
$0<p_{-}\leq 1$ and $p_{-}\leq p_{+}< 2.$
 Then the following assertions hold.

\begin{enumerate}
\item[{\rm(i)}]
If $\varepsilon>d(\frac1{p_-}-\frac1{p_{+}})$ and $2M>d(\frac2{p_-}-\frac1{2}-\frac1{p_+}),$
then
\[
 h^{p(\cdot)}_{L,\mathrm{mol},M,\varepsilon}(\bx)
 \subseteq
 h^{p(\cdot)}_L(\bx).
\]

\item[{\rm(ii)}]
If \(\varepsilon>0\) and \(M\in\nn\), then
\[
 h^{p(\cdot)}_L(\bx)
 \subseteq
 h^{p(\cdot)}_{L,\mathrm{mol},M,\varepsilon}(\bx).
\]
\end{enumerate}
\end{theorem}
\begin{proof}
Take $p:=p_-$, $q_0=2$, and $s_0\in(0,\,\min\{1,{p}_-\}]$ such that $2M \in(d [\frac{2}{{s}_0}-\frac12-\frac{1}{q_0}],\fz)\cap\mathbb{N}$ and $\varepsilon\in(\frac{d}{s_0}-\frac{d}{q_0},\fz)$.
From \cite{shyy17}, it follows that the BQBF space $X(\bx)=L^{p(\cdot)}(\bx)$ satisfies Assumption \ref{a11} for the above exponents.
Moreover, by \cite{shyy17}, we conclude that
$$\lf[\lf(\lf[L^{p(\cdot)}(\bx)\r]^{\frac{1}{s_0}}\r)^{'}\r]^{\frac{1}{(2/s_0)^{'}}}=
L^{\frac{(p(\cdot)/s_0)^{'}}{(2/s_0)^{'}}}(\bx),$$
where, for any $x\in\bx$, $\frac{1}{(p(x)/s_0)^{'}}+\frac{1}{p(x)/s_0}=1$.
Therefore, from the assumption that $0<{p}_-\leq {p}_+<2$ and \cite{shyy17}, we deduce that the BQBF space $X(\bx)=L^{p(\cdot)}(\bx)$ satisfies Assumption \ref{a12} for above exponents. Therefore, all the assumption of Theorem \ref{mol} are satisfied, this implies the desired result. Therefore, we finish the proof of Theorem \ref{mol-VHS}.
\end{proof}
\begin{theorem}\label{eq0-VHS}
Let \(L\) be a non-negative self-adjoint operator on \(L^2(\bx)\) satisfying
Assumption \ref{a2.5}. Let \(p(\cdot):\bx\to(0,\infty)\) be a globally
log-H\"older continuous function such that
$0<p_{-}\leq 1$ and $p_{-}\leq p_{+}< 2.$
 If
$
 \inf\sigma(L)>0,
$
then
\[
 h^{p(\cdot)}_L(\bx)=H^{p(\cdot)}_L(\bx)
\]
algebraically and topologically, where \(\sigma(L)\) denotes the spectrum of
\(L\) on \(L^2(\bx)\).
\end{theorem}

\begin{theorem}\label{eq1-VHS}
Let \(L\) be a non-negative self-adjoint operator on \(L^2(\bx)\) satisfying
Assumption \ref{a2.5}. Let \(p(\cdot):\bx\to(0,\infty)\) be a globally
log-H\"older continuous function such that
$0<p_{-}\leq 1$ and $p_{-}\leq p_{+}< 2.$
 Then
\[
 h^{p(\cdot)}_L(\bx)=H^{p(\cdot)}_{L+I}(\bx)
\]
algebraically and topologically.
\end{theorem}

\begin{theorem}\label{Riesz-VHS}
Let \(p(\cdot):\rr^d\to(0,\infty)\) be a globally
log-H\"older continuous function such that
$0<p_{-}\leq 1$ and $p_{-}\leq p_{+}< 2.$ Let
\(L:=-\operatorname{div}(A\nabla)\) be as above, and let \(T:=L+\kappa I\) with $\kappa>0$.
Then the local Riesz transform
$
        \mathcal R_{L,\kappa}^{\rm loc}
$
admits a bounded extension from \(h_L^{{p(\cdot)}}(\rr^d)\) to
\(L^{p(\cdot)}(\rr^d;\mathbb C^d)\). More precisely, there exists a positive constant
\(C\) such that, for all \(f\in h_L^{p(\cdot)}(\rr^d)\cap L^2(\rr^d)\),
\[
 \left\|
 \mathcal R_{L,\kappa}^{\rm loc}f
 \right\|_{L^{p(\cdot)}(\rr^d;\mathbb C^d)}
 \leq
 C\|f\|_{h_L^{p(\cdot)}(\rr^d)}.
\]
\end{theorem}

\begin{theorem}
\label{hor-VHS}Let \(L\) be a non-negative self-adjoint operator on \(L^2(\rr^d)\) satisfying
Assumption \(\ref{a2.5}\), and \(p(\cdot):\rr^d\to(0,\infty)\) be a globally
log-H\"older continuous function such that
$0<p_{-}\leq 1$ and $p_{-}\leq p_{+}< 2.$
Put \(T:=L+I\). Then there exists \(\beta_0>0\) such that, for every
\(\beta>\beta_0\) and every bounded Borel function \(F\) on
\([0,\infty)\) satisfying
\[
        \|F\|_{L^\infty([0,\infty))}+\|F\|_{\mathcal H^\beta}<\infty,
\]
the operator \(F(L+I)\), initially defined on \(L^2(\rr^d)\) by the spectral
theorem, extends uniquely to a bounded operator on \(h_L^{p(\cdot)}(\rr^d)\). Moreover,
\[
 \|F(L+I)f\|_{h_L^{p(\cdot)}(\rr^d)}
 \leq
 C
 \left[
 \|F\|_{L^\infty([0,\infty))}+\|F\|_{\mathcal H^\beta}
 \right]
 \|f\|_{h_L^{p(\cdot)}(\rr^d)}
\]
for all \(f\in h_L^{p(\cdot)}(\rr^d)\).
\end{theorem}

The proofs of Theorems \ref{eq0-VHS}, \ref{eq1-VHS}, \ref{Riesz-VHS}, and \ref{hor-VHS} are similar to that of Theorem \ref{mol-VHS}; therefore, we omit the details.

\begin{remark}
To the best of our knowledge, Theorems \ref{mol-VHS}--\ref{eq1-VHS} coincide with the corresponding results in \cite{abdr20}, whereas Theorems \ref{Riesz-VHS} and \ref{hor-VHS} are entirely new.

\end{remark}

\subsection{Local mixed-norm Hardy spaces}\label{MHS}

In this subsection, as an example in the Euclidean setting, we apply
Theorems \ref{mol}, \ref{eq0}, and \ref{eq1} to local mixed-norm Hardy
spaces associated with \(L\). We begin by recalling the definition of
mixed-norm Lebesgue spaces.

For a given vector
\[
 \vec r:=(r_1,\ldots,r_d)\in(0,\infty)^d,
\]
the mixed-norm Lebesgue space \(L^{\vec r}(\mathbb R^d)\) is defined to be
the set of all measurable functions \(f\) on \(\mathbb R^d\) such that
\[
 \|f\|_{L^{\vec r}(\mathbb R^d)}
 :=
 \left\{
 \int_{\mathbb R}
 \cdots
 \left[
 \int_{\mathbb R}
 |f(x_1,\ldots,x_d)|^{r_1}\,dx_1
 \right]^{r_2/r_1}
 \cdots dx_d
 \right\}^{1/r_d}
 <\infty,
\]
with the usual modifications when \(r_i=\infty\) for some
\(i\in\{1,\ldots,d\}\). Throughout this subsection, we write
\[
 r_-:=\min\{r_1,\ldots,r_d\}
 \quad\hbox{and}\quad
 r_+:=\max\{r_1,\ldots,r_d\}.
\]

It is worth pointing out that \(L^{\vec r}(\mathbb R^d)\), with
\(\vec r\in(0,\infty)^d\), is a ball quasi-Banach function space. However,
even when \(\vec r\in[1,\infty]^d\), \(L^{\vec r}(\mathbb R^d)\) may fail to
be a Banach function space; see, for instance, \cite[Remark 7.21]{ZYYW2021}.
The study of mixed-norm Lebesgue spaces goes back to H\"ormander
\cite{H1960} and Benedek and Panzone \cite{BP1961}. Further results on
mixed-norm Lebesgue spaces and related mixed-norm function spaces can be
found in \cite{CGN2019,CGN2017,CGN2017-2,HLYY2019} and the references
therein.

Taking
\[
 X(\mathbb R^d):=L^{\vec r}(\mathbb R^d),
\]
the local Hardy space \(h_L^X(\mathbb R^d)\) becomes the local mixed-norm
Hardy space associated with \(L\). In this case, we write
$$h_L^{\vec r}(\mathbb R^d):=h_L^X(\mathbb R^d), \quad
h^{\vec r}_{L,\mathrm{mol},M,\varepsilon}(\mathbb R^d)
:=h^X_{L,\mathrm{mol},M,\varepsilon}(\mathbb R^d),
\quad {\rm{and}} \quad
H^{\vec r}_{L+I}(\mathbb R^d):=H^X_{L+I}(\mathbb R^d).
$$
Applying Theorems \ref{mol}, \ref{eq0}, and \ref{eq1} with
\(X(\mathbb R^d)=L^{\vec r}(\mathbb R^d)\), we obtain the following
conclusions.

\begin{theorem}\label{mol-MNHS}
Let \(L\) be a non-negative self-adjoint operator on
\(L^2(\mathbb R^d)\) satisfying Assumption \ref{a2.5}. Let
$
 \vec r:=(r_1,\ldots,r_d)\in(0,\infty)^d
$ with
$
 0<r_{-}\leq 1,
 $ $
 r_{-}\leq r_{+}< 2.
$
Then the following assertions hold.

\begin{enumerate}
\item[{\rm(i)}]
If
$\varepsilon>d(\frac1{r_-}-\frac1{2})$ and
$2M>d(\frac2{r_-}-1),$
then
\[
 h^{\vec r}_{L,\mathrm{mol},M,\varepsilon}(\mathbb R^d)
 \subseteq
 h^{\vec r}_L(\mathbb R^d).
\]

\item[{\rm(ii)}]
If \(\varepsilon>0\) and \(M\in\nn\), then
\[
 h^{\vec r}_L(\mathbb R^d)
 \subseteq
 h^{\vec r}_{L,\mathrm{mol},M,\varepsilon}(\mathbb R^d).
\]
\end{enumerate}
\end{theorem}
\begin{proof}
Take $p:={r}_-$, $q_0:=2$, and $s_0\in(0,\,\min\{1,{r}_-\}]$ such that $2M \in(d [\frac{2}{s_0}-1],\fz)\cap\mathbb{N}$ and $\varepsilon\in(\frac{d}{s_0}-\frac{d}{q_0},\fz)$.
From \cite{shyy17}, it follows that the BQBF space $X(\rr^d)=L^{\vec{r}}(\rr^d)$ satisfies Assumption \ref{a11} for above exponents.
Moreover, by \cite{shyy17}, we conclude that
$$\lf[\lf(\lf[L^{\vec{r}}(\rr^d)\r]^{\frac{1}{s_0}}\r)^{'}\r]^{\frac{1}{(2/s_0)^{'}}}=
L^{\frac{(\vec{r}/s_0)^{'}}{(2/s_0)^{'}}}(\rr^d),$$
where, for any $x\in\rr^d$, $\frac{1}{(r_j/s_0)^{'}}+\frac{1}{r_j/s_0}=1$.
Therefore, from the assumption that $0<r_-\leq r_+\leq2$ and \cite{shyy17}, we deduce that the BQBF space $X(\rr^d)=L^{\vec{r}}(\rr^d)$ satisfies Assumption \ref{a12} for above exponents. Therefore, all the assumption of Theorem \ref{mol} are satisfied, this implies the desired result. Therefore, we finish the proof of Theorem \ref{mol-MNHS}.
\end{proof}
\begin{theorem}\label{eq0-MNHS}
Let \(L\) be a non-negative self-adjoint operator on
\(L^2(\mathbb R^d)\) satisfying Assumption \ref{a2.5}. Let
$
 \vec r:=(r_1,\ldots,r_d)\in(0,\infty)^d
$ with
$
 0<r_{-}\leq 1
 $,  $
 r_{-}\leq r_{+}< 2.
$
If
$
 \inf\sigma(L)>0,
$
then
\[
 h^{\vec r}_L(\mathbb R^d)=H^{\vec r}_L(\mathbb R^d)
\]
algebraically and topologically, where \(\sigma(L)\) denotes the spectrum of
\(L\) on \(L^2(\mathbb R^d)\).
\end{theorem}

\begin{theorem}\label{eq1-MNHS}
Let \(L\) be a non-negative self-adjoint operator on
\(L^2(\mathbb R^d)\) satisfying Assumption \ref{a2.5}. Let
$
 \vec r:=(r_1,\ldots,r_d)\in(0,\infty)^d$
with
$
 0<r_{-}\leq 1
 $,  $
 r_{-}\leq r_{+}< 2.
$
Then
\[
 h^{\vec r}_L(\mathbb R^d)=H^{\vec r}_{L+I}(\mathbb R^d)
\]
algebraically and topologically.
\end{theorem}

\begin{theorem}\label{Riesz-MHS}
Let \(\vec{r}\in(0,\fz)^d\). Assume that
$0<r_{-}\leq 1$ and $r_{-}\leq r_{+}< 2.$ Let
\(L:=-\operatorname{div}(A\nabla)\) be as above, and let \(T:=L+\kappa I\) with $\kappa>0$.
Then the local Riesz transform
$
       \mathcal R_{L,\kappa}^{\rm loc}
$
admits a bounded extension from \(h_L^{\vec{r}}(\rr^d)\) to
\(L^{\vec{r}}(\rr^d;\mathbb C^d)\). More precisely, there exists a positive constant
\(C\) such that, for all \(f\in h_L^{\vec{r}}(\rr^d)\cap L^2(\rr^d)\),
\[
 \left\|
 \mathcal R_{L,\kappa}^{\rm loc}f
 \right\|_{L^{\vec{r}}(\rr^d;\mathbb C^d)}
 \leq
 C\|f\|_{h_L^{\vec{r}}(\rr^d)}.
\]
\end{theorem}

\begin{theorem}
\label{hor-MHS}Let \(L\) be a non-negative self-adjoint operator on \(L^2(\rr^d)\) satisfying
Assumption \(\ref{a2.5}\), and
 \(\vec{r}\in(0,\fz)^d\). Assume that
$0<r_{-}\leq 1$ and $r_{-}\leq r_{+}< 2.$
Put \(T:=L+I\). Then there exists \(\beta_0>0\) such that, for every
\(\beta>\beta_0\) and every bounded Borel function \(F\) on
\([0,\infty)\) satisfying
\[
        \|F\|_{L^\infty([0,\infty))}+\|F\|_{\mathcal H^\beta}<\infty,
\]
the operator \(F(L+I)\), initially defined on \(L^2(\rr^d)\) by the spectral
theorem, extends uniquely to a bounded operator on \(h_L^{\vec{r}}(\rr^d)\). Moreover,
\[
 \|F(L+I)f\|_{h_L^{\vec{r}}(\rr^d)}
 \leq
 C
 \left[
 \|F\|_{L^\infty([0,\infty))}+\|F\|_{\mathcal H^\beta}
 \right]
 \|f\|_{h_L^{\vec{r}}(\rr^d)}
\]
for all \(f\in h_L^{\vec{r}}(\rr^d)\).
\end{theorem}

The proofs of Theorems  \ref{eq0-MNHS}, \ref{eq1-MNHS}, \ref{Riesz-MHS} and
\ref{hor-MHS} are similar to that of Theorem \ref{mol-MNHS}, we omit the details here.
\begin{remark}
To the best of our knowledge, Theorems \ref{mol-MNHS}--\ref{hor-MHS} are completely new.

\end{remark}
\begin{remark}
We should emphasize that,  our results in Sections \ref{s3} and \ref{s6x} can also be applied to another specific function spaces, for examples, Lorentz-Hardy spaces, Orlicz-slice Hardy spaces and so on, we refer the reader to \cite{shyy17,wyy20}. To limit the length of this paper, we leave these details to interested readers.

\end{remark}


\textbf{Acknowledgements.}
X. Liu was supported by the Joint Innovation Fund Project of Lanzhou Jiaotong
University and Beijing Jiaotong University (No. LH2026013), and the Foundation
for Innovative Fundamental Research Group Project of Gansu Province (No.
25JRRA805). W. Wang was supported by the China Postdoctoral Science
Foundation (No. 2024M754159) and the Postdoctoral Fellowship Program of the
CPSF (Nos. GZB20230961, GZC20261666)


%

\medskip

\noindent Xiong Liu \\
\noindent Address: School of Mathematics and Physics, Lanzhou Jiaotong University,
Lanzhou 730070, P. R. China \\
\noindent{E-mail}: liuxmath@126.com    \\

\noindent Wenhua Wang   \\
\noindent Address: Institute for Advanced Study in Mathematics, Harbin Institute of Technology, Harbin 150001, P. R. China \\
\noindent{E-mail}: whwangmath@whu.edu.cn    \\

\medskip

\noindent Tiantian Zhao \\
\noindent Address: Institute for Advanced Study in Mathematics, Harbin Institute of Technology, Harbin 150001, China \\
\noindent{E-mail}:zhaotiantian@hit.edu.cn     \\

%

\medskip

\noindent Hongliang Zuo \\
 Address: School of Mathematical and Statistics, Henan Normal University, Xinxiang 453007, P. R. China  \\
\noindent{E-mail}: zuohongliang@htu.edu.cn \\

\end{document}